\documentclass[11pt]{amsart}
\usepackage[a4paper,hmargin=28mm,vmargin=27mm]{geometry}
\usepackage[T1]{fontenc}
\usepackage{lmodern,amsmath,amssymb,amsthm,mathrsfs,enumitem}
\input{ot1lmr.fd}
\DeclareFontShape{OT1}{lmr}{m}{n}{<->rm-lmr10}{}
\usepackage{xcolor,needspace}
\usepackage{hyperref}
\hypersetup{colorlinks=true,linkcolor=blue!55!black,
  urlcolor=blue!55!black,citecolor=blue!55!black,
  pdftitle={Hodge semiampleness and b-semiampleness},
  pdfauthor={Zhengyu Hu, Kefeng Liu, Yang Shen}}
\newcommand{\C}{\mathbb C}
\newcommand{\Q}{\mathbb Q}
\newcommand{\Z}{\mathbb Z}
\newcommand{\R}{\mathbb R}
\newcommand{\cO}{\mathcal O}
\newcommand{\cI}{\mathcal I}
\newcommand{\an}{\mathrm{an}}
\newcommand{\red}{\mathrm{red}}
\DeclareMathOperator{\Hom}{Hom}
\DeclareMathOperator{\im}{im}
\DeclareMathOperator{\coker}{coker}
\DeclareMathOperator{\Spec}{Spec}
\DeclareMathOperator{\End}{End}
\DeclareMathOperator{\Gr}{Gr}
\DeclareMathOperator{\rk}{rk}
\theoremstyle{plain}
\newtheorem{theorem}{Theorem}[section]
\newtheorem{proposition}[theorem]{Proposition}
\newtheorem{lemma}[theorem]{Lemma}
\newtheorem{corollary}[theorem]{Corollary}
\theoremstyle{definition}
\newtheorem{definition}[theorem]{Definition}
\newtheorem{example}[theorem]{Example}
\theoremstyle{remark}
\newtheorem{remark}[theorem]{Remark}
\numberwithin{equation}{section}
\setlist[enumerate]{label=\textup{(\arabic*)},leftmargin=*}
\title{Semiampleness of Hodge line bundles and b-semiampleness of the moduli part}
\author{Zhengyu Hu}
\address{Mathematical Sciences Research Center, Chongqing University of Technology,
Chongqing 400054, China}
\email{zhengyuhu16@gmail.com}
\author{Kefeng Liu}
\address{Mathematical Sciences Research Center, Chongqing University of Technology,
Chongqing 400054, China;
Department of Mathematics, University of California at Los Angeles,
Los Angeles, CA 90095-1555, USA}
\email{liu@math.ucla.edu}
\author{Yang Shen}
\address{Mathematical Sciences Research Center, Chongqing University of Technology,
Chongqing 400054, China}
\email{syliuguang2007@163.com}
\date{November 9, 2026}

\begin{document}
\begin{abstract}
We give proofs of Hodge semiampleness and b-semiampleness
without using o-minimality. We establish a geometric
semiampleness criterion for line bundles on proper algebraic
spaces. Let $(X,D)$ be a log smooth pair with $X$ a proper
algebraic space. For an integral polarizable pure variation
of Calabi--Yau type on $X\setminus D$ with unipotent local
monodromy, let $L$ be its extended highest Hodge line.
If $L$ is numerically integrable and has torsion combinatorial
monodromy, then $L$ is semiample. We also prove semiampleness
of extended Griffiths lines. We construct normal projective
period compactifications on which a positive power of the
Griffiths line extends to an ample line bundle.

As an application, we prove that the moduli b-divisor
$\mathbf M$ of a projective lc-trivial fibration over a
normal algebraic variety is b-semiample, provided that
the generic boundary is effective.
Suppose further that the geometric generic pair is klt and
its boundary is big, with nonzero coefficients in a fixed
finite set $I\subset(0,1)\cap\Q$. Then $c\mathbf M$ is b-Cartier and b-free for
an integer $c>0$ depending only on the relative dimension
and $I$. We extend b-semiampleness to analogous projective
fibrations over normal complex analytic bases, assuming a
global relative polarization and an invertible log
pluricanonical sheaf pulled back from the base.
\end{abstract}
\maketitle
\tableofcontents

\section{Introduction}\label{sec:introduction}

We work over $\C$. Let $X$ be smooth projective, let $D$ be an SNC divisor, and
put $U=X\setminus D$. Consider a polarized variation of Hodge structures
on $U$ with unipotent local monodromy. Suppose its extended highest
Hodge bundle is a line $L$. By the curvature formula, $L$ is nef.
To prove semiampleness, we seek a morphism
\[
 q:X\longrightarrow Y,\qquad L^a\simeq q^*A,
 \qquad a>0,\quad A\text{ ample}.
\]
The map $q$ must contract every curve $C$ with $L\cdot C=0$.
The difficulty is to construct it near these curves, especially
in the boundary. Under the hypotheses of Theorem~\ref{intro:hodge},
we first show that, for every connected chain $F$ of such curves,
there are an analytic neighbourhood $W\supset F$ and an integer
$m>0$ such that
\[
 L^{\otimes m}|_W \simeq \mathcal O_W.
\]
We prove this neighbourhood triviality in
Proposition~\ref{bdy:class-frame}; compare the local construction
in \cite[Proposition 3.13]{BFMT}.
We then adapt the inductive proof of a geometric semiampleness
criterion to obtain a globally generated power of $L$. This gives
projective compactifications
of period images and proves b-semiampleness of the moduli part in
the canonical bundle formula. Bakker--Filipazzi--Mauri--Tsimerman
\cite{BFMT} establish these semiampleness and compactification results.
We give independent proofs of Hodge semiampleness and b-semiampleness
using classical Hodge theory and algebraic geometry, without using
o-minimality.

Griffiths related curvature to the differential of the period map
and proved projectivity of compact period images
\cite[Theorems 5.2 and 9.7]{Griffiths70}. Schmid's theory \cite{Schmid}
controls the extension of the Hodge filtration across the boundary.
Bakker--Brunebarbe--Tsimerman \cite{BBT} proved algebraicity of period
images and ampleness of their Griffiths line by o-minimal GAGA.
Green--Griffiths--Laza--Robles \cite{GGLR} studied completion of
period images and positivity of the augmented Hodge line.
Green--Griffiths--Robles \cite{GGR} constructed a topological
completion and proved neighbourhood triviality for powers of the
extended Griffiths line.

Our approach is influenced by the work of Liu and Shen
\cite[\S\S3--5]{LS}, who prove algebraicity of ordinary period
images using infinitesimal neighbourhoods and formal algebraization.
We use neighbourhood trivializations of extended Hodge lines
to lift sections through nonreduced thickenings and to descend
sections along finite diagrams. We prove the required statements in
Appendix~\ref{sec:section-descent}, without invoking the period-image
algebraicity theorem of \cite{LS}.

We follow the CY-minimal quotient construction in \cite{BFMT}
and the strategy used there for locally stable models and
pluricanonical representations. We prove the required neighbourhood
triviality, finiteness of relations and comparison of line bundles
here, without invoking the semiampleness theorems of that paper.

A log smooth pair consists of a smooth space and a reduced simple normal crossings
divisor. An integral polarizable pure variation is of
\emph{Calabi--Yau type} (CY type) if its highest nonzero Hodge filtration piece
is a line.

Numerical integrability requires bigness on smooth tests with
generically immersive induced period map, including boundary tests.
Torsion combinatorial monodromy requires a line that is torsion
on every normalized component of a nodal curve to be torsion on
the whole curve. Definition~\ref{prep:assumptions} gives the precise
conditions. The following is the Hodge theorem of
\cite[Theorem 1.6]{BFMT}.

\begin{theorem}[Hodge semiampleness; Theorem~\ref{thm:hodge-proper}]
\label{intro:hodge}
Let $(X,D)$ be a log smooth pair with $X$ a proper algebraic space.
Let $\mathbb V$ be an integral polarizable pure variation of Calabi--Yau type on
$X\setminus D$, with unipotent local monodromy. Let $L$ be the Schmid
extension of its highest Hodge line. If $L$ is numerically integrable
and has torsion combinatorial monodromy, then $L$ is semiample.
\end{theorem}

For smooth projective $X$, Theorem~\ref{thm:main} identifies the fibres
of the connected contraction with the equivalence classes generated
by curves of Hodge degree zero. Theorem~\ref{thm:normal-source} also gives
descent to a normal proper algebraic space $X$ carrying a line
$\mathcal L$. On a smooth projective log resolution
$\pi:\widetilde X\to X$, we require a Hodge line $\widetilde L$
satisfying the same hypotheses and an isomorphism
\[
 \pi^*\mathcal L\simeq\widetilde L^{\otimes r},\qquad r>0.
\]
No klt, lc or rational-singularity assumption on $X$ is needed.
Without descent of the line, Corollary~\ref{cor:hodge-b-line}
gives a b-semiample Hodge b-divisor.

The first application concerns period compactifications, for which
the Baily--Borel compactification \cite{BB66} is the classical model
in the locally symmetric case. We start with an algebraic variety
or algebraic space $Y$ carrying a period map. For an integral polarizable pure
variation, its Griffiths line is the tensor product of the
determinants of the nontrivial Hodge filtration steps. Its curvature
detects the differential of the full period map. Its Schmid
extension on a smooth projective log pair with unipotent local
monodromy is semiample by Corollary~\ref{cor:griffiths}. The resulting
contraction gives the following compactification, without a neatness
assumption.

\begin{theorem}[Normal period compactifications; Theorem~\ref{thm:normal-level-compactification}]
\label{intro:period}
Let $Y$ be a normal quasi-projective variety with a locally liftable
quasi-finite period map $Y^{\an}\to\Gamma\backslash D$ for polarized
integral pure Hodge structures, where $\Gamma$ is a discrete subgroup
of the integral isometry group of the polarization. Then $Y$ has a normal projective compactification
$\overline Y$ carrying an ample line whose restriction to $Y$ is a
positive descended power of the Griffiths line, with the algebraic
structure defined in Section~\ref{sec:applications}.
\end{theorem}

The full statement describes the section ring by growth at the
boundary and proves independence of the level and log model.
Corollary~\ref{cor:normal-functoriality} gives functoriality at neat
level; Section~\ref{sec:nonnormal-descent} treats conductor descent.

The second application is the canonical bundle formula. Beginning
with Kodaira's formula for elliptic surfaces \cite{Kod63}, it
expresses the canonical class of a fibration in terms of the base,
the singular fibres, and the variation of the fibres. For an
lc-trivial fibration $f:(V,\Delta)\to B$, the formula is
\[
 K_V+\Delta\sim_{\Q}f^*(K_B+B_B+M_B).
\]
Here $B_B$ is the discriminant divisor and $M_B$ is the moduli part.
The moduli divisors on birational models form a b-divisor
$\mathbf M$. The b-semiampleness conjecture
\cite[Conjecture 7.13.1]{PS09} asks for a birational model of $B$
on which a positive Cartier multiple of the moduli part is globally
generated.

Kawamata's positivity theorem \cite{Kaw98} and Ambro's work
\cite{Amb04,Amb05} established the fundamental positivity properties
of the moduli part. Fujino--Gongyo \cite{FG14} extended b-nefness
and abundance from the klt to the lc setting. Hu
\cite[Theorems 1.2 and 1.5]{Hu21} proved log abundance and
compatibility with restriction to strata on a good dlt model.
Earlier semiampleness cases include relative dimension one, with
effective adjunction in the work of Kawamata and Prokhorov--Shokurov
\cite{Kaw97,PS09}, and fibrations by abelian varieties or K3 surfaces
\cite[Theorem 1.2]{Fujino03}. We first state the algebraic case.

\begin{corollary}[B-semiampleness; Corollary~\ref{cor:cbf}]
\label{intro:cbf}
Let $f:(V,\Delta)\to B$ be a projective algebraic lc-trivial fibration,
where $B$ is a normal algebraic variety and $\Delta$ is a $\Q$-divisor
on $V$. Assume that $\Delta_\eta\ge0$, where $\eta$ is the generic
point of $B$. Then its moduli b-divisor is b-semiample.
\end{corollary}

The base need not be proper. After an alteration and possible
adjunction to a minimal lc stratum, we realize a multiple of the
moduli part as a power of a cyclic-cover Hodge line. We compare the
line bundles across the boundary. Moduli b-abundance and finite
pluricanonical representations then verify the hypotheses of
Theorem~\ref{intro:hodge}.

For normal projective varieties, Corollary~\ref{cor:lc-fano-image}
shows that contractions preserve lc Calabi--Yau type and lc Fano
type. Example~\ref{ex:weak-fano-image} shows that the analogous
assertion for lc weak Fano type fails even for a birational
contraction from a smooth threefold.

Let $\bar\eta$ be a geometric generic point of $B$. If the pair
$(V_{\bar\eta},\Delta_{\bar\eta})$ is klt and $\Delta_{\bar\eta}$ is
big on $V_{\bar\eta}$, with nonzero coefficients in a fixed finite
set $I\subset(0,1)\cap\Q$, Corollary~\ref{cor:effective-big-boundary}
gives an integer $c=c(n,I)>0$ such that $c\mathbf M$ is b-Cartier
and b-free, where $n$ is the relative dimension. This recovers
\cite[Corollary 7.8]{BFMT}. By Hacon--Xu boundedness \cite{HX15},
we need only finitely many parameter families. We prepare these
families once, so the indices and covering degrees used in norm
descent do not depend on the test fibration.

Theorem~\ref{thm:analytic-cbf} extends b-semiampleness to analogous
projective fibrations over normal complex analytic bases, under
additional polarization and pluricanonical pullback hypotheses.
We reduce to algebraic parameter families and compare their moduli
lines with that of the original fibration across the boundary.

We construct global sections by induction on dimension. The next
criterion isolates the geometric requirements: neighbourhood
trivializations for descent and transverse bigness for lifting sections.

For a closed algebraic equivalence relation $R$ on $X_{\red}$,
an integral closed subspace $Z\subset X$ is \emph{transverse to $R$}
if
\[
 \dim\bigl(Z\cap R(z)\bigr)=0\quad\text{for general }z\in Z.
\]
A neighbourhood is saturated if it is a union of whole classes.

\Needspace{10\baselineskip}
\begin{theorem}[Semiampleness criterion; Theorem~\ref{thm:frame-criterion}]
\label{intro:criterion}
Let $X$ be a proper complex algebraic space, and let $L$ be a line
bundle on $X$. Let $R\subset X_{\red}\times X_{\red}$ be a closed
algebraic equivalence relation. Assume the following:
\begin{enumerate}
\item there is an integer $a>0$ such that, for every $R$-class $F$,
the line $L^a$ is trivial on a saturated analytic open neighbourhood
of $F$ in $X^{\an}$;
\item $L|_Z$ is big for every integral closed subspace $Z\subset X$ transverse
to $R$.
\end{enumerate}
Then $L$ is semiample.
\end{theorem}

A frame of a line bundle on an open set is a section that generates
it as an invertible module. The frames in (1) are defined on the full
space, including its nilpotents. We use them to lift sections through nonreduced divisors
and to glue sections across the conductor. For Hodge lines, the
proof separates two cases. If the line is not big, complete period
fibres reduce the problem to a smaller-dimensional base. If it is
big, Kodaira vanishing lifts sections from an effective divisor.
We use minimal Hodge structures to prove that the gluing relation
is finite. Torsion combinatorial monodromy supplies the invariant
covector tensors from which we construct the neighbourhood frames.

Section~\ref{sec:preliminaries} recalls the Hodge conventions,
complete ordinary period fibres and finite descent.
Section~\ref{sec:preparation} proves the geometric criterion and
Hodge semiampleness, including descent to normal algebraic spaces.
Sections~\ref{sec:applications} and~\ref{sec:cbf} give the period
compactifications and canonical-bundle-formula applications.
Appendix~\ref{sec:section-descent} proves the section-lifting lemmas.

\medskip
\noindent\textbf{AI disclosure.}
AI tools assisted the authors with exploring arguments,
checking calculations and references, and revising the exposition.
The authors are responsible for the final manuscript.

\section{Preliminaries}\label{sec:preliminaries}

\subsection{Hodge lines and boundary variations}

Let $X$ be a smooth projective variety, let $D\subset X$ be a
reduced simple normal crossings divisor, and write $j:U=X\setminus D\hookrightarrow X$.
Let $\mathbb V$ be a polarized integral pure variation of Hodge structures
on $U$. Assume that its local monodromies are unipotent and that
$F^{m+1}=0$ and $\operatorname{rank}F^m=1$.
We write $(\mathcal V,F)$ for its Deligne--Schmid extension and set
$L=F^m\mathcal V$. For the ordinary period arguments, we pass to
neat monodromy by Lemma~\ref{intro:level}. For the limiting Hodge
structures and neighbourhood frames, we require only unipotent
local monodromy.

We write $X_I$ for the locally closed strata of $D$, including
$X_{\varnothing}=U$, and set $D_I=(\overline{X_I}\setminus X_I)_{\red}$.
On $X_I$, let $N_i$ be the normal monodromy logarithms and set
\[
 C_I=V_{\mathbb Q}/\sum_{i\in I}\im N_i,
 \qquad V_I=(C_I,W^I,F).
\]
Here $W^I$ is induced by the normal-cone monodromy weight filtration,
centred at the weight of $\mathbb V$, and $F$ is the image of the
extended filtration. Its lattice is the image of $V_{\mathbb Z}$ in
$C_I$. We record the limit compatibility used below.

\begin{lemma}[Boundary quotients and face limits]\label{prep:face-quotients}
The data $V_I$ define a canonical admissible graded-polarizable
variation on $X_I$. Fix normal coordinates at $x\in X_J$, where
$J=I\sqcup K$, and identify a nearby rational fibre $H$ in a logarithmic
frame. The quotient maps are
\[
 H\xrightarrow{\pi_I}C_I
 \xrightarrow{\pi_{IJ}}C_J,\qquad
 C_J=C_I\big/\sum_{k\in K}\im\overline N_k,
 \qquad\pi_J=\pi_{IJ}\pi_I.
\]
Here $\overline N_k$ is induced by $N_k$. We denote the face weight
filtrations on $H$ by $\widetilde W^I,\widetilde W^J$ and its extended
Hodge filtration at $x$ by $\widetilde F_x$. Thus $W^I=\pi_I\widetilde W^I$.
The limit of $V_I$ in the remaining normal directions has
\[
 M\left(\sum_{k\in K}b_k\overline N_k,W^I\right)
       =\pi_I\widetilde W^J,\qquad b_k>0,
\]
and Hodge filtration $\pi_I\widetilde F_x$ in these coordinates.
For morphisms of these admissible local variations, taking face
limits is exact and the induced maps are strict. In particular, it
commutes with taking the cokernel of
$\bigoplus_{i\in I}H(1)\xrightarrow{\sum N_i}H$.
For a curve $z_k=t^{a_k}u_k(t)$, the additional factor
$\exp(\sum c_kN_k)$, with $c_k$ determined by logarithms
of the units $u_k(0)$, acts trivially after projection to $C_J$.
\end{lemma}
\begin{proof}
The face degeneration before taking coinvariants is an admissible
graded-polarizable mixed variation by
\cite[Theorem~2.9, Proposition~2.10 and the following paragraph]{CK1989}.
For $i\in I$, the maps $N_i:H(1)\to H$ are morphisms on this
face-limit variation. Their common cokernel is $V_I$.
Successive face limits remain admissible by
\cite[Corollary~5.5.4]{Kashiwara1986}. The face identity and exactness
in \cite[Propositions~5.2.5--5.2.6]{Kashiwara1986} give
\[
 \begin{aligned}
 M\left(\sum_{k\in K}b_kN_k,\widetilde W^I\right)
   &=\widetilde W^J,\\
 M\left(\sum_{k\in K}b_k\overline N_k,\pi_I\widetilde W^I\right)
   &=\pi_I\widetilde W^J.
 \end{aligned}
\]
Strictness gives the Hodge filtration $\pi_I\widetilde F_x$.
Notice that $F^mH(1)=F^{m+1}H=0$.

Changing a defining equation for $D_i$ changes the local limit
filtration by an exponential of $N_i$. This exponential is the
identity on $C_I$, so the local constructions glue canonically;
their rational local system is
$(i_I^{-1}j_*\mathbb V_{\Q}^{\vee})^{\vee}$, where
$i_I:X_I\hookrightarrow X$. The same argument gives the last
assertion after projection to $C_J$. Parallel transport preserves
the image lattice. Finally, the Deligne extension and its Hodge
subbundles algebraize on projective $X$; their constant-rank residue
cokernels and induced filtrations are algebraic on $X_I$.
\end{proof}

In the category of admissible rational mixed variations, we define
\[
 M_I=V_I^{\mathrm{min}}=V_I/K_I,
 \qquad K_I=\max\{K\subset V_I:F^mK=0\}.
\]
The maximum exists: sums of subobjects with zero highest piece still
have zero highest piece, by strictness, and the ranks are bounded.
A rational mixed Hodge structure or admissible mixed variation with
one-dimensional highest Hodge piece is \emph{CY-minimal} if it has
no nonzero subobject with zero highest piece, taking subobjects in
the respective category. Pointwise minimality refers to this condition
on a fibre.
Every nonzero subobject of $M_I$ contains its highest line. In particular,
its lowest nonzero weight piece $P_I=W_{k_I}M_I$ is pure and contains
that line. For a pure CY variation $P$, write $P^{\mathrm{tr}}$ for the
unique summand containing the highest line and having no nonzero
subvariation with zero highest piece. Semisimplicity identifies
$P^{\mathrm{tr}}$ with $P^{\mathrm{min}}$. In particular, $P_I$
is this summand of the lowest graded piece: a complementary
summand would have zero highest piece.
Our lattice convention is
\begin{equation}\label{prep:lattice}
 (M_I)_{\mathbb Z}=\im((V_I)_{\mathbb Z}\to(M_I)_{\mathbb Q}),
 \qquad (P_I)_{\mathbb Z}=(M_I)_{\mathbb Z}\cap(P_I)_{\mathbb Q}.
\end{equation}
We give $M(x)=(V_{I,x})^{\mathrm{min}}$ the image lattice from
$V_{I,x}$; a realization in $\operatorname{Gr}^W V_I$ need not give
the same lattice. On a subvariety, we restrict first and then take
the minimal quotient, since minimality can change under restriction.

\begin{lemma}\label{prep:minimum}
Let $q:A\twoheadrightarrow B$ be a surjective morphism of rational
mixed Hodge structures, or of admissible graded-polarizable mixed
variations. Suppose that it induces an isomorphism
$F^mA\simeq F^mB$ on the highest Hodge lines. Then
$A^{\mathrm{min}}\simeq B^{\mathrm{min}}$, and this isomorphism identifies
the image lattices obtained from the same original lattice.
If a variation on a connected base has one pointwise-minimal fibre,
then it is CY-minimal in the category of variations.
For a CY-minimal mixed Hodge structure $M$, restriction to $F^mM_{\mathbb C}$
is injective on $\operatorname{End}_{\mathrm{MHS}}(M)$.
\end{lemma}
\begin{proof}
Strictness gives $F^m\ker q=0$. If $K_A$ and $K_B$ denote the maximal
subobjects with zero highest piece, then
\[
 K_A=q^{-1}(K_B),\qquad A/K_A\xrightarrow{\sim}B/K_B.
\]
Both quotients carry the image of the same original lattice.
A nonzero subvariation with zero highest piece has a nonzero fibre
with zero highest piece at every point, contradicting the asserted
pointwise minimality.
Finally, an endomorphism killing the highest line has image with
zero highest piece, hence has zero image by minimality.
\end{proof}

\begin{lemma}\label{prep:restricted-pure-part}
Let $T$ be a smooth connected variety and let $h:T\to X_I$ be a
morphism. The lowest pure part of $(h^*V_I)^{\mathrm{min}}$ is
naturally isomorphic to $(h^*P_I)^{\mathrm{tr}}$ as a rational variation
of Hodge structures. This isomorphism identifies their highest Hodge lines.
\end{lemma}
\begin{proof}
We set $Q=(h^*V_I)^{\mathrm{min}}$. By
Lemma~\ref{prep:minimum}, this is also $(h^*M_I)^{\mathrm{min}}$.
For $k=k_I$, strictness gives $W_{k-1}Q=0$ and
\[
 W_kQ=\im(h^*P_I\to Q)\ne0,
\]
since this image contains the highest line. The kernel of
$h^*P_I\to Q$ has zero highest piece. Conversely, every subvariation
of $h^*P_I$ with zero highest piece is a subvariation of $h^*M_I$
and is killed in $Q$. Hence $W_kQ=(h^*P_I)^{\mathrm{min}}$, which
semisimplicity identifies with $(h^*P_I)^{\mathrm{tr}}$.
\end{proof}

\begin{lemma}\label{prep:line}
The line $L$ is nef. Schmid extension commutes with unipotent
logarithmic pullback, tensor powers, and highest-line-preserving
quotients, preserving the comparisons on the interior.
Let $A$ be the highest-line extension of a pure CY variation on a
smooth projective pair $(T,E)$ with unipotent local monodromy.
Write
\[
 \theta_{\mathrm{top}}:T_{T\setminus E}
       \longrightarrow\Hom(F^m,F^{m-1}/F^m)
\]
for its highest-line Higgs map. If $\theta_{\mathrm{top}}$ is
generically injective, then $A$ is big.
On a smooth projective curve, $\deg A=0$ forces the highest-line
Higgs field to vanish.

For an integral polarizable pure variation with unipotent local
monodromy on $T\setminus E$, the extended Griffiths line is nef
and is big if its period map is generically immersive. On a smooth
projective curve, degree zero of this line forces the full Higgs
field to vanish.
\end{lemma}
\begin{proof}
Because the monodromies are unipotent, the residues are nilpotent.
We therefore need no integral-eigenvalue corrections under pullback
or tensor operations \cite{Deligne,Schmid,CKS}.
For the curvature assertions, we set $d=\dim T$. Griffiths' curvature formula
\cite[Theorem 5.2]{Griffiths70} identifies the null directions of
$c_1(A,h)$ with $\ker\theta_{\mathrm{top}}$.
The Schmid metric extends with semipositive curvature and zero Lelong
numbers, and $A$ is nef \cite[Theorem 1.1 and Corollary 1.2]{FF}.
For its closed positive curvature current $\Theta\in c_1(A)$,
the volume inequality \cite[Theorem 1.2(i)]{Boucksom} gives
\[
 \operatorname{vol}(A)\ge\int_T\Theta_{\mathrm{ac}}^d
   =\int_{T\setminus E}c_1(A,h)^d.
\]
If $\theta_{\mathrm{top}}$ is generically injective, the last
integral is positive, so $A$ is big. On a curve, $\int_T\Theta=\deg A$; a positive
current of degree zero vanishes, so $\theta_{\mathrm{top}}=0$.
For the Griffiths line, the curvature is the sum of the curvature
forms of the Hodge filtration determinants; its null directions are
the kernel of the full Higgs field. The preceding volume and
degree arguments therefore give bigness and the degree-zero assertion
for this line as well.
\end{proof}

\begin{lemma}[Schmid-line comparison]\label{prep:schmid-comparison}
Let $\overline T$ be a smooth projective variety and let
$b:\overline T\to X$ be a morphism. Let $T^\circ\subset\overline T$
be an open subset such that $\overline T\setminus T^\circ$ is an SNC
divisor and $b(T^\circ)\subset X_I$. Suppose that on $T^\circ$ there
is a diagram of admissible variations
\[
 b^*V_I\twoheadrightarrow M_T\hookleftarrow Q_T,
 \qquad Q_T\text{ pure},
\]
whose arrows induce isomorphisms on the highest piece $F^m$.
If the boundary monodromies are unipotent, the resulting comparison
of highest lines extends to
\[
 b^*L\simeq\overline{F^mQ_T}
 \quad\text{on }\overline T.
\]
\end{lemma}
\begin{proof}
At a general point of a boundary prime, we choose a transverse disc
with parameter $t$. If its centre maps to $X_J$, where $J=I\cup K$, then
\[
 z_k=t^{a_k}u_k(t),\qquad a_k>0,\quad u_k(0)\ne0,\qquad
 N_T=\sum_{k\in K}a_k\overline N_k.
\]
The limit filtration on $C_I$ is induced by the ambient Deligne
filtration. Strictness for the normal monodromies gives
\[
 F^m\sum_{i\in I}\im(N_i:V(1)\to V)
 =\sum_{i\in I}N_iF^{m+1}V=0.
\]
At the disc centre the ambient highest line still maps nontrivially
to the boundary coinvariant line. The coordinate units contribute
only the invertible holomorphic matrix
$\exp(\sum_k\log(u_k)\overline N_k)$; they produce neither a zero
nor a pole. For $K=\varnothing$ this is the usual unipotent extension.

We set $K_T=\ker(b^*V_I\to M_T)$ and $C_T=M_T/Q_T$. In unipotent
Deligne frames, a flat morphism is a constant matrix, so canonical
extension is exact. Exactness of admissible limits \cite{CKS,SZ}
and Hodge strictness \cite[Th\'eor\`eme 2.3.5]{DeligneHII} give
\[
 F^m(b^*V_I)_{\lim}\simeq F^m(M_T)_{\lim}
 \simeq F^m(Q_T)_{\lim}.
\]
Indeed, $F^mK_T=F^mC_T=0$ on the punctured disc, and these Hodge
pieces keep their ranks under extension. The line comparison has
order zero along every boundary prime; normality of $\overline T$
then extends it and its inverse across codimension two.
\end{proof}

\begin{definition}\label{prep:assumptions}
Let $(X,D)$, $\mathbb V$ and $L$ be as above.
\begin{enumerate}
\item We say that $L$ is \emph{numerically integrable} if there exists
a logarithmic modification $\mu:(X',D')\to(X,D)$ as in
Lemma~\ref{prep:log-pullback} with the following property.
Write $L'=\mu^*L$, and let $P'_I$ be the pure variation associated
above to the stratum $X'_I$ of $(X',D')$.
For every stratum $X'_I$, every log smooth pair $(T,E)$ with $T$
projective, and every morphism $g:T\to\overline{X'_I}$ generically
finite onto its image such that $g(T\setminus E)\subset X'_I$, set
\[
 P_T=\bigl((g|_{T\setminus E})^*P'_I\bigr)^{\mathrm{tr}}.
\]
If the associated period map $\Phi_{P_T}$ is generically immersive,
then $g^*L'$ is big.

\item We say that $L$ has \emph{torsion combinatorial monodromy} if
the following holds for every morphism $c:C\to X$ from a proper
nodal curve. Let $\nu:\coprod_i C_i\to C$ be the normalization,
and write $\nu_i=\nu|_{C_i}$. If there are integers $a_i>0$ with
\[
 (\nu_i^*c^*L)^{\otimes a_i}\simeq\mathcal O_{C_i}
 \qquad\text{for every }i,
\]
then there is an integer $a>0$ such that
\[
 (c^*L)^{\otimes a}\simeq\mathcal O_C.
\]
\end{enumerate}
\end{definition}

Numerical integrability is called \emph{integrability} in
\cite[Definition 2.17]{BFMT}. The tests in (1) include smooth models
of all boundary subvarieties, and the transcendental part is taken
after pullback. By Lemma~\ref{prep:log-pullback}, the condition
persists on every higher logarithmic model, so we may use a common
refinement.

For a nodal test in (2), suppose that all $\nu_i^*c^*L$ are torsion.
The canonical finite connections on these lines glue by the node
identifications of $c^*L$. We take a common power trivial on every
$C_i$. The remaining monodromies are products of the gluing scalars
around cycles of the dual graph. They are all roots of unity if and
only if a further power is trivial on $C$. Thus (2) requires the
glued connection to have finite monodromy.

If the test line is not big, numerical integrability forces the
ordinary period fibres to have positive dimension. We use
combinatorial torsion to construct frames near zero-degree curve
chains in Section~\ref{sec:boundary}.

\subsection{Complete ordinary period fibres}

By Lemma~\ref{intro:level}, we choose a finite level inside a sufficiently
small principal congruence subgroup. Image lattices and saturated
sublattices inherit this level. Their monodromy lies in neat arithmetic
groups, and the ordinary period maps take values in the corresponding
quotients.

\begin{lemma}\label{lem:proper-period-extension}
Let $B$ be a smooth quasi-projective variety carrying a polarized
integral pure variation of Hodge structures. Assume that its monodromy
lies in a neat arithmetic group. Choose a smooth projective
compactification $\overline B$ such that
$\overline B\setminus B=\bigcup_jD_j$ is an SNC divisor, and let
$T_j$ denote the monodromy around $D_j$. The variation extends to
\[
B^+=\overline B\setminus
      \bigcup_{T_j\text{ of infinite order}}D_j,
\]
and its extended period map is proper. Every fibre is a closed
projective algebraic subspace of $\overline B$.
\end{lemma}
\begin{proof}
Since the monodromy group is neat, a finite-order $T_j$ is the
identity. We obtain the extension and its properness from Griffiths'
extension theorem \cite[Theorem 9.5, Propositions 9.10--9.11]{Griffiths70}.
Each fibre is compact analytic in $B^+$, hence closed analytic in
$\overline B$. It is therefore algebraic by Chow's theorem.
\end{proof}

\subsection{Pullback and descent}

\begin{lemma}\label{sg:generically-finite-descent}
\label{lem:pullback-descent}\label{lem:birational-sections}
Let $p:(T,E)\to(X,D)$ be a morphism of smooth log pairs with
$p(T\setminus E)\subset X\setminus D$. Let $L$ be the highest Hodge
line of a variation on $X\setminus D$ with unipotent local monodromy.
Then the pullback comparison on the interiors extends to
\[
p^*\overline L\simeq\overline{p^*L}.
\]
If $\tau:T\to X$ is proper, surjective and generically finite between
normal separated algebraic spaces of finite type, then
\[
\tau^*\mathcal L\text{ semiample}\quad\Longrightarrow\quad
\mathcal L\text{ semiample}
\qquad(\mathcal L\in\operatorname{Pic}(X)).
\]
If $\tau$ is birational, then $\tau_*\cO_T=\cO_X$ and pullback
identifies the graded section rings; generation of any fixed power
is equivalent on the two spaces.
\end{lemma}
\begin{proof}
The first assertion is the unipotent Deligne--Schmid comparison
of Lemma~\ref{prep:line}; coordinate units contribute an invertible
holomorphic matrix, so the comparison on the interior is preserved.

For the second assertion, we work componentwise and retain the dominating
components. We consider the Stein factorization
\[
T\xrightarrow{\beta}X'\xrightarrow{h}X,\qquad
\beta_*\cO_T=\cO_{X'},\qquad h\text{ finite},\quad d=\deg h.
\]
By the projection formula, generation of $\tau^*\mathcal L^m$
implies generation of $h^*\mathcal L^m$. For each $x\in X$, we choose
a section $s$ nonzero at every point of the finite fibre $h^{-1}(x)$.
Its norm over the normal base \cite[Tags 0BD3 and 0BCZ]{Stacks} is a section
of $\mathcal L^{md}$ nonzero at $x$. The exponent is independent
of $x$, and we take a common multiple over all components. If $\tau$
is birational, then $h$ is an isomorphism because $X$ is normal.
The projection formula then identifies all sections and their products,
with no change of power.
\end{proof}

\begin{lemma}[Descent along a nef reduction]\label{lem:nef-frame-descent}
Let $T$ be a smooth projective variety and let $N$ be a nef line
bundle. Suppose that the nef reduction of $N$ has dimension
$k<\dim T$. Assume that there is an integer $a>0$ such that
$N^a$ is trivial on an analytic neighbourhood of every connected
reduced closed subspace $F\subset T$ with $N|_F\equiv0$.
Then there are projective morphisms
\[
 T\xleftarrow{\ b\ }P\xrightarrow{\ \pi\ }B,
 \qquad \dim B=k,
\]
where $b$ is birational, $B$ is smooth, $\pi$ is flat, and
\[
 \pi_*\cO_P=\cO_B,\qquad b^*N^a\simeq\pi^*M
 \quad\text{for some }M\in\operatorname{Pic}(B).
\]
If $b^*N$ is semiample on a smooth projective multisection of $\pi$,
then $N$ is semiample.
\end{lemma}

\begin{proof}
We take the graph of the nef reduction \cite[Theorem 2.1]{NefReduction}.
We flatten it over a projective modification of its base, resolve that
base, and take the flat main transform. This gives $P,B,b,\pi$
as above. The generic fibre is geometrically integral and
$b^*N$ is numerically trivial on it. Hence $P$ is integral, every
fibre has pure dimension $d=\dim T-k$, and
$\pi_*\cO_P=\cO_B$ by Stein factorization and normality of $B$.
The flat dimension formula applies here; see
\cite[Tags 0BCM and 00ON]{Stacks}.

We take a $\pi$-ample line bundle $H$. Flatness gives
\[
 b^*N\cdot H^{d-1}\cdot[P_v]=0\qquad(v\in B).
\]
Since $b^*N$ is nef, its mixed degree on every component of $P_v$
is nonnegative, and hence zero. Given a curve on such a component,
we include it in a complete intersection of sufficiently large
powers of $H$. Nefness forces its degree to be zero. Hence
$b^*N|_{P_v}\equiv0$ for every $v$.
The reduced image $F=b(P_v)_{\red}$ is connected and
$N|_F\equiv0$. We pull back the assumed frame near $F$. Since $\pi$ is proper,
it is defined on $\pi^{-1}U$ for some analytic neighbourhood $U$
of $v$. Using $\pi_*\cO_P=\cO_B$, we descend the transition
functions to $B^{\an}$. They define a line bundle which is algebraic
by proper GAGA. We denote it by $M$; then
\[
 M\in\operatorname{Pic}(B),\qquad \pi^*M\simeq b^*N^a.
\]

Suppose that $b^*N|_W$ is semiample on a smooth projective
multisection $W\to B$. Then $(\pi|_W)^*M$ is semiample, so $M$
is semiample by norm descent. Finally, $b_*\cO_P=\cO_T$ because
$P$ is integral and $T$ is normal. By the projection formula,
the sections generating a power of $b^*N$ come from $T$ and
generate the same power of $N$.
\end{proof}

\begin{lemma}[Connected-fibre descent on nodal curves]
\label{nh:connected-descent}
Let $q:W\to X$ be a projective surjection of projective varieties
with connected fibres, and let $L\in\operatorname{Pic}(X)$.
Suppose that for every morphism $d:D\to W$ from a connected
projective nodal curve, $d^*q^*L$ is torsion whenever its restriction
to each normalized component of $D$ is torsion. Then for every
morphism $c:C\to X$ from such a curve, $c^*L$ is torsion whenever
its restriction to each normalized component of $C$ is torsion.
\end{lemma}
\begin{proof}
Let $c:C\to X$ be a morphism from a connected projective nodal
curve such that $P=c^*L$ is torsion on each normalized component.
After replacing $P$ by a power, these restrictions are trivial.
Writing $C_i$ for the normalized components of $C$, we choose
normalized integral multisections of $W\times_X C_i\to C_i$.
At each node, we join chosen lifts of the two branches by complete
curve chains in the connected projective fibre of $q$. We glue
these curves, inserting constant $\mathbb P^1$ components to separate
repeated attachment points. The line is trivial on each normalized
component of the resulting curve $D$, and the maps satisfy
\[
 d:D\to W,\qquad \pi:D\to C,\qquad q\circ d=c\circ\pi.
\]
The dual graph of $D$ contains a subdivision of the dual graph of $C$.
Transport through each inserted chain is the identity on the fibre
$P_p$. Hence a closed path $\gamma$ in the latter graph has a lift
$\widetilde\gamma$ with
\[
 \operatorname{Hol}_{P}(\gamma)
   =\operatorname{Hol}_{\pi^*P}(\widetilde\gamma).
\]
The right side is a root of unity by the hypothesis. Since the graph
has finitely generated first homology, one power kills all these
holonomies and trivializes $P$.
\end{proof}

\begin{lemma}\label{prep:log-pullback}
Let $\mu:(X',D')\to(X,D)$ be a projective modification of smooth
log pairs. Assume that $D'$ is SNC, that
$(\mu^{-1}D)_{\red}\subset D'$, and that the local monodromy around
each additional component of $D'$ is trivial. Let $T$ be a smooth
open subset of a stratum of $(X',D')$ such that $\mu(T)\subset X_I$.
Write $V'_T$ for the boundary variation of this stratum restricted
to $T$. There is an epimorphism
\[
 q_T:V'_T\twoheadrightarrow h^*V_I,
 \qquad (V'_T)^{\mathrm{min}}
       \simeq(h^*V_I)^{\mathrm{min}},\qquad h=\mu|_T.
\]
The map $q_T$ induces an isomorphism on the highest Hodge lines.
The displayed isomorphism identifies the image lattices from
$V_{\mathbb Z}$, and the extended highest line on $X'$ is $\mu^*L$.
Numerical integrability and torsion combinatorial monodromy pass
to $(X',D')$.
\end{lemma}
\begin{proof}
We write the modification locally as
\[
 z_i=u_i\prod_jw_j^{a_{ij}},\qquad
 N'_j=\sum_i a_{ij}N_i.
\]
For every $i\in I$ some $a_{ij}$ is positive. An interior element
of the new cone therefore maps to an interior element of the old
cone. The cone-weight identity gives $W'=h^*W^I$ on the ambient
limit, and
\[
 \sum_j\im N'_j\subset\sum_{i\in I}\im N_i.
\]
Taking coinvariants gives $q_T$. Coordinate units act trivially
on $h^*V_I$, and strictness shows that $q_T$ preserves the highest
line. Lemma~\ref{prep:minimum} identifies its minimal-quotient kernels:
\[
 K_{V'_T}=q_T^{-1}K_{h^*V_I},\qquad
 K(V'_{T,x})=q_{T,x}^{-1}K(V_{I,h(x)}).
\]
Both minimal quotients therefore carry the image of the same original
lattice. If $I=\varnothing$, the new monodromies vanish. The line
comparison is unipotent Schmid functoriality.

Let $a:Z\to X'$ be a smooth test as in Definition~\ref{prep:assumptions},
and put $v=\mu\circ a$ and $A=\overline{v(Z)}$. On the dense open
of $Z$ mapping into a stratum, Lemma~\ref{prep:restricted-pure-part}
identifies the lowest pure parts after restriction and minimization.
Over a dense smooth open of $A$, dominant pullback preserves
minimality: the unique maximal subvariation with zero highest piece
descends along a generic finite Galois multisection.
The induced pure period map $\Phi_Z$ therefore factors through $A$.
If it is generically immersive, then
\[
 \dim Z=\operatorname{rank}d\Phi_Z\le\dim A\le\dim Z.
\]
Hence $v$ is generically finite onto its image, and the original
numerical condition makes $a^*\mu^*L=v^*L$ big.
For combinatorial torsion, we compose each nodal test with $\mu$.
\end{proof}

\begin{lemma}\label{intro:level}
With $(X,D)$, $\mathbb V$ and $L$ as above, there is a proper
surjective generically finite morphism $\pi:\widetilde X\to X$, with
$\widetilde X$ smooth projective and
$\widetilde D=(\widetilde X\setminus\pi^{-1}U)_{\red}$ SNC, such that
the pulled-back variation has neat monodromy and Schmid highest line
$\widetilde L\simeq\pi^*L$. Numerical integrability and torsion
combinatorial monodromy pass to $\widetilde L$; if $\widetilde L$
is semiample, then $L$ is semiample.
\end{lemma}
\begin{proof}
We take a finite Galois cover $U'\to U$ at a sufficiently small neat
congruence level. We then normalize $X$ in $\C(U')$ and resolve the
boundary to obtain
\[
 \pi:\widetilde X\longrightarrow X,\qquad
 \widetilde L=\pi^*L.
\]
Lemma~\ref{lem:pullback-descent} gives the line-bundle comparison.
Numerical integrability passes to this model by the test condition
in Definition~\ref{prep:assumptions}; combinatorial torsion follows
by composing each nodal test with $\pi$. Semiampleness descends
along $\pi$ by Lemma~\ref{sg:generically-finite-descent}.
\end{proof}

\section{Semiampleness of Hodge lines}\label{sec:preparation}

We first prove a geometric semiampleness criterion. For Hodge
lines, closedness of the zero-degree relation is not known at the
start of the proof. We therefore adapt the induction underlying
the criterion, using neighbourhood frames and proving finiteness
of the conductor relations within the Hodge argument.

\subsection{A semiampleness criterion}\label{sec:frame-criterion}

We give algebraic relations on closed points their reduced structures.
An integral closed subspace $Z\subset X$ is \emph{transverse} to
$R$ if $\dim(Z\cap R(z))=0$ for general $z\in Z$.

\begin{theorem}\label{thm:frame-criterion}
Let $X$ be a proper complex algebraic space, and let $L$ be a line
bundle on $X$. Let $R\subset X_{\red}\times X_{\red}$ be a closed
algebraic equivalence relation. Assume the following:
\begin{enumerate}
\item there is an integer $a>0$ such that, for every $R$-class $F$,
the line $L^a$ is trivial on a saturated analytic open neighbourhood
of $F$ in $X^{\an}$;
\item $L|_Z$ is big for every integral closed subspace $Z\subset X$ transverse
to $R$.
\end{enumerate}
Then $L$ is semiample.
\end{theorem}

The neighbourhoods in (1) carry the full structure sheaf of $X^{\an}$,
including its nilpotents.

\begin{proof}
\emph{Step 1. Induction on the dimension.}
For an integral curve $C\subset X$ contained in an $R$-class, we
have $L\cdot C=0$ by (1). Otherwise $C$ is transverse because $R$
is closed, and $L\cdot C>0$ by (2). Hence $L$ is nef and
\begin{equation}\label{eq:frame-zero-curves}
 L\cdot C=0\quad\Longleftrightarrow\quad
 C\text{ is contained in an }R\text{-class}.
\end{equation}
We replace $L$ by the power in (1) and keep $X$ fixed. To apply
induction to divisors and conductors with their full scheme
structures, we prove the following stronger statement: for every
proper algebraic space $T$ and morphism $a:T\to X$, the line
$N=a^*L$ is semiample. We argue by induction on $\dim T$, using
$R_T=(a\times a)^{-1}R$ on $T_{\red}$.
Pulling back the neighbourhood frames proves (1) for $(T,N,R_T)$.
For (2), an integral subvariety transverse to $R_T$ maps generically
finitely to a subvariety transverse to $R$, so the pullback line
is big. The zero-dimensional case is immediate.

\emph{Step 2. The normal case.}
Suppose that $T$ is normal and integral. Birational descent of
section rings and Chow's lemma \cite[Tag 088U]{Stacks} allow us to
replace $T$ by a smooth projective resolution. If $N$ is big, we choose
an integer $b>0$, a line bundle $H$ on $T$, and an effective Cartier
divisor $E$ such that
\[
 N^b\simeq H\otimes\cO_T(E).
\]
We take $H$ to be globally generated, with both $H$ and
$H\otimes K_T^{-1}$ ample. Induction applies to $E$ with its full
scheme structure. For $m\ge b$, Kodaira vanishing gives
\[
 H^1(T,N^m(-E))=H^1(T,N^{m-b}\otimes H)=0.
\]
For sufficiently divisible $m$, sections generating $N^m|_E$
therefore lift to $T$. The defining section of $E$, multiplied by
sections of $H$, generates $N^b$ away from $E$. Taking common
powers gives generation everywhere.

If $N$ is not big, the general $R_T$-class is positive-dimensional
by (2). To obtain zero-degree curves through general points of $T$,
we choose an irreducible component $Q\subset R_T$ dominating $T$
under $p_1$ with positive-dimensional general fibres. Since
\[
 \dim Q>\dim T\ge\dim p_2(Q),
\]
the general $p_2$-fibre is positive-dimensional. At a general point of $Q$, we choose a complete
curve in this fibre through that point. Its image under $p_1$ is a
curve through a general point of $T$, contained in one $R_T$-class;
indeed, $p_1$ embeds every $p_2$-fibre into $T$.

The nef reduction theorem \cite[Theorem 2.1]{NefReduction} thus
gives a base of dimension less than $\dim T$. Every connected
numerically trivial subspace of $T$ maps into one $R$-class, by
\eqref{eq:frame-zero-curves}; the frame in (1) pulls back to its
neighbourhood. We apply Lemma~\ref{lem:nef-frame-descent} and choose a smooth
projective multisection $W$ of its flat morphism. Since
$\dim W<\dim T$ and $W$ still maps to $X$, the pullback of $N$
to $W$ is semiample by induction. The lemma gives
semiampleness on $T$.

\emph{Step 3. The conductor.}
Suppose that $T$ is reduced. Let $\nu:\widetilde T\to T$ be its
normalization, let $C\subset T$ be the full conductor, and set
$\widetilde C=C\times_T\widetilde T$. The pullback of $N$ to
$\widetilde T$ is semiample by the normal case, and $N|_C$ is
semiample by dimension induction. Its pullback to $\widetilde C$
is therefore semiample as well.
After taking a common power, we obtain connected contractions
$\widetilde T\to Y_{\widetilde T}$, $C\to Y_C$, and
$\widetilde C\to Y_{\widetilde C}$. The schemes
$Y_{\widetilde T},Y_C,Y_{\widetilde C}$ are projective and carry
ample line bundles $A_{\widetilde T},A_C,A_{\widetilde C}$.
Pullback of sections gives maps
\[
 Y_{\widetilde C}\rightrightarrows
  B:=Y_{\widetilde T}\sqcup Y_C.
\]
The two maps pull back $A_{\widetilde T}$ and $A_C$, respectively,
to $A_{\widetilde C}$ through their identifications with powers of $N$.
Since $A_{\widetilde C}$ is ample, both maps are finite.

By curve chains and \eqref{eq:frame-zero-curves}, each connected
fibre of these contractions maps into one $R$-class.
The inverse image of an $R$-class is proper, and the pulled-back
line bundle is numerically trivial there by (1). Ampleness therefore
makes its images in $Y_{\widetilde T}$ and $Y_C$ finite.
The relation on $B$ that identifies points belonging to the same
original class is thus closed and finite.

Let $k$ be a uniform
bound for the cardinalities of its fibres. Any two points in a
generated class can be joined by at most $k-1$ arrows, taken in
either direction. The generated relation is therefore a finite
union of these composites and the diagonal, hence closed and finite.

A frame near the corresponding original $R$-class descends to
frames on $B$. Their pullbacks agree
on $Y_{\widetilde C}$, with its full scheme structure, including
its nilpotents.

Lemma~\ref{lem:finite-polarized-diagram} supplies generating sections
on $Y_{\widetilde T}$ and $Y_C$ agreeing on $Y_{\widetilde C}$.
Their pullbacks descend through the conductor exact sequence
\[
 0\longrightarrow N^m\longrightarrow
 \nu_*\nu^*N^m\oplus N^m|_C
 \longrightarrow \nu_*N^m|_{\widetilde C}\longrightarrow0.
\]
These sections generate on $T$, proving the reduced case.

\emph{Step 4. Nilpotent thickenings.}
For arbitrary $T$, the line $N|_{T_{\red}}$ is semiample by the reduced
case. Each fibre of its connected contraction maps into an $R$-class
by \eqref{eq:frame-zero-curves}. We pull back the frame in (1) to a
neighbourhood of the whole fibre, with its nilpotent structure, and apply
Lemma~\ref{lem:nilpotent-frame-lifting} to prove semiampleness on $T$.

\end{proof}

If the classes are connected, they are exactly the fibres of the
connected semiample contraction $q:X\to Y$: the frames make $L$
numerically trivial on each class, and curve chains together with
\eqref{eq:frame-zero-curves} give the converse.

\subsection{Minimal Hodge structures and finite relations}

We follow \cite[Lemma 2.13]{BFMT} in the next lemma. Under its
hypotheses, a CY-minimal mixed variation becomes constant after a
finite \'etale cover if its lowest pure part is isotrivial.
We then use pure period maps to control integral
isomorphism types and the finite relations needed for descent.

\begin{lemma}\label{prep:fixed-part}
Let $T$ be a connected smooth algebraic variety and let $N$ be an
admissible variation of mixed Hodge structures on $T$. Assume that
$N$ is graded-polarizable and CY-minimal. Write $P=W_kN$ for its
lowest nonzero weight piece. If $P$ is isotrivial, then $N$ becomes
constant after a finite étale cover of $T$.
\end{lemma}
\begin{proof}
We pass to a finite Galois étale cover on which $P$ is constant.
The maximal subvariation with zero highest piece is Galois-invariant
and descends, so $N$ remains minimal on this cover. Hodge symmetry
and $F^m\operatorname{Gr}_{>k}^WN=0$ give
\[
 \begin{gathered}
 (p,q)\text{ in }\operatorname{Gr}_{>k}^WN
 \ \Longrightarrow\ p,q\le m-1,\quad p+q\ge k+1,\\
 F^{m-k-1}W_{-k-1}N^\vee=0,
 \qquad F^{m-k}N^\vee\simeq F^{m-k}P^\vee.
 \end{gathered}
\]
A flat section $v$ of $F^{m-k}P^\vee$ has a unique lift
$\widetilde v\in F^{m-k}N^\vee$, and transversality gives
\[
 \nabla\widetilde v\in
 F^{m-k-1}W_{-k-1}N^\vee\otimes\Omega_T^1=0.
\]
Uniqueness makes the lift single-valued. The admissible mixed
fixed-part theorem \cite[Proposition 4.19]{SZ}, with the reduction
to curves in \cite[proof of Theorem 1]{Andre}, gives
\[
 K=H^0(T,N^\vee_{\Q})_T\subset N^\vee,\qquad
 F^{m-k}P^\vee\subset\im(K_{\C}\to P^\vee_{\C}).
\]
Here taking invariants commutes with $\Q\subset\C$. Since $P$ is
pure and minimal, the map $K\to P^\vee$ is surjective.
Dualizing gives a quotient $N\twoheadrightarrow K^\vee$ that
induces an isomorphism on the highest Hodge line. Minimality then
implies $K=N^\vee$. The evaluation isomorphism
preserves $W,F$, so $N$ is constant.
\end{proof}

We use the argument of Narasimhan--Nori \cite{NN} for polarizations
of fixed discriminant, as in \cite[Lemma 2.22]{BFMT}.

\begin{lemma}\label{prep:pure-isom}\label{lem:fixed-hodge-locus}
Let $P$ be a polarized integral variation of pure Hodge structures
on a smooth quasi-projective variety $T$. Let $H$ be an integral pure
Hodge structure with the same weight and rank as $P$. Then
\[
 Z_H(P)=\{t\in T:P_t\simeq H\text{ as integral Hodge structures}\}
\]
is closed algebraic.
\end{lemma}
\begin{proof}
We work on a connected component of $T$ and may assume that
$Z_H(P)\ne\varnothing$. We transport a fibre polarization $Q_0$
to $H$ and set
\[
 \Lambda=H_{\mathbb Z},\qquad
 B=\End_{\mathrm{Hdg}}(H_{\Q}),\qquad S=B^\theta,
\]
where $\theta$ is the $Q_0$-adjoint involution.
Since $H$ is polarizable, $B$ is semisimple.
Every integral polarization with the same discriminant as $Q_0$
has the form $Q_a(x,y)=Q_0(ax,y)$, where
\[
 a\in\Xi:=\{a\in S:Q_0(a\Lambda,\Lambda)\subset\mathbb Z\},
 \qquad \det_H(a)=1.
\]
The reductive group
\[
 G=\ker(\det_H:B^\times\longrightarrow\mathbb G_m)
\]
acts on $S$ on the right by $a\cdot g=\theta(g)ag$.
The subgroup $G_{\mathbb Z}:=G(\Q)\cap\operatorname{GL}(\Lambda)$
preserves the lattice $\Xi$, since this action is pullback of
bilinear forms.

To apply \cite[Theorem 6.9]{BHC}, we consider
\[
 F=\{a\in S_{\C}:\det_H(a)=1\}.
\]
To check that its $G_{\C}$-orbits are closed, we take $a\in F$.
Since $a$ is invertible, it has a polynomial square root
$b\in\C[a]$. Then $\theta(b)=b$,
$a=\theta(b)b$, and $\det_H(b)=\pm1$.
Thus $F$ is a single orbit of
$G^\pm:=\{g\in B_{\C}^\times:\det_H(g)=\pm1\}$.
The subgroup $G_{\C}$ is normal of index at most two in $G^\pm$,
so its orbits on $F$ are finite in number and have the same
dimension. They are closed: $F$ is closed in $S_{\C}$, and the
boundary of a nonclosed orbit would contain an orbit of smaller
dimension.

The cited theorem gives finitely many $G_{\mathbb Z}$-orbits
in $F\cap\Xi$. Hence the integral polarizations under consideration
have finitely many orbits under
$\operatorname{Aut}_{\mathrm{Hdg}}(\Lambda)$.
We choose representatives $Q_1,\ldots,Q_s$.

We now pass to a finite \'etale cover $T'\to T$ whose monodromy lies
in a neat congruence subgroup of the integral isometry group.
For each $(H,Q_j)$ whose lattice is isometric to the fibre lattice,
its markings give finitely many points of the arithmetic period
quotient, since the congruence subgroup has finite index.
By Lemma~\ref{lem:proper-period-extension}, the inverse image of
each point is closed algebraic in $T'$. Their union is
$Z_H(P|_{T'})$. Its image under the finite map $T'\to T$ is
$Z_H(P)$, proving the assertion.
\end{proof}

\begin{lemma}\label{prep:generic-minimal}
Let $T$ be a connected smooth algebraic variety and let $N$ be an
admissible rational variation of mixed Hodge structures on $T$.
Assume that $N$ is graded-polarizable and CY-minimal. Then $N_t$ is
minimal outside a countable union of proper closed analytic subsets
in flat charts.
Suppose in addition that $N$ carries an integral lattice and that,
for some integral mixed Hodge structure $H$,
$(N_t)^{\mathrm{min}}\simeq H$ for every $t\in T$, where each minimal
quotient carries the image lattice. Then $N$ is locally constant as
a mixed Hodge structure, has finite monodromy, and $N_t\simeq H$
for every $t\in T$.
\end{lemma}
\begin{proof}
In a flat chart, we let $0\ne K\subset N_{\Q}$ be a rational subspace.
We set
\[
 d=\dim K,\qquad w=\sum_k k\dim\Gr_k^WK,\qquad
 0\ne\delta_K\in\bigwedge^dK.
\]
Since $W$ is flat,
$\delta_K\in W_w\bigwedge^dN\setminus W_{w-1}\bigwedge^dN$.
We claim that $K$ is a sub-MHS exactly when
\[
 w\in2\Z\quad\text{and}\quad
 \delta_K\in F^{w/2}\bigwedge^dN_{\C}.
\]
The forward implication follows by taking determinants.
For the converse, $\delta_K$ defines a Hodge morphism
$\Q(-w/2)\to\bigwedge^dN$, and
\[
 K=\ker\left(N\xrightarrow{v\mapsto v\wedge\delta_K}
       \bigwedge^{d+1}N\otimes\Q(w/2)\right).
\]
Hence the locus where $K$ is a sub-MHS is closed analytic.

Suppose this locus contains an open set. By analytic continuation
on the universal flat cover, $K$ is then a sub-MHS everywhere.
Its Hodge numbers are locally constant. If $F^mK=0$ at one point,
it is zero everywhere, and
\[
 K^{\rm var}:=\sum_{\gamma\in\pi_1(T)}\gamma K
\]
is a nonzero subvariation with zero highest piece. The sum is
finite-dimensional, so it contradicts minimality of $N$.
Therefore every nonminimal fibre belongs to a proper incidence
locus for some rational $K$. There are only countably many such
$K$, and we can use a countable flat atlas. This establishes the first
assertion.

For the second assertion, we take a minimal fibre.
Its isomorphism with $H$ gives $\operatorname{rk}N=\dim H$. Hence every quotient
$N_t\to(N_t)^{\mathrm{min}}\simeq H$ is an isomorphism, also on
the image lattices. In a flat chart, the filtration takes values
in the countable integral isomorphism orbit of $H$. Being
holomorphic, it must be constant.

It remains to check monodromy. On the lowest pure part, monodromy
preserves the Hodge metric and its lattice, so its image
is finite. Its character on the highest line therefore has
finite image. By Lemma~\ref{prep:minimum}, this character is
faithful on the monodromy of $N$, which is consequently finite.
\end{proof}

\begin{lemma}[Constancy along zero-degree curves]\label{prep:zero-curve}
Let $x,y\in X$ be joined by a chain of complete curves, each of
$L$-degree zero. Then $M(x)\simeq M(y)$ as integral mixed Hodge
structures, with the image lattices of \eqref{prep:lattice}.
\end{lemma}
\begin{proof}
It suffices to compare points on one complete curve $C$ with
$L\cdot C=0$. Suppose its generic point lies in $X_I$.
On the normalization of $C$, let $M_C$ be the minimal quotient
of the restricted variation, with lowest pure part $P_C$.
By the Schmid comparison and curvature formula, the highest
line is flat. Its character is finite by
\cite[Corollaire 4.2.8(iii)(b)]{DeligneHII}.
After a finite cover, the fixed-part theorem gives a constant
rational subvariation containing that line. Minimality forces
this subvariation to be $P_C$. Lemma~\ref{prep:fixed-part}
then shows that $M_C$ is constant. Its fibres are minimal.

We must also compare with the minimal quotient at a missing
point. If this point maps to $X_{I\cup K}$, we put
\[
 N_C=\sum_{k\in K}a_kN_k,\qquad a_k>0,\qquad A=C_I/\im N_C.
\]
Give $C_I$ its admissible curve-limit mixed Hodge structure and
$A$ the induced quotient structure. The limiting quotient to the
constant variation $M_C$ kills $N_C$. Lemma~\ref{prep:face-quotients}
therefore gives strict highest-line-preserving maps
\[
 A\twoheadrightarrow(M_C)_0,\qquad
 A\twoheadrightarrow C_{I\cup K}.
\]
For the second map, the coordinate-unit factors disappear on
$C_{I\cup K}$. The highest line survives because
\[
 F^m\im\!\left(\bigoplus N_k:\bigoplus V(1)\to V\right)=0.
\]
Taking minimal quotients yields
\[
 (M_C)_0\simeq A^{\mathrm{min}}
          \simeq(V_{I\cup K})^{\mathrm{min}}.
\]
These are integral identifications: both lattices are the image
of the original $V_{\mathbb Z}$ in $A^{\mathrm{min}}$, namely
\[
 ((M_C)_0)_{\mathbb Z}
 =\im(V_{\mathbb Z}\to A^{\mathrm{min}})
 =((V_{I\cup K})^{\mathrm{min}})_{\mathbb Z}.
\]
We can therefore compose the comparisons at successive points
of the curve chain.
\end{proof}

We first prove finiteness at general points and then on a saturated
closed complement.

\begin{lemma}[Finiteness of the generated equivalence relation]\label{lem:hodge-finite-diagram}
For $1\le i\le s$, let $a_i:T_i\to X$ be projective morphisms and
let $h_i:T_i\to B_i$ be projective surjections with connected fibres,
where each $B_i$ is a projective scheme. Suppose
that $A_i$ is an ample line bundle on $B_i$ and
\[
 h_i^*A_i\simeq a_i^*L^{m_i},\qquad m_i>0.
\]
Put $B=\coprod_iB_i$. Let $S\subset B_{\red}\times B_{\red}$ be
closed, contain the diagonal, and be invariant under exchanging the
two factors. Assume that both projections $S\to B_{\red}$ are finite.
Assume that, for every $(b_i,b_j)\in S$, the images
$a_i(h_i^{-1}(b_i))$ and $a_j(h_j^{-1}(b_j))$ can be joined by
a chain of complete curves of $L$-degree zero in $X$. Then the equivalence relation
generated by $S$ is closed and finite.
\end{lemma}
\begin{proof}
\emph{Step 1. Fixed Hodge types.}
A connected $h_i$-fibre is joined by curve chains on which $a_i^*L$
has degree zero. Hence Lemma~\ref{prep:zero-curve} defines
\[
 \tau(b)=[M(a_i(t))]\qquad(t\in h_i^{-1}(b)).
\]
This does not depend on $t$. By the assumption on $S$,
$\tau$ is constant along every $S$-chain.
The hypotheses and this description persist under base change
of $B$, so we may normalize, reduce, or restrict it to a closed
subspace when needed.

On each integral component $B_i$, we choose a smooth projective
multisection $u_i:W_i\to B_i$ with evaluation $v_i:W_i\to X$.
After resolving the boundary, we take a dense open $B_i^\circ$
over which $u_i$ is finite and $W_i^\circ$ maps to one stratum
$X_{I_i}$. On $W_i^\circ$, we set
\[
 N_i=(v_i^*V_{I_i})^{\mathrm{min}},\qquad
 P_i=\text{the lowest pure part of }N_i.
\]
By Lemma~\ref{prep:schmid-comparison}, the highest line extends
to $v_i^*L$, and $(v_i^*L)^{m_i}\simeq u_i^*A_i$.

For a fixed integral pure CY-minimal Hodge structure $H$ of
the appropriate weight and rank, the locus $Z_H(P_i)$ is closed
algebraic by Lemma~\ref{prep:pure-isom}. In fact, it is finite.
Otherwise we choose a curve in a positive-dimensional component.
On a connected smooth piece, the filtration takes values in a
countable integral orbit, hence is constant, and the polarized
integral monodromy is finite. On the smooth projective completion
of the curve, its extended highest line therefore has degree zero.
This contradicts
\[
 (v_i^*L)^{m_i}\simeq u_i^*A_i,
\]
since the curve meets the locus where $u_i$ is finite and $A_i$
is ample.

\emph{Step 2. Orbits at very general points.}
We define
\[
 S_0=S,\qquad S_{j+1}=S_j\circ S_j,\qquad
 \langle S\rangle=\bigcup_{j\ge0}S_j.
\]
Each $S_j$ has finite projections. On every maximal-dimensional
component of $B$, we remove the projections of the nondominant
components of all the $S_j$. We also remove the complements of
$B_i^\circ$ and their inverse images under dominant components
of these iterates, together with the finite images and inverse
images of the nonminimal loci in
Lemma~\ref{prep:generic-minimal}.
Away from the branch loci, these are countably many nowhere
dense analytic incidence sets in flat charts. By Baire's theorem,
their complement is nonempty.

For a point $b$ in this complement, every orbit point lifts to
some $w\in W_i^\circ$ with $(N_i)_w$ minimal of integral type
$\tau(b)$. Its lowest pure part has type $P(\tau(b))$, with the
image lattice in \eqref{prep:lattice}. Hence
\[
 \langle S\rangle(b)\subset
 \bigcup_i u_i\bigl(Z_{P(\tau(b))}(P_i)\bigr).
\]
Step~1 makes the right side finite. Hence all orbits at these
very general points are finite.

\emph{Step 3. The closed relation.}
We argue by induction on $d=\dim B$, the case $d=0$ being clear.
We normalize $B$, pull back $S$, and include all pairs lying over
the same point of the original $B$. The sets
\[
 U_j=\{b:S_{j+1}(b)=S_j(b)\}
\]
are constructible and increasing. Since $S$ contains the
diagonal, equality here means that the orbit has stabilized.
The points found in Step~2 belong to $\bigcup_jU_j$.
If no $U_j$ were dense in a maximal-dimensional component, their
closures, together with the excluded incidence sets, would
contradict Baire's theorem. Some common $U_j$ therefore contains
a dense open of every maximal-dimensional component.

Let $R^{\mathrm{dom}}$ be the $d$-dimensional part of $S_j$,
together with the diagonals of the smaller-dimensional components.
This is a finite equivalence relation. Indeed, its dominant
projections to the normal $d$-dimensional components are
universally open, so the relevant fibre products have pure
dimension $d$. Transitivity on the dense open where the iterates
stabilize then extends by closure; see the proof of
\cite[Lemma 28]{KolFinite}.

Let $Z_0$ consist of the smaller-dimensional components of $B$
and both projections of all components of $S$ of dimension less
than $d$. We set
\[
 Z=R^{\mathrm{dom}}(Z_0).
\]
This is closed and has dimension less than $d$. It is also
$S$-saturated: dominant arrows lie in $R^{\mathrm{dom}}$, and
both ends of every other arrow lie in $Z_0$.
We apply induction on $Z$ to the families
\[
 T_i\times_{B_i}(Z\cap B_i)\longrightarrow Z\cap B_i
\]
and their maps to $X$. We obtain a closed finite generated
relation $R_Z$. Since $R^{\mathrm{dom}}\subset S_j$, its arrows
over $Z$ are represented by finite $S$-chains, which remain in
$Z$ by saturation. Hence $R^{\mathrm{dom}}|_Z\subset R_Z$, and
\[
 \langle S\rangle=R^{\mathrm{dom}}\cup R_Z.
\]
This establishes closedness and finiteness on the normalization.
Taking its image under the finite normalization map proves
the assertion on the original $B$.
\end{proof}

\subsection{Frames near numerically trivial subspaces}\label{sec:boundary}

A neighbourhood frame requires more than an abstract identification
of the pointwise minimal Hodge structures: their duals must agree
across strata as subspaces of the ambient sheaf. We establish this
compatibility first, then construct an invariant covector tensor
and extend it to a neighbourhood.

\begin{lemma}\label{bdy:class-system}
Let $F\subset X$ be a connected reduced projective subspace that is
a union of SNC strata, and write $i:F\hookrightarrow X$. Assume that
$M(x)\simeq M(y)$ as integral mixed Hodge structures for all $x,y\in F$.
The duals of the quotient maps to $M(x)$ identify the spaces
$M(x)^\vee$ with subspaces of $(j_*\mathbb V_{\Q}^\vee)_x$.
These subspaces form a local system
\[
 \mathbb K_F\subset i^{-1}j_*\mathbb V_{\Q}^{\vee}.
\]
Its dual $\mathbb M_F$ carries a locally constant mixed Hodge
structure, with highest line $L|_F$. On each connected stratum
$T\subset F$, it is the minimal variation $(V_I|_T)^{\mathrm{min}}$
and has finite monodromy.
\end{lemma}
\begin{proof}
On a connected stratum $T\subset F\cap X_I$, we put
$N_T=(V_I|_T)^{\mathrm{min}}$. Lemma~\ref{prep:generic-minimal}
identifies the quotient maps $(N_T)_y\to M(y)$ as integral
isomorphisms for all $y\in T$ and shows that $N_T$ is locally
constant with finite monodromy.

In an SNC polydisc, let $T_J\subset\overline{T_I}$ be incident
strata in $F$, with $J=I\cup K$. Every new normal monodromy on
$N_{T_I}$ is finite and unipotent, hence is the identity. Its
Deligne extension has zero normal residues, and its flat filtration
has the same pointwise-minimal limit as its fibres on $T_I$.
The limiting quotient kills each $N_k$, $k\in K$, and factors as
\[
 q_x:C_J=C_I/\sum_{k\in K}\im\overline N_k
       \twoheadrightarrow (N_{T_I})_x,\qquad x\in T_J.
\]
By Lemma~\ref{prep:face-quotients}, the limiting quotient is a
morphism of mixed Hodge structures, with source weight filtration
induced by $\widetilde W^J$. Its factorization through the strict
quotient $C_I\twoheadrightarrow C_J$ therefore makes $q_x$ a
morphism of mixed Hodge structures.
Since $F^m\sum\im N_k=0$, the map $q_x$ preserves the highest
line; coordinate-unit factors act trivially on $(N_{T_I})_x$.
Strictness gives $F^m\ker q_x=0$. Conversely, any sub-MHS of $C_J$
with zero highest piece has zero image in the pointwise-minimal
quotient $(N_{T_I})_x$. Hence $\ker q_x=K(C_J)$, identifying this
quotient with $M(x)$ and preserving the image lattice from $V_{\mathbb Z}$.

In an SNC coordinate neighbourhood of $x$, every incident stratum specializes to
the deepest one. Since the quotient maps come from the same ambient
coinvariant projection, the equality of their kernels identifies
the dual subspaces under generization with
\[
 M(x)^\vee\subset
       \bigcap_{j\in J}\ker N_j^\vee=(j_*\mathbb V_{\Q}^\vee)_x.
\]
Flat invariant sections on the polydisc therefore give a constant
subsheaf on its intersection with $F$. On overlaps, both subsheaves have stalk
$M(y)^\vee$ with the same embedding in $(j_*\mathbb V_{\Q}^\vee)_y$;
they consequently agree and glue to $\mathbb K_F$. Strictness of
limits identifies their weight and Hodge filtrations and their
highest line with $L|_F$.
\end{proof}

\begin{lemma}\label{bdy:curve-detection}
Suppose that $F$ is either a connected smooth projective variety
or a connected union of components of an SNC divisor in a smooth
projective variety. Then there are a connected proper nodal curve
$C$ and an algebraic map $g:C\to F$ such that
\[
 g_*:\pi_1(C^{\an})\twoheadrightarrow\pi_1(F^{\an}).
\]
\end{lemma}
\begin{proof}
If $F$ is smooth, we take a smooth complete-intersection curve and
apply the weak Lefschetz theorem. In dimension one we take $F$
itself, and for a point we use a constant $\mathbb P^1$.

Now write $F=\bigcup_i F_i$ as an SNC union. We take one point in
each connected component of every nonempty intersection $F_i\cap F_j$.
Bertini and weak Lefschetz give a connected smooth complete-intersection
curve $C_i\subset F_i$ through all the chosen points on $F_i$, with
$\pi_1(C_i)\twoheadrightarrow\pi_1(F_i)$; if $\dim F_i=1$, we take
$C_i=F_i$. The groupoid form of van Kampen for an SNC union shows
that $\pi_1(F)$ is generated by the component groups and paths
through the chosen overlap points. Any such path within $F_i$ can
be replaced, by a homotopy that fixes its two ends, by a path
in $C_i$: we adjust a chosen path using the surjection on
fundamental groups.

For each distinct chosen point $y$, we attach a $\mathbb P^1$ mapped
constantly to $y$ to the marked points of all $C_i$ through $y$,
using distinct attaching points on $\mathbb P^1$. This gives a
connected proper nodal curve mapping to $F$. The generators just
described lie in the image of its fundamental group.
\end{proof}

\begin{lemma}[A highest-line character and a neighbourhood frame]
\label{bdy:character-frame}
Let $F\subset X^{\an}$ be a compact connected reduced analytic
subspace, with inclusion $i:F\hookrightarrow X^{\an}$. Let
$\mathbb K\subset i^{-1}j_*\mathbb V_{\Q}^\vee$ be a local system,
and set $\mathbb M=\mathbb K^\vee$. Assume that $\mathbb M$ carries
a locally constant mixed Hodge structure and that the quotient dual
to the given inclusion identifies $L|_F$ with $F^m\mathbb M_{\mathbb C}$,
where $F^{m+1}\mathbb M_{\mathbb C}=0$ and
$\operatorname{rank}F^m\mathbb M_{\mathbb C}=1$.
If the monodromy character on this line has finite image, then
$L^n$ is holomorphically trivial on a neighbourhood of $F$ for some
integer $n>0$.
\end{lemma}
\begin{proof}
We first construct on $F$ a flat covector tensor that pairs
nontrivially with a power of $L$. On a fibre $M$, Deligne's
functorial bigrading \cite[Lemme 1.2.11]{DeligneHII} gives a
monodromy-equivariant projection
\[
 e:M_{\C}=\bigoplus_{p,q}I^{p,q}\longrightarrow F^mM_{\C},
 \qquad\ker e=\bigoplus_{p<m}I^{p,q}.
\]
We take $0\ne\lambda\in(F^mM_{\C})^\vee$ and an integer $n$
that kills the highest-line character. Then
\[
 \eta=(\lambda e)^{\otimes n}\in\operatorname{Sym}^nM_{\C}^\vee
\]
is monodromy invariant and restricts to the nonzero tensor
$\lambda^{\otimes n}$ on $(F^mM_{\C})^{\otimes n}$.
Since $\eta$ is invariant, it defines a section of
$\operatorname{Sym}^n\mathbb K_{\C}$ and hence of $i^{-1}\mathscr F$,
where $\mathscr F=j_*\operatorname{Sym}^n\mathbb V_{\C}^\vee$.

Since $F$ is compact, sheaf continuity gives the isomorphism
\[
 \varinjlim_{F\subset W}H^0(W,\mathscr F)
       \simeq H^0(F,i^{-1}\mathscr F).
\]
Indeed, we choose finitely many local representatives covering $F$
and shrink their neighbourhoods until the representatives agree
on overlaps. We thus extend the invariant tensor to one neighbourhood
$W$ of the whole $F$. In unipotent Deligne frames it has zero residues,
so it extends holomorphically across $D$. Pairing it with $L^n$,
we obtain $L^n|_W\to\cO_W$, nonzero at every point of $F$.
After shrinking $W$, this map is an isomorphism. We use the bigrading
only on the locally constant MHS on $F$ and extend only the resulting
flat covector tensor off $F$.
\end{proof}

\begin{proposition}\label{bdy:class-frame}
Let $X$ be a smooth projective variety and let $D\subset X$ be
an SNC divisor. Let $\mathbb V$ be an integral polarizable variation
of pure Hodge structures of CY type on $X\setminus D$, with unipotent
local monodromy. Let $L$ be the Schmid extension of its highest
Hodge line, and assume that $L$ has torsion combinatorial monodromy. Put
\[
 b_R=\operatorname{lcm}\{d\ge1:\varphi(d)\le R\},
 \qquad R=\operatorname{rank}\mathbb V,
\]
where $\varphi$ is Euler's totient function. If $F\subset X$ is
connected, reduced and projective, and $L|_F\equiv0$, then $L^{b_R}$
is trivial on an analytic neighbourhood of the whole $F$.
\end{proposition}
\begin{proof}
If $F\ne X$, we principalize its ideal together with $D$:
\[
 \mu:(X',D')\longrightarrow(X,D),\qquad
 F'=(\mu^{-1}F)_{\red}.
\]
The divisor $D'$ is SNC and $F'$ is a union of its components.
Additional boundary components inside $X\setminus D$ have trivial
monodromy. If $F=X$, we take $X'=X$ and use all its strata.
In either case the highest line is $L'=\mu^*L$ by
Lemma~\ref{prep:log-pullback}. The fibres of $\mu$ are connected,
so $F'$ is connected, and $L'|_{F'}\equiv0$ by the projection formula.
We join points of $F'$ by complete zero-degree curve chains.
By Lemma~\ref{prep:zero-curve}, their pointwise minimal MHS are
integrally isomorphic. Lemma~\ref{bdy:class-system} gives $\mathbb M_{F'}$ together with
the embedding of its dual in the ambient sheaf on the modified pair.

To prove that the highest-line character has finite image on all
of $F'$, we use Lemma~\ref{bdy:curve-detection} to choose a nodal
curve $g:C\to F'$ surjective on fundamental groups.
A dense open of each normalized
component maps into one stratum, where $\mathbb M_{F'}$ has finite
monodromy. The local system extends over the complete normalized
component, whose fundamental group is a quotient of that of the
open curve. Its highest line therefore has finite monodromy and
is the torsion line $g^*L'$ there. The induced connection is the
canonical finite connection on that torsion line: on a common
connected finite cover, the ratio of two flat trivializations is
a global holomorphic unit on a proper connected curve, hence constant.

The comparisons of Lemma~\ref{bdy:class-system} agree at the nodes
with the fibre identifications of $g^*L'$. Combinatorial torsion,
applied to the composed nodal test $C\to X$, gives finite holonomy
on $C$. Since $g_*$ is surjective, the highest-line character
$\chi:\pi_1(F')\to\C^*$ has finite image. Each value $\chi(\gamma)$
is an eigenvalue of a rational monodromy matrix on $\mathbb M_{F'}$.
If $\chi(\gamma)$ has order $d$, its cyclotomic polynomial divides
the characteristic polynomial of that matrix. Hence
\[
 \operatorname{ord}\chi(\gamma)=d
 \quad\Longrightarrow\quad
 \varphi(d)\le\operatorname{rank}\mathbb M_{F'}\le R.
\]
Since $\varphi(d)^2\ge d/2$, only finitely many orders occur.
Their least common multiple divides $b_R$, so $\chi^{b_R}=1$;
compare \cite[Corollary 4.14]{Flo14}.

Lemma~\ref{bdy:character-frame}, with $n=b_R$, gives a
trivialization of $(\mu^*L)^{b_R}$ on a neighbourhood $W'$ of $F'$.

Since $\mu$ is proper, we can choose a neighbourhood $W\supset F$
with $\mu^{-1}W\subset W'$. Normality of $X$ implies
$\mu_*\cO_{X'}=\cO_X$ analytically. By the projection formula,
we descend both the frame and its inverse and obtain
$L^{b_R}|_W\simeq\cO_W$.
\end{proof}

\subsection{The contraction and descent to normal spaces}
\label{sec:singular-sources}\label{sec:normal-sources}

\begin{theorem}[Projective Hodge contraction]\label{thm:main}
Let $(X,D)$ be a smooth projective pair with simple normal crossings
boundary, and let $\mathbb V$ be an integral polarizable pure variation
of Calabi--Yau type on $X\setminus D$, with unipotent local monodromy.
Let $L$ be the Schmid extension of its highest Hodge line.
If $L$ is numerically integrable and has torsion combinatorial monodromy,
there exist a normal projective variety $Y$, an integer $a>0$, and a
projective surjection $q:X\to Y$ such that
\[
 q_*\cO_X=\cO_Y,\qquad L^a\simeq q^*A,
 \qquad A\in\operatorname{Pic}(Y)\text{ ample}.
\]
Its fibres are exactly the equivalence classes generated by complete
integral curves $C\subset X$ with $L\cdot C=0$. In particular, $L$ is
semiample.
\end{theorem}
\begin{proof}
\emph{Step 1. Induction on the dimension.}
We pass to a log model satisfying Definition~\ref{prep:assumptions}
and to neat monodromy by Lemma~\ref{intro:level}.
Closedness of the zero-degree relation on $X$ is not yet available,
so we adapt the induction in the proof of
Theorem~\ref{thm:frame-criterion}. For this fixed $(X,D)$ and
$\mathbb V$, we prove that $N=a^*L$ is semiample for every projective
morphism $a:T\to X$, including nonnormal and nonreduced $T$,
by induction on $\dim T$. We establish finiteness of the conductor
relation separately in Step~3. The assertion in dimension zero
is immediate.

\emph{Step 2. The normal case.}
If $T$ is normal integral, section-ring descent allows us to
resolve $T$ and its boundary. On a dense open $T^\circ$ mapping
to $X_I$, put
\[
 a^\circ=a|_{T^\circ},\qquad
 P_T=((a^\circ)^*P_I)^{\mathrm{tr}}.
\]
By Lemmas~\ref{prep:restricted-pure-part} and
\ref{prep:schmid-comparison}, its highest line extends to $N$.
Write $\Phi_T$ for its period map. If $\Phi_T$ is generically
immersive, then $a$ is generically finite onto its image, since
$\Phi_T$ is constant along the smooth fibres of $a$. Numerical
integrability therefore makes $N$ big.

Suppose that $N$ is not big. Then
$\operatorname{rank}d\Phi_T<\dim T$, and
Lemma~\ref{lem:proper-period-extension} completes the period
fibres inside $T$. Through a general point there is thus a
positive-dimensional complete subvariety on which the pure
variation becomes constant after a finite cover. All its curves
have $N$-degree zero, so the nef reduction has dimension
less than $\dim T$.
For every connected reduced $F\subset T$ with $N|_F\equiv0$,
we have $L|_{a(F)_{\red}}\equiv0$. Proposition~\ref{bdy:class-frame}
trivializes $N^{b_R}$ near $F$. Lemma~\ref{lem:nef-frame-descent}
now reduces semiampleness to a smooth projective multisection,
whose dimension is smaller. Induction applies.

If $N$ is big, Step~2 of Theorem~\ref{thm:frame-criterion}
applies without change: induction gives generation on the
effective Cartier divisor with its multiplicities, and Kodaira
vanishing lifts the sections. This establishes the normal case.

\emph{Step 3. The conductor.}
We use the construction in Step~3 of
Theorem~\ref{thm:frame-criterion} for reduced $T$.
The normal case and dimension induction give semiampleness on
$\widetilde T$ and the full conductor $C$, respectively.
After a common power, the contractions give finite maps
\[
 p_1,p_2:Y_{\widetilde C}\rightrightarrows
 B=Y_{\widetilde T}\sqcup Y_C,
\]
because the ample lines on $B$ pull back to the ample line on
$Y_{\widetilde C}$. Lemma~\ref{lem:hodge-finite-diagram} proves
that the generated relation is closed and finite.
Let $A$ denote the ample line bundle on $B$. The prescribed
comparisons with powers of $N$ induce an isomorphism
\[
 \alpha:p_1^*A\xrightarrow{\sim}p_2^*A.
\]

We check the other requirement of that construction, namely
compatible frames near a whole finite class. Let $F\subset X$
be the union of the images of its contraction fibres.
These images are joined by the images of the corresponding
fibres over $Y_{\widetilde C}$, so $F$ is connected projective.
Lifting curves to the contraction fibres gives $L|_F\equiv0$.
We take a common power of the lines and their prescribed comparisons
and keep the same notation. Proposition~\ref{bdy:class-frame} gives
a frame $s$ of the corresponding power of $L$ near $F$. Properness gives an
analytic neighbourhood $U\subset B^{\an}$ of the finite class
whose inverse images in $\widetilde T$ and $C$ map into this
neighbourhood of $F$. Since both contractions have direct image
$\cO$, the pullback of $s$ and its inverse descend to a frame $e$
of $A^{\an}|_U$.

Set
\[
 W=p_1^{-1}(U)\cap p_2^{-1}(U)\subset Y_{\widetilde C}^{\an}.
\]
On the inverse image of $W$ under the contraction of $\widetilde C$,
both $\alpha(p_1^*e)$ and $p_2^*e$ pull back to the same frame induced
by $s$. This contraction also has direct image $\cO$, so
\[
 \alpha(p_1^*e)=p_2^*e\qquad\text{on }W.
\]
The equality holds on the full analytic structure of $W$, including
its nilpotents. The finite relation allows us to shrink $U$ to a
saturated neighbourhood of the class.
Lemma~\ref{lem:finite-polarized-diagram} and the conductor exact
sequence from Step~3 of Theorem~\ref{thm:frame-criterion} now give
generation on $T$.

\emph{Step 4. Nilpotent thickenings.}
For arbitrary $T$, each reduced fibre of the contraction of
$T_{\red}$ maps to a connected projective numerically trivial
subspace of $X$. The frame from Proposition~\ref{bdy:class-frame}
pulls back to a neighbourhood in $T$ with its nilpotents.
Lemma~\ref{lem:nilpotent-frame-lifting} completes the induction.
Taking $T=X$ and applying
Lemma~\ref{sg:generically-finite-descent} through the log
modification and finite-level cover proves semiampleness
on the original $X$.

\emph{Step 5. The fibres of the contraction.}
Let $q:X\to Y$ be the connected semiample contraction. Then
\[
 q_*\cO_X=\cO_Y,\qquad L^a\simeq q^*A,
 \qquad Y\text{ normal projective},\quad A\text{ ample}.
\]
A complete curve is contracted exactly when its $L$-degree
is zero. Write $R_L$ for the equivalence relation generated by
these curves. Since each connected projective fibre is joined by
curve chains, these fibres are exactly the $R_L$-classes.
\end{proof}

\begin{corollary}[Neighbourhood frames on $X$]
\label{cor:original-frame-bound}
Under the hypotheses of Theorem~\ref{thm:main}, every $R_L$-class
has a saturated analytic neighbourhood on which $L^{b_R}$ is trivial,
with $b_R$ as in Proposition~\ref{bdy:class-frame}. No neatness
assumption is needed.
\end{corollary}
\begin{proof}
The contraction in Theorem~\ref{thm:main} gives
$R_L=(X\times_YX)_{\red}$ on the original $X$.
Its reduced fibres are connected, projective and numerically trivial.
We apply Proposition~\ref{bdy:class-frame} using the original rank $R$.
If $W$ is the resulting neighbourhood of a fibre, we replace it by
$X\setminus R_L(X\setminus W)$. Since the projections of the relation
are proper, this is a saturated open neighbourhood of the same fibre.
\end{proof}

The bound $b_R$ gives neighbourhood frames; it does not assert
that $L^{b_R}$ is globally generated. We next descend semiampleness
to proper log pairs and to normal algebraic spaces carrying a line
bundle whose pullback is a power of the Schmid line.

\begin{theorem}[Hodge semiampleness]\label{thm:hodge-proper}\label{cor:proper}
Let $(X,D)$ be a log smooth pair with $X$ a proper algebraic space.
Let $\mathbb V$ be an integral polarizable pure variation of Calabi--Yau type on
$X\setminus D$, with unipotent local monodromy. Let $L$ be the Schmid
extension of its highest Hodge line. If $L$ is numerically integrable
and has torsion combinatorial monodromy, then $L$ is semiample.
\end{theorem}
\begin{proof}
By Chow's lemma and resolution, we choose a smooth projective log model
$\pi:(\widetilde X,\widetilde D)\to(X,D)$ with $\widetilde L=\pi^*L$.
The two Hodge hypotheses hold upstairs by the conditions on smooth and nodal tests
and Lemma~\ref{prep:log-pullback}. We apply Theorem~\ref{thm:main}
and descend semiampleness by Lemma~\ref{lem:birational-sections}.
\end{proof}

\begin{theorem}[Normal algebraic spaces]\label{thm:normal-source}
Let $X$ be an integral normal proper algebraic space, let
$U\subset X$ be a smooth dense open subset, and let $\mathbb V$ be
an integral polarizable variation of pure Hodge structures of CY
type on $U$. Choose a proper birational morphism
$\pi:\widetilde X\to X$ with $\widetilde X$ smooth projective such that
\[
\widetilde D=(\widetilde X\setminus\pi^{-1}U)_{\red}
\]
is an SNC divisor. Assume that the pulled-back variation has unipotent
local monodromy and that its Schmid highest line $\widetilde L$ is
numerically integrable with torsion combinatorial monodromy. Let
$\mathcal L\in\operatorname{Pic}(X)$ and fix an isomorphism, for some
integer $r>0$:
\begin{equation}\label{eq:normal-source-comparison}
\pi^*\mathcal L\simeq\widetilde L^r.
\end{equation}
Then $\mathcal L$ is semiample. For some integer $a>0$ and an ample
line bundle $A$, its connected contraction satisfies
\[
q:X\to Y,\qquad q_*\cO_X=\cO_Y,\qquad
\mathcal L^a\simeq q^*A,
\qquad Y\text{ normal projective},\quad A\text{ ample}.
\]
Its fibres are exactly the classes generated by complete curves of
$\mathcal L$-degree zero. If $X$ is projective, $q$ is projective.
The contraction is independent of the chosen resolution satisfying these
hypotheses. No neatness, rational-singularity, klt or lc condition on $X$
is required.
\end{theorem}
\begin{proof}
Theorem~\ref{thm:hodge-proper} makes $\widetilde L$ semiample. Since
$\pi_*\cO_{\widetilde X}=\cO_X$, the comparison gives
\[
 H^0(X,\mathcal L^m)=H^0(\widetilde X,\widetilde L^{rm}).
\]
Thus $\mathcal L$ is semiample, and its section ring defines the
asserted contraction independently of the model. A curve is
contracted exactly when its $\mathcal L$-degree is zero. Each fibre
lifts to a connected projective subspace of $\widetilde X$, so its
points are joined by images of curve chains. This proves the fibre
description. If $X$ is projective, so is $q$.
\end{proof}

For a rank-one reflexive sheaf $\mathcal H$, the same conclusion
applies to an invertible reflexive power $\mathcal H^{[r]}$ with a
fixed isomorphism $\pi^*\mathcal H^{[r]}\simeq\widetilde L^r$.
$\Q$-Cartierness alone does not supply this isomorphism.

\begin{remark}\label{rem:model-independence}
We can check the comparison on one log smooth model; see also
\cite[Remark 4.3]{BFMT}. On a common smooth projective log model,
unipotent functoriality gives
\[
 \widetilde X_1\xleftarrow{p_1}T\xrightarrow{p_2}\widetilde X_2,
 \qquad p_1^*\widetilde L_1\simeq p_2^*\widetilde L_2.
\]
Lemma~\ref{lem:birational-sections} then descends the comparison to
either model. Numerical integrability pulls back by
Lemma~\ref{prep:log-pullback}; it descends by lifting each test to a
generically finite multisection and descending bigness. Combinatorial
torsion pulls back directly. To descend it, we join the lifted branches
of a nodal test by curve chains in the projective exceptional fibres
and apply Lemma~\ref{nh:connected-descent}. The line is constant on
these fibres. Unipotence descends along the strict transforms of the
original boundary divisors.

For nonprojective $X$, the contraction is proper but need not be projective:
a constant variation on a smooth complete nonprojective variety has
the contraction to a point.
\end{remark}

\begin{corollary}\label{cor:normal-alteration}
Let $X$ be an integral normal proper algebraic space and let
$\mathcal L\in\operatorname{Pic}(X)$. Let $\tau:T\to X$ be a proper
surjective generically finite morphism, with $T$ smooth projective,
and let $E\subset T$ be a reduced SNC divisor. Suppose that $T\setminus E$
carries an integral polarizable variation of pure Hodge structures
of CY type with unipotent local monodromy. Write $L_T$ for its
Schmid highest line, and assume that $L_T$ is numerically integrable
and has torsion combinatorial monodromy. If
\[
\tau^*\mathcal L\simeq L_T^r\qquad(r>0),
\]
then $\mathcal L$ is semiample.
\end{corollary}
\begin{proof}
We apply Theorem~\ref{thm:hodge-proper} to $L_T$, then descend
the semiampleness of $\tau^*\mathcal L\simeq L_T^r$ by
Lemma~\ref{sg:generically-finite-descent}.
\end{proof}

\begin{corollary}[Hodge b-divisors]\label{cor:hodge-b-line}
Let $\mathbb V$ and $\pi:\widetilde X\to X$ satisfy the hypotheses
of Theorem~\ref{thm:normal-source} concerning the variation and its
log model; no line bundle on $X$ is assumed. Choose a nonzero rational
section of $\widetilde L$, and let $H_{\widetilde X}$ be its divisor,
so that $\cO(H_{\widetilde X})\simeq\widetilde L$.
The Cartier closure $\mathbf H$ of $H_{\widetilde X}$ is a
b-Cartier, b-semiample Hodge b-divisor,
independent of the log model up to the principal b-divisor associated
with the rational section. If $\mathbf H$ descends as a $\Q$-Cartier
b-divisor to a normal model $X_0$, its trace $H_{X_0}$ is semiample.
The trace being $\Q$-Cartier alone does not assert descent of $\mathbf H$.
\end{corollary}
\begin{proof}
By unipotent Schmid pullback, we have $L_T=p^*\widetilde L$ on every
higher log model $p:T\to\widetilde X$. Using the same rational section,
we obtain $H_T=p^*H_{\widetilde X}$. We compare the Cartier closures
on a common model and apply Theorem~\ref{thm:hodge-proper} to prove
b-semiampleness. When the b-divisor descends, its trace is semiample
by Lemma~\ref{lem:birational-sections}.
\end{proof}

\section{Period compactifications}\label{sec:applications}

Let $S$ be a reduced complex analytic space. A holomorphic map
$\Phi:S\to\Gamma\backslash D$
is \emph{locally liftable} if every point of $S$ has an analytic
neighbourhood on which $\Phi$ lifts holomorphically to $D$.
We call such a map horizontal if $\Phi\circ u$ is horizontal
for every holomorphic map $u:T\to S$ with $T$ smooth, including maps
with $u(T)\subset S_{\mathrm{sing}}$. Equivalently, the pullback of
$\Phi$ to a resolution of $S$ is horizontal: resolve a dominating
multisection of its fibre product with $T$ and use generic smoothness.
All period maps in this section are understood to be horizontal.

\subsection{Griffiths lines}

Let $X$ be smooth and projective, let $D\subset X$ be a reduced SNC divisor,
and let $\mathbb V$ be an integral polarizable pure variation on
$X\setminus D$ with unipotent local monodromy along $D$.
Let $a$ and $b$ be the least and greatest indices $p$ for which
$\operatorname{Gr}_F^p\mathcal V\ne0$, and put
\[
\Lambda(\mathbb V)=\bigotimes_{p=a+1}^{b}\det F^p\mathcal V,
\qquad
G(\mathbb V)=\bigotimes_{p=a+1}^{b}
\bigwedge^{\operatorname{rk}F^p}\mathbb V.
\]
The tensor variation $G(\mathbb V)$ is integral, polarizable, pure,
and of CY type, with highest line $\Lambda$. We verify the numerical
integrability and torsion combinatorial monodromy required by
Theorem~\ref{thm:main}.

\begin{lemma}\label{app:boundary-determinants}
On a boundary stratum $X_I$, let
\[
 H_I=\bigoplus_k\operatorname{Gr}_k^{W^I}\Psi_I\mathbb V.
\]
Here $\Psi_I\mathbb V$ is the limiting mixed variation in local
normal coordinates. Its weight-graded pieces are independent of
these coordinates and glue on $X_I$.
Let $\overline Z$ be smooth projective, let $Z\subset\overline Z$
be dense open with SNC complement, and let
$\overline u:\overline Z\to X$ satisfy $\overline u(Z)\subset X_I$.
Write $u:Z\to X_I$ for the restriction. If $u^*H_I$ has unipotent
monodromy along $\overline Z\setminus Z$, then
\[
 \overline u^*\Lambda_X\simeq\overline{\Lambda(u^*H_I)}.
\]
The determinant on the right is formed with the same filtration
indices $a+1,\ldots,b$ as $\Lambda(\mathbb V)$. The identity commutes
with further unipotent logarithmic pullback.
\end{lemma}
\begin{proof}
On $Z$, strictness for $W^I$ gives the comparison by taking
determinants of the sequences for $F^pW^I_k$ and multiplying over
$p=a+1,\ldots,b$.

We check the extension at a general point of each boundary prime.
Choose a transverse disc $\mathbb D$ and local normal coordinates
along $X_I$. Let $E$ be the pullback to $\mathbb D^*$ of the full
face-limit variation $\Psi_I\mathbb V$. It is admissible, as are
its weight subobjects and quotients
\cite[Theorem~2.9 and Proposition~2.10]{CK1989}.
Its monodromy is unipotent because it is unipotent on every
$\operatorname{Gr}^{W^I}_kE$. Canonical extension is therefore
exact on the underlying flat bundles. Exactness of admissible
limits and Hodge strictness give exact sequences of vector bundles
on $\mathbb D$:
\[
0\longrightarrow F^p\overline{W^I_{k-1}E}
\longrightarrow F^p\overline{W^I_kE}
\longrightarrow F^p\overline{\operatorname{Gr}^{W^I}_kE}
\longrightarrow0.
\]
Here we use \cite[Propositions~5.2.5--5.2.6 and Corollary~5.5.4]
{Kashiwara1986} and \cite[Th\'eor\`eme 2.3.5]{DeligneHII}.
Successive face limits identify $F^\bullet\overline E$ with the
pullback of the ambient extended Hodge filtration. In local
expressions $z_j=t^{a_j}u_j(t)$, the units contribute only the
invertible holomorphic matrix $\exp(\sum_j\log(u_j)N_j)$.

Taking determinants and multiplying over $p$ therefore extends
the prescribed comparison as an isomorphism on $\mathbb D$.
The comparison and its inverse have no pole along any boundary
prime and extend across codimension two by normality of
$\overline Z$. Canonical extension and these filtration maps
commute with further unipotent logarithmic pullback.
\end{proof}

\begin{proposition}\label{app:griffiths-integrability}
The highest line of $G(\mathbb V)$ is numerically integrable.
\end{proposition}
\begin{proof}
For the CY variation $G(\mathbb V)$, let $P_I^G$ be the lowest pure
boundary variation of Section~\ref{sec:preliminaries}. Fix a test
$g:(T,E)\to(\overline{X_I},D_I)$ as in
Definition~\ref{prep:assumptions}, and set
\[
 T^\circ=T\setminus E,\qquad g^\circ=g|_{T^\circ},\qquad
 P_T=((g^\circ)^*P_I^G)^{\mathrm{tr}},\qquad
 H_T=(g^\circ)^*H_I.
\]
For $r_p=\operatorname{rk}F^p\mathbb V$, put
\[
 G_{\mathrm{gr}}(H_T)
 =\bigotimes_{p=a+1}^{b}\bigwedge^{r_p}H_T.
\]
The monodromy weight filtration commutes with tensor and exterior
powers, so $P_T$ is a pure subquotient of this graded variation.
For a tangent vector $v$,
\[
 \theta_{H_T}(v)=0
 \ \Longrightarrow\ \theta_{G_{\mathrm{gr}}(H_T)}(v)=0
 \ \Longrightarrow\ \theta_{P_T}(v)=0.
\]
If $\Phi_{P_T}$ is generically immersive, $\theta_{H_T}$ is
generically injective. The curvature of $\Lambda(H_T)$ is the
sum of the Griffiths curvatures of the pure weight pieces; its
null space is $\ker\theta_{H_T}$. The volume inequality in
Lemma~\ref{prep:line} therefore gives
\[
 \overline{\Lambda(H_T)}\text{ big},\qquad
 g^*\Lambda_X\simeq\overline{\Lambda(H_T)}
\]
by Lemma~\ref{app:boundary-determinants}. This is the required
numerical integrability.
\end{proof}

The next argument follows the norm-one and polarization arguments of
\cite[Theorems 5.21--5.22]{GGR} and
\cite[Theorem 2.24 and Lemma 2.25]{BFMT}; we give the branch comparisons
at the nodes explicitly.

\begin{lemma}\label{app:griffiths-torsion}
The extended Griffiths line $\Lambda_X$ has torsion combinatorial monodromy.
\end{lemma}

\begin{proof}
\emph{Step 1. Comparison at the nodes.}
After a Tate shift, we may assume that the Hodge indices of $\mathbb V$
lie in $[0,n]$. Let $C$ be a connected proper nodal curve and
$g:C\to X$ a morphism such that $\deg g^*\Lambda_X=0$.

Since $\Lambda_X$ is nef, its pullback has degree zero on every
normalized component. We write $\widehat W^I$ for the monodromy weight
filtration on the full nearby-cycle space, before taking the
coinvariant quotient $C_I$. On the dense open of a component mapping
to $X_I$, we put
\[
 E_I=\bigoplus_k\Gr_k^{\widehat W^I}\Psi\mathbb V.
\]
By positivity and Lemma~\ref{app:boundary-determinants}, the full
Higgs field of $E_I$ vanishes. Its pure graded pieces therefore have
constant period map and finite monodromy.

At a missing point mapping to $X_J$, $J\supset I$, we choose
a boundary polydisc and write
$N_C=\sum_{j\in J\setminus I}a_jN_j$, $a_j>0$, on the full
nearby-cycle vector space.
The induced monodromy on $\Gr^{\widehat W^I}$ is finite and unipotent,
hence trivial. The relative monodromy filtration has its only jump
on $\Gr_k^{\widehat W^I}$ at $k$, so it equals $\widehat W^I$.
The cone identity \cite{CKS,SZ} therefore gives
\[
 \widehat W^I=M(N_C,\widehat W^I)=\widehat W^J.
\]
Admissible limit strictness identifies the pure Hodge filtrations
on the associated gradeds. Coordinate-unit factors act trivially
there, since each $N_j$ lowers $\widehat W^J$ by two.

In a common nearby-fibre lattice $V_{\mathbb Z}$, the equal weight
filtrations give the same free graded lattices
\[
 E_{k,\mathbb Z}=(\widehat W_k\cap V_{\mathbb Z})/
                     (\widehat W_{k-1}\cap V_{\mathbb Z});
\]
the intersections are saturated. These comparisons identify both branches at a node with the same
graded Hodge structure. Changes of logarithmic lift act by normal
monodromy, trivially on the common graded space, so the comparisons
agree on overlaps. They define an integral local system
$E=\bigoplus_k E_k$ on $C$, each $E_k$ carrying a locally constant
pure Hodge structure with polarizable fibres.
The determinant comparisons come from the exact sequences of these
filtrations on the Schmid fibre at each node. Hence
\[
 g^*\Lambda_X=\bigotimes_{p=1}^{n}\det F^pE
\]
identifies the node gluing of $g^*\Lambda_X$. On each complete normalized
component this flat connection is the canonical connection of the
torsion line.

\emph{Step 2. Absolute value of the character.}
For a monodromy transformation $\gamma$, we write
$a_{p,q}=\det(\gamma|_{E^{p,q}})$ and denote by $\chi(\gamma)$ its
action on the displayed line. Hodge symmetry and integrality give
\[
 \chi(\gamma)=\prod_{p,q}a_{p,q}^{p},\qquad
 |\chi(\gamma)|^2
 =\prod_{p,q}a_{p,q}^{p+q}
 =\prod_k\det(\gamma|_{E_k})^k=1.
\]
Indeed, the last determinants are $\pm1$ and the left side is positive.

\emph{Step 3. Finiteness of the character.}
The tensor $T=\bigotimes_{p=1}^{n}\bigwedge^{\rk F^pE}E$ is split
by weights, and its unique highest Hodge line lies in a pure summand
$T_w$. Let $P$ be the smallest rational Hodge quotient of $T_w$
in which this line remains nonzero. By pure semisimplicity, $P$ is simple. Monodromy
preserves this intrinsic quotient and the image of
$T_{\mathbb Z}\cap T_{w,\mathbb Q}$. We fix a rational polarization
$Q$ of $P$.
Let $\ell\subset P_{\C}$ be the highest Hodge line. Hodge
orthogonality pairs $\ell$ only with its one-dimensional conjugate,
and this pairing is unchanged because $|\chi(\gamma)|^2=1$.
The rational Hodge morphism therefore satisfies
\[
 Q-\gamma^*Q:P\longrightarrow P^\vee(-w),\qquad
 (Q-\gamma^*Q)|_\ell=0.
\]
Since $P$ is simple, this morphism is zero. Hence every $\gamma$
preserves $Q$.
The group of integral Hodge isometries of $P$ is finite, so the
character $\chi$ has finite image. A power of the componentwise flat
frames therefore glues over $C$, proving the required torsion.
\end{proof}

\begin{corollary}\label{cor:griffiths}
Let $(X,D)$ be a log smooth pair, where $X$ is a separated algebraic
space of finite type. Let $\mathbb V$ be an integral polarizable pure variation on
$X\setminus D$ with unipotent local monodromy.
Let $\Lambda_X$ be the Schmid extension of its Griffiths line.
Then $\Lambda_X$ is algebraic and semiample. If $X$ is projective,
the connected contraction $q:X\to Z$ of $\Lambda_X$ is projective,
and a complete integral curve $C\subset X$ is contracted by $q$
if and only if $\deg_C\Lambda_X=0$.
\end{corollary}

\begin{proof}
For $X$ smooth projective with SNC boundary, we apply
Theorem~\ref{thm:main} to $G(\mathbb V)$ using
Proposition~\ref{app:griffiths-integrability} and
Lemma~\ref{app:griffiths-torsion}. Tensor compatibility identifies
the extended highest line with $\Lambda_X$ and gives the stated
contraction.

For arbitrary $X$, we take a smooth quasi-projective \'etale chart
$T\to X$ and a smooth projective SNC compactification of $(T,D_T)$.
Finite unipotent reduction at infinity and Schmid's theorem give
coherent extended Hodge filtrations. Invariant direct images, with
Deligne's residue convention, give coherent filtered extensions on
the original compactification; proper GAGA algebraizes them.
Over $T$, their Griffiths determinant gives an algebraic line
$\Lambda_T$ whose analytification is the prescribed Schmid line,
since the monodromy along $D_T$ is unipotent.
For the two projections $p_1,p_2:T\times_XT\to T$, the
regular-singular Riemann--Hilbert correspondence algebraizes the
canonical comparisons of connections \cite{Deligne,Schmid}.
These preserve the Hodge filtrations and give
\[
 p_1^*\Lambda_T\simeq p_2^*\Lambda_T
 \quad\text{on }T\times_XT.
\]
The comparisons satisfy the cocycle condition and descend $\Lambda_T$
to $\Lambda_X$. The same argument algebraizes the unipotent
logarithmic pullback comparisons.

Chow's lemma and resolution give a smooth quasi-projective proper
birational log model of $X$. Compactifying it and applying finite-level
preparation and resolution gives a smooth projective log pair
$(\widetilde X,\widetilde D)$ with unipotent boundary monodromy and
a proper generically finite surjection $p:\widetilde X^\circ\to X$
from an open subset. If $\widetilde\Lambda$ is the Schmid Griffiths
line on $\widetilde X$, unipotent pullback gives
\[
 p^*\Lambda_X\simeq\widetilde\Lambda|_{\widetilde X^\circ}.
\]
We apply the projective case and descend semiampleness by
Lemma~\ref{sg:generically-finite-descent}.
\end{proof}

\begin{corollary}[Griffiths lines on normal algebraic spaces]
\label{cor:singular-griffiths}
Let $X$ be a normal proper algebraic space, let $U\subset X$
be a smooth dense open, and let $\mathbb V$ be an integral
polarizable pure variation on $U$. Choose a proper birational
$\pi:\widetilde X\to X$ with $\widetilde X$ smooth projective
and $\widetilde X\setminus\pi^{-1}U$ an SNC divisor. Assume
unipotent boundary monodromy upstairs, and let $\widetilde\Lambda$
be the Schmid-extended Griffiths line. If
\[
 \mathcal L\in\operatorname{Pic}(X),\qquad
 \pi^*\mathcal L\simeq\widetilde\Lambda^r\quad(r>0),
\]
then $\mathcal L$ is semiample, with the contraction described in
Theorem~\ref{thm:normal-source}.
In particular, an everywhere locally liftable period map on normal proper
$X$, with neat monodromy, has an algebraic semiample Griffiths line.
If the period map is quasi-finite, this line is ample and $X$ is projective.
\end{corollary}
\begin{proof}
Corollary~\ref{cor:griffiths} and the assumed pullback identity
give semiampleness by Lemma~\ref{lem:birational-sections}.
For an everywhere defined period map, proper GAGA algebraizes
the holomorphic Griffiths line; its pullback is the Schmid line
with no boundary correction. If the period map is quasi-finite,
Griffiths curvature gives positive degree on every complete
integral curve. A positive-dimensional proper fibre contains a
curve, as seen on a projective modification. The semiample
contraction therefore has finite fibres, so it is finite and
the line is ample.
\end{proof}

\begin{remark}\label{rem:singular-extension}
The open restriction does not determine an extension: both
$\cO_{\mathbb P^1}$ and $\cO_{\mathbb P^1}(-\infty)$ restrict to the
trivial line on $\mathbb A^1$, but only the former is semiample.
The identity \eqref{eq:normal-source-comparison}, including
its interior comparison, excludes such boundary twists on a resolution.
Likewise, quasi-unipotence does not give the unipotent pullback formula
for an arbitrary Deligne extension; the comparison in
Corollary~\ref{cor:normal-alteration} is required.
\end{remark}

Griffiths proved projectivity of period images with complete smooth
base \cite[Theorem 9.7]{Griffiths70}. The following statement records
the descent of the Griffiths line on a reduced proper algebraic space.

\begin{theorem}[Period maps defined on reduced proper spaces]
\label{sg:everywhere-period}
Let $X$ be a reduced proper algebraic space and let
\[
\Phi:X^{\an}\longrightarrow\Gamma\backslash\mathcal D
\]
be an everywhere defined locally liftable period map for integral
polarizable pure Hodge structures, with $\Gamma$ neat.
Its Griffiths line is algebraic and semiample on $X$.
If $\Phi$ is quasi-finite, the Griffiths line is ample and $X$ is projective.
No normality or seminormality assumption is required.
\end{theorem}
\begin{proof}
Proper GAGA algebraizes the Griffiths line $\Lambda_X$ and
the closed analytic relation
\[
 R=\{(x,x'):\Phi(x)=\Phi(x')\}\subset X^{\an}\times X^{\an}.
\]
A local frame of the Griffiths line near a period point pulls back
to a frame on a saturated neighbourhood of the corresponding
$R$-class.
For an integral subspace $Z\subset X$ transverse to $R$, we take
a smooth projective resolution. Its period map is generically
immersive, so the pulled-back line is nef and big by Griffiths
curvature. Bigness descends to $Z$, and
Theorem~\ref{thm:frame-criterion} proves that $\Lambda_X$ is semiample,
including across the full conductor.

If $\Phi$ is quasi-finite, positivity on every complete curve makes
the semiample contraction finite, as in
Corollary~\ref{cor:singular-griffiths}. Hence $\Lambda_X$ is ample
and $X$ is projective.
\end{proof}

\begin{example}[The K3 Hodge line]\label{cor:k3-line}
In Theorem~\ref{sg:everywhere-period}, assume that the variation has
weight two and Hodge numbers $(1,h^{1,1},1)$. Let $\lambda=F^2E$
be its highest line and $\delta=\det E$. Then
\[
 \Lambda=\delta\otimes\lambda^2,\qquad
 \delta^2\simeq\cO_X,\qquad \Lambda^2\simeq\lambda^4.
\]
Indeed, the polarization identifies $E/F^1E\simeq\lambda^\vee$
and $E\simeq E^\vee$, so $\det F^1E=\delta\otimes\lambda$;
proper GAGA algebraizes these identities. By
Theorem~\ref{sg:everywhere-period}, $\lambda$ is semiample and has
the same connected contraction as $\Lambda$. The identities also
hold for Schmid extensions on a smooth log pair. On a nonnormal
boundary model, they must hold on the full conductor to apply
Theorem~\ref{sg:conductor-target}.
\end{example}

\subsection{Normal compactifications and extension}

\begin{lemma}
\label{app:polarized-extension}
Let $X$ be smooth projective, let $Z$ be normal projective, and let
$q:X\to Z$ be a projective surjection. For an open subset $Y\subset Z$,
put $U=q^{-1}(Y)$ and assume that $X\setminus U$ is a reduced SNC divisor.
Let $Y$ carry a locally liftable pure period map whose pullback to $U$
has unipotent monodromy along $X\setminus U$. Let $\mathcal L_Y$
be its highest or Griffiths line, and let $\mathcal L_X$ be the
corresponding Schmid extension on $X$. Suppose that an integer $a>0$
and an ample line $A$ on $Z$ admit isomorphisms
\[
 \mathcal L_X^a\simeq q^*A,\qquad A|_Y\simeq\mathcal L_Y^a
\]
which agree after pullback to $U$.
Let $T$ be smooth, algebraic or complex analytic, let $E\subset T$
be a reduced SNC divisor, and let $h:T\setminus E\to Y$ be a map whose
pulled-back variation has unipotent monodromy along $E$. Denote its
Schmid highest or Griffiths line by $\mathcal L_T$.
If $h$ extends rationally, respectively meromorphically, to $Z$, it
extends uniquely to a morphism, respectively holomorphic map,
$\overline h:T\to Z$, and the open comparison extends to
$\overline h^*A\simeq\mathcal L_T^a$.
\end{lemma}
\begin{proof}
We work \'etale-locally on $T$ in the algebraic case and locally
in the analytic case. We resolve the graph and boundary, leaving
$T\setminus E$ unchanged:
\[
 p:T'\to T\text{ projective modification},\qquad
 \psi:T'\to Z.
\]
Here $T'$ is smooth and $(p^{-1}E)_{\red}$ is SNC.
Locally on $T'$, we choose a component of $T'\times_ZX$ dominating
$T'$, intersect it with general relative hyperplanes, and resolve
its boundary. We obtain
\[
 \rho:W\to T'\text{ proper generically finite},\qquad
 \eta:W\to X,\qquad q\eta=\psi\rho,
\]
with $W$ smooth and the boundary SNC. Over $T\setminus E$ both
variations are pulled back from that on $Y$, even if the image lies in
$Y_{\rm sing}$. Unipotent Schmid functoriality gives
\[
 \rho^*\psi^*A\simeq\eta^*\mathcal L_X^a
       \simeq\rho^*p^*\mathcal L_T^a.
\]
We apply the trace divided by the degree of the finite Stein factor
of $\rho$ to extend the interior identification locally on $T'$.
Since $T'$ is normal, the extension is regular, respectively
holomorphic. The same construction extends the inverse; the two
composites equal the identity on the dense interior.
The local isomorphisms agree there and hence glue to
\begin{equation}\label{app:compact-line-comparison}
 \psi^*A\simeq p^*\mathcal L_T^a.
\end{equation}
Every curve in a $p$-fibre has $\psi^*A$-degree zero and is contracted
by $\psi$. These fibres are connected projective and connected
by curve chains, so $\psi$ is constant on each fibre. The image
of $(p,\psi)$ in $T\times Z$ is proper over $T$ and has one point over
each point of $T$. Its projection to $T$ is finite and birational,
hence an isomorphism because $T$ is normal. We obtain $\overline h$.
Using $p_*\mathcal O_{T'}=\mathcal O_T$ and the projection formula,
we descend the line comparison to $T$. These local extensions
agree on the dense interior and hence on their overlaps.
\end{proof}

We use the moderate-growth convention of \cite[Definition 5.1]{BFMT}.
On a smooth log model with boundary coordinates $z_1,\ldots,z_r$,
a holomorphic section $s$ has moderate growth if, for every
$\varepsilon>0$, locally
\[
 \|s\|\le C_\varepsilon\prod_{i=1}^r|z_i|^{-\varepsilon}.
\]
This sub-polynomial condition differs from Deligne's polynomial
moderate growth, which allows poles.

For a normal quasi-projective $Y$ carrying a locally liftable period
map to a neat quotient, we give the holomorphic Griffiths line an
algebraic structure as follows. We choose a projective resolution $r:U\to Y$
and a smooth projective SNC compactification $X\supset U$. Neatness
makes the boundary monodromy unipotent, and proper GAGA algebraizes
the Schmid line $\Lambda_X$. Set
\[
\Lambda_Y=r_*(\Lambda_X|_U).
\]
Proper direct-image comparison and $r_*\cO_U=\cO_Y$ identify
$\Lambda_Y^{\an}$ with the holomorphic Griffiths line: on $U^{\an}$
the latter pulls back to $\Lambda_X|_U$ by Schmid functoriality.
Thus $\Lambda_Y$ is an algebraic line bundle and
$r^*\Lambda_Y\simeq\Lambda_X|_U$. A common smooth logarithmic refinement
and proper direct-image comparison make this algebraic structure and
these comparisons independent of the choices.

\begin{theorem}\label{thm:normal-compactification}
Let $Y$ be a normal quasi-projective variety with a locally liftable
quasi-finite period map for integral polarizable pure Hodge structures
and neat monodromy. Write $\Lambda_Y$ for its Griffiths line with the
algebraic structure just described.
There exist a normal projective compactification $Y\subset Y^{\rm BB}$,
an integer $a>0$, and an ample line $A\in\operatorname{Pic}(Y^{\rm BB})$
such that $A|_Y\simeq\Lambda_Y^a$.
For a projective resolution $r:U\to Y$ and a smooth projective SNC
compactification $X\supset U$,
\[
Y^{\rm BB}\simeq\operatorname{Proj}
\bigoplus_{m\ge0}H^0(X,\Lambda_X^{am}),
\]
and the induced morphism $X\to Y^{\rm BB}$ is the connected
contraction of $\Lambda_X$. The compactification and the ample rational
polarization $A^{1/a}$ are independent of the choices.
Let $(T,E)$ be a log smooth pair with $T$ a separated algebraic space
of finite type. Every algebraic morphism $h:T\setminus E\to Y$
extends uniquely to $\overline h:T\to Y^{\rm BB}$, with
\[
\overline h^*A\simeq\Lambda_T^a,
\]
where $\Lambda_T$ is the Schmid-extended Griffiths line of the
variation pulled back by $h$.
\end{theorem}

\begin{proof}
\emph{Step 1. The compactification.}
Neatness and the monodromy theorem give unipotent boundary monodromy.
Corollary~\ref{cor:griffiths} gives a connected contraction
$q:X\to Z$, with $Z$ normal projective and
$\Lambda_X^a\simeq q^*A$ for an ample line $A$.
Quasi-finiteness and Griffiths curvature make $\Lambda_X$ big,
so $q$ is birational.

Suppose a complete irreducible curve $C\subset X$ meets $U$ and
has $\Lambda_X$-degree zero. The positive extended curvature current
on its normalization \cite{FF} has zero mass, hence vanishes.
The period map is therefore constant on $C\cap U$, and
quasi-finiteness gives $r(C\cap U)=\{y\}$. Since $r^{-1}(y)$ is
proper, $C\subset r^{-1}(y)\subset U$. Conversely, curves in an
$r$-fibre have degree zero. Curve chains in the connected projective
fibres show that $U$ is saturated and that $r$ and $q|_U$ have the
same fibres.

Since $U$ is saturated, $q(U)$ is open in $Z$. The image
of $(r,q|_U)$ in $Y\times q(U)$ projects properly and bijectively
to both factors. These projections are finite and birational; both
factors are normal, so they are isomorphisms. Hence $Y\simeq q(U)$.
The projection formula gives $A|_Y\simeq\Lambda_Y^a$,
and we set $Y^{\rm BB}=Z$.

\emph{Step 2. The section ring.}
On a common smooth projective model, unipotent Schmid pullback
identifies the two lines. Lemma~\ref{lem:birational-sections} therefore
identifies their graded section rings in common divisible degrees,
and hence their Proj and ample rational polarizations.
For $s\in H^0(Y^{\rm an},\Lambda_Y^{am})$, the growth condition
is equivalent to extension on the log model:
\[
s\in H^0_{\rm mg}(Y,\Lambda_Y^{am})
\quad\Longleftrightarrow\quad
r^*s\text{ extends to }\Lambda_X^{am}\text{ on }X.
\]
Indeed, in a boundary polydisc with coordinates $z_1,\ldots,z_r$,
the unipotent norm estimates \cite{Schmid,CKS} give, for a local
nonvanishing Schmid frame $e$ and some $C,N>0$,
\[
 \|e\|+\|e\|^{-1}\le C\prod_{i=1}^r(1+|\log|z_i||)^N.
\]
To obtain a uniform estimate on the punctured polydisc, we use finitely
many angular sectors and the finitely many orderings of the quantities
$-\log|z_i|$; the sectorwise estimates give the displayed bound after
increasing $C$ and $N$. Hence an extended section has sub-polynomial norm. Conversely, if
$r^*s=fe$ has sub-polynomial norm, then for every $\varepsilon>0$
the holomorphic coefficient satisfies
$|f|\le C_\varepsilon\prod_i|z_i|^{-\varepsilon}$ on a smaller
polydisc. Cauchy's formula eliminates all negative Laurent
coefficients, so $f$ extends holomorphically. This also proves that
the growth condition is independent of the chosen log model.
Then, by proper GAGA and $r_*\mathcal O_U=\mathcal O_Y$,
\[
\bigoplus_{m\ge0}H^0(X,\Lambda_X^{am})
\simeq\bigoplus_{m\ge0}H^0_{\rm mg}(Y,\Lambda_Y^{am}).
\]

Finally, an algebraic $h:T\setminus E\to Y$ defines a rational map to $Z$.
The pulled-back monodromy is unipotent by
neatness. Lemma~\ref{app:polarized-extension} gives the unique
extension and the line comparison.
\end{proof}

\begin{corollary}\label{cor:compactification}
Let $Y$ be a normal separated algebraic space of finite type carrying
a locally liftable pure CY period map to a neat quotient. Choose a
proper birational $r:U\to Y$, with $U$ smooth quasi-projective, and a
smooth projective log compactification $X\supset U$. Let $\mathcal L_X$
be the Schmid highest line, and equip the holomorphic highest line
on $Y$ with the algebraic structure
$\mathcal L_Y=r_*(\mathcal L_X|_U)$, as above.
Assume that $\mathcal L_X$ is numerically integrable and has torsion
combinatorial monodromy. For every smooth connected curve $C^\circ$
and every nonconstant algebraic map $c:C^\circ\to Y$, let
$\overline C$ be its smooth projective completion and let
$\mathcal L_{\overline C}$ be the Schmid extension of the pulled-back
highest Hodge line. Assume that
\[
\deg_{\overline C}\mathcal L_{\overline C}>0
\]
for every such $c$.

Then the connected contraction $q:X\to Y^{\rm BB}$ of $\mathcal L_X$
is a normal projective compactification of $Y$, with
\[
 \mathcal L_X^a\simeq q^*A,\qquad A|_Y\simeq\mathcal L_Y^a,
 \qquad a>0,\quad A\text{ ample}.
\]
The section-ring description, independence of choices and algebraic
extension property of Theorem~\ref{thm:normal-compactification}
hold with $\Lambda$ replaced by the highest line $\mathcal L$.
In particular, the extended map satisfies
$\overline h^*A\simeq\mathcal L_T^a$, where $\mathcal L_T$ is the
Schmid highest line of the pulled-back variation.
The hypotheses hold for the Griffiths line of a locally liftable
quasi-finite period map.
\end{corollary}

\begin{proof}
Theorem~\ref{thm:main} gives the connected contraction
$q:X\to Z$, with $Z$ normal projective and $q^*A\simeq\mathcal L_X^a$
for an ample line $A$. For a complete curve $C\subset X$ meeting
$U$, the curvewise positivity hypothesis and Schmid pullback give
\[
 \mathcal L_X\cdot C=0\quad\Longrightarrow\quad
 r(C\cap U)=\{y\}\quad\Longrightarrow\quad C\subset r^{-1}(y)\subset U;
\]
the last implication uses properness of $r$. Conversely, the line
has degree zero on every $r$-fibre. The saturation, graph and
common-model arguments of Theorem~\ref{thm:normal-compactification}
now identify $Y=q(U)\subset Z$ and prove independence and the
Proj description.
Lemma~\ref{app:polarized-extension} gives extension and its line
comparison. For the Griffiths line, quasi-finiteness makes
the curvature positive on a nonempty open of every nonconstant
curve; its positive extended current gives positive degree.
Corollary~\ref{cor:griffiths} supplies the remaining hypotheses.
\end{proof}

\begin{corollary}\label{cor:normal-functoriality}
For $i=1,2$, let $Y_i^{\rm BB}$ be the compactification of $Y_i$
in Theorem~\ref{thm:normal-compactification} or
Corollary~\ref{cor:compactification}, and let $\mathbb V_i$ be the
corresponding variation. Let $g:Y_1\to Y_2$ be an algebraic morphism
with an isomorphism $\mathbb V_1\simeq g^*\mathbb V_2$.
Assume that both compactifications are defined by Griffiths lines,
or both by highest Hodge lines. Then $g$ extends uniquely to a finite map
\[
g^{\rm BB}:Y_1^{\rm BB}\longrightarrow Y_2^{\rm BB}.
\]
After choosing a common sufficiently divisible power $a$, the ample
extensions satisfy $(g^{\rm BB})^*A_2\simeq A_1$. These extensions
respect identities and composition.
\end{corollary}
\begin{proof}
On a smooth projective log model $q_1:X_1\to Y_1^{\rm BB}$, with
$r_1:U_1\to Y_1$ over the interior, the extension property gives
\[
h:X_1\longrightarrow Y_2^{\rm BB},\qquad
h|_{U_1}=gr_1,\qquad h^*A_2\simeq\mathcal L_{X_1}^a=q_1^*A_1.
\]
Ampleness of $A_2$ makes $h$ constant on each connected projective
$q_1$-fibre. Since $(q_1)_*\cO_{X_1}=\cO_{Y_1^{\rm BB}}$,
$h=g^{\rm BB}q_1$ and $(g^{\rm BB})^*A_2\simeq A_1$. Ampleness
of $A_1$ makes this proper map finite. Uniqueness follows from
density and separatedness, and gives compatibility with composition.
\end{proof}

\begin{lemma}\label{app:singular-meromorphic-extension}
Let $Y$ be a reduced quasi-projective variety with a locally liftable
quasi-finite period map to a neat quotient for integral polarizable
pure Hodge structures, and let
$\overline Y$ be a projective compactification. Every holomorphic map
\[
h:\Delta^p\times(\Delta^*)^q\longrightarrow Y^{\rm an}
\]
extends meromorphically to $\Delta^{p+q}$ with values in
$\overline Y^{\rm an}$.
\end{lemma}
\begin{proof}
We first consider a nonconstant map $\gamma:\Delta^*\to Y^{\rm an}$.
Let $Z\subset Y$ be the reduced Zariski closure of its image. We
resolve the closure of $Z$ in $\overline Y$:
\[
\pi:\widetilde Z\longrightarrow\overline Z,\qquad
V=\pi^{-1}Z,\qquad \widetilde Z\setminus V\text{ SNC}.
\]
The map $\gamma$ lifts to $\widetilde\gamma:\Delta^*\to V$.
Indeed, it lifts off a discrete set, and the dominant component of the
pulled-back graph is proper and birational over the smooth curve
$\Delta^*$; its normalization is $\Delta^*$.
The lift is Zariski dense in $\widetilde Z$. The pulled-back period map
is immersive at some point of $V$, since the original map has
zero-dimensional fibres and $\pi$ is birational.

By \cite[Theorem 1.4]{Deng}, there are a Finsler metric $h_F$,
positive definite on a dense Zariski open $V^0\subset V$, and a
K\"ahler form $\omega$ on $\widetilde Z$ such that
\[
\mathrm d\mathrm d^c\log|\widetilde\gamma'|_{h_F}^2
\geq\widetilde\gamma^*\omega.
\]
Zariski density gives $\widetilde\gamma(\Delta^*)\cap V^0\ne\varnothing$,
so $|\widetilde\gamma'|_{h_F}\not\equiv0$.
The Big Picard criterion \cite[Theorem 5.4]{Deng} extends
$\widetilde\gamma$ holomorphically to $\Delta$.
Composing with $\pi$ extends $\gamma$ to $\overline Y$.
Constant maps have the same conclusion.

Following \cite[proof of Proposition 5.2]{Deng}, we embed
$\overline Y$ in $\mathbb P^N$ and apply Siu's extension theorems.
Near a smooth point of the missing coordinate divisor, each
transverse punctured-disc map extends to $\mathbb P^N$ by the
one-variable case above. Siu's separate extension
theorem \cite[p.~442, $(*)$]{Siu75} gives meromorphic extension there.
His codimension-two extension theorem \cite[Theorem 1]{Siu75} then
extends across the intersections of the coordinate divisors.
The extended graph lies in $\Delta^{p+q}\times\overline Y$:
this product is closed in $\Delta^{p+q}\times\mathbb P^N$ and
contains the graph on the dense original domain.
\end{proof}

\begin{corollary}\label{cor:normal-picard}
Let $Y^{\rm BB}$ be the Griffiths compactification of
Theorem~\ref{thm:normal-compactification}, or the highest-line
compactification of Corollary~\ref{cor:compactification} when its period
map is quasi-finite. Let $\mathcal L$ denote the corresponding Hodge
line and $A$ its ample extension in degree $a$. Every holomorphic map
\[
h:\Delta^p\times(\Delta^*)^q\longrightarrow Y^{\rm an}
\]
extends uniquely to a holomorphic map
$\overline h:\Delta^{p+q}\to(Y^{\rm BB})^{\rm an}$. Moreover,
\[
\overline h^*A\simeq\mathcal L_{\Delta^{p+q}}^a
\]
extends the given comparison on the punctured polydisc.
\end{corollary}

\begin{proof}
Lemma~\ref{app:singular-meromorphic-extension} gives a meromorphic
extension to $\Delta^{p+q}$. Neatness makes the pulled-back local
monodromies unipotent. Applying Lemma~\ref{app:polarized-extension}
to the polarized log model of the compactification gives the
holomorphic extension and line comparison; density gives uniqueness.
\end{proof}

The meromorphic extension uses the full period map; the ample highest
or Griffiths line then contracts exceptional fibres. Descent across a
nonnormal conductor is treated in Section~\ref{sec:nonnormal-descent}.

We next remove neatness by finite-level descent. For a test map from
a log smooth space $T$, write $\Lambda_T^{[a]}$ for the descent of the
$a$th power of the unipotent Schmid line on local finite covers
trivializing the chosen neat level. The proof below establishes this
descent, including at ramification points, and its independence of the
local covers. This notation does not mean the $a$th power of an
arbitrarily chosen Deligne extension on $T$.

\begin{theorem}[Normal compactifications without neatness]
\label{thm:normal-level-compactification}
Let $Y$ be a normal quasi-projective variety with a locally liftable
quasi-finite period map $\Phi:Y^{\an}\to\Gamma\backslash D$ for
polarized integral pure Hodge structures, where $\Gamma$ is a discrete
subgroup of the integral isometry group of the fixed polarization.
There are a normal projective compactification $\overline Y$, an
integer $a>0$ and an ample line $A$ whose restriction is the descended
$a$th Griffiths power $\Lambda_Y^{[a]}$. For this power,
\[
\bigoplus_{m\ge0}H^0(\overline Y,A^m)
 =\bigoplus_{m\ge0}H^0_{\rm mg}(Y,\Lambda_Y^{[am]}),
\qquad
\overline Y=\operatorname{Proj}
 \bigoplus_{m\ge0}H^0_{\rm mg}(Y,\Lambda_Y^{[am]}).
\]
Moderate growth is measured by the invariant Hodge metric as above.
A section has moderate growth if and only if its pullback to a neat
level cover extends to the corresponding Schmid line on a smooth
projective log model. The growth condition,
the compactification and its ample rational polarization $A^{1/a}$
are independent of these choices and of the divisible exponent.

Every holomorphic map $h:\Delta^p\times(\Delta^*)^q\to Y^{\an}$
extends uniquely to a
holomorphic map $\overline h:\Delta^{p+q}\to\overline Y^{\an}$.
For a log smooth pair $(T,E)$ with $T$ a separated algebraic space
of finite type, every algebraic map $h:T\setminus E\to Y$ extends
uniquely to an algebraic morphism $\overline h:T\to\overline Y$.
In both cases the descended Schmid line $\Lambda_T^{[a]}$ exists
and its comparison on $T\setminus E$ extends to
\[
\overline h^*A\simeq\Lambda_T^{[a]},
\]
where $T=\Delta^{p+q}$ in the analytic case.
\end{theorem}
\begin{proof}
\emph{Step 1. A neat cover.}
We intersect $\Gamma$ with a normal neat congruence subgroup of the
integral isometry group to obtain a normal neat subgroup
$\Gamma_0\subset\Gamma$ of finite index. We set $G=\Gamma/\Gamma_0$
and normalize the reduced analytic fibre product
\[
Y^{\an}\times_{\Gamma\backslash D}(\Gamma_0\backslash D).
\]
Locally a lift to $D$ writes this fibre product as a finite union of
graphs over $Y^{\an}$. Since $Y$ is normal, normalizing the fibre product
gives disjoint copies of the local base. We obtain a finite
\'etale cover, algebraized by Riemann existence as
\[
\pi:Y_0\longrightarrow Y,\qquad Y_0/G=Y.
\]
The cover may be disconnected. We take a $G$-equivariant smooth
projective log model $X_0$ of $Y_0$. By unipotent Schmid functoriality,
its algebraic Griffiths line $\Lambda_{X_0}$ is $G$-linearized.
Theorem~\ref{thm:normal-compactification} therefore gives an
equivariant contraction
\[
q_0:X_0\longrightarrow\overline Y_0,
\qquad q_0^*A_0=\Lambda_{X_0}^{b}.
\]
with equivariant line comparison after choosing a sufficiently
divisible Proj twisting sheaf. The quotient
$\overline Y=\overline Y_0/G$ is normal projective and contains
$Y_0/G=Y$ as an open subvariety.

\emph{Step 2. The polarization and section ring.}
After replacing $b$ by a multiple, $A_0$ descends to a line $A$ on
$\overline Y$: a power divisible by $|G|$ kills the stabilizer
characters on its fibres. Invariant local frames then give descent,
and finite-surjective descent of ampleness makes $A$ ample.
Let $\overline\pi:\overline Y_0\to\overline Y$ be the quotient map
and let $a$ be the final exponent. In each degree, descent gives
\[
 H^0(\overline Y,A^m)
   \simeq H^0(\overline Y_0,\overline\pi^*A^m)^G.
\]
By Theorem~\ref{thm:normal-compactification}, the right side is
$H^0_{\rm mg}(Y,\Lambda_Y^{[am]})$. A common neat refinement and
log model identify the rings in common divisible degrees, proving
independence of the compactification and rational polarization.

\emph{Step 3. Extension of maps and Schmid comparison.}
For a punctured-polydisc map, we pull back $Y_0\to Y$. A finite cover of
$\Delta^p\times(\Delta^*)^q$ is trivialized by a sufficiently divisible
coordinate-power cover, since its fundamental group is $\Z^q$.
We take a lift on that cover and extend it to $\overline Y_0$
by Corollary~\ref{cor:normal-picard}. Its composite with
$\overline Y_0\to\overline Y$ is invariant under the coordinate deck
group on the dense punctured open, hence everywhere. It descends to
the original polydisc. Density proves uniqueness. Write
$\sigma:S'\to S$ for the coordinate-power cover and
$\overline h_0:S'\to\overline Y_0$ for the extended lift. The neat
comparison gives
\[
\sigma^*\overline h^*A\simeq\overline h_0^*A_0
                    \simeq\Lambda_{S'}^a,
\]
where $A_0$ has been replaced by the power descending to $A$.
These comparisons respect the descent data on the punctured open.
Functoriality of the unipotent Schmid extension extends the deck
comparisons, and equality on the dense open makes the displayed
isomorphism equivariant everywhere. Its left side is pulled back from
$S$, so every stabilizer at a ramification point acts trivially on
the right side as well. Finite-group descent therefore gives a line
$\Lambda_S^{[a]}$ and the comparison
$\overline h^*A\simeq\Lambda_S^{[a]}$. Common coordinate-power
refinements identify these descents on overlaps and show that they
are independent of the local covers.

For an algebraic map $T\setminus E\to Y$, we apply this analytic
extension locally along $E$. The local extensions glue uniquely. The closure
of the algebraic graph in $T\times\overline Y$ analytifies to its
holomorphic graph. Its projection to $T$ is therefore an isomorphism,
so the extension is algebraic. This argument is \'etale-local on $T$.
\end{proof}

\subsection{Conductor descent}\label{sec:nonnormal-descent}

Let $X$ be a reduced proper algebraic space, let
$\nu:X^\nu\to X$ be its finite normalization, and let
$\mathcal L\in\operatorname{Pic}(X)$. Set
\[
\mathfrak c=\operatorname{Ann}_{\cO_X}(\nu_*\cO_{X^\nu}/\cO_X),
\quad D=V(\mathfrak c),\quad D^\nu=X^\nu\times_XD,
\quad\nu_D:D^\nu\to D.
\]
The subschemes $D,D^\nu$ retain their full scheme structures.

\begin{proposition}[Conductor descent of sections]
\label{prop:conductor-sections}
For every $m\ge0$,
\[
H^0(X,\mathcal L^m)=
H^0(X^\nu,\nu^*\mathcal L^m)
\mathop{\times}_{H^0(D^\nu,\nu^*\mathcal L^m|_{D^\nu})}
H^0(D,\mathcal L^m|_D).
\]
Consequently, $\mathcal L$ is semiample exactly when, for some $m>0$,
the projection of this fibre product to the first factor contains a
finite base-point-free system on $X^\nu$.
\end{proposition}
\begin{proof}
On an \'etale affine chart $R\subset S=\operatorname{Norm}(R)$,
\[
\mathfrak c=\{r\in R:rS\subset R\},\qquad
R=S\times_{S/\mathfrak c}R/\mathfrak c.
\]
After trivializing $\mathcal L^m$, we glue these identities and take
global sections.
For a compatible section $s$,
\[
(\nu^*s)(x^\nu)\ne0\iff s(\nu(x^\nu))\ne0.
\]
Surjectivity of $\nu$ proves the generation criterion.
\end{proof}

\begin{theorem}[Finite descent on the conductor image]
\label{sg:conductor-target}
Let $a$ be a positive integer, let $Y^\nu$ be a projective scheme,
and let $A\in\operatorname{Pic}(Y^\nu)$ be ample. Suppose that a
morphism $q^\nu:X^\nu\to Y^\nu$ is equipped with an isomorphism
\[
\nu^*\mathcal L^a\simeq(q^\nu)^*A.
\]
Let $E\subset Y^\nu$ be a closed subscheme through which
$q^\nu|_{D^\nu}$ factors, say $u:D^\nu\to E$.
Assume that a proper algebraic space $T$, a finite surjection $f:E\to T$,
a morphism $h:D\to T$, and $B\in\operatorname{Pic}(T)$ satisfy
\[
fu=h\nu_D,\qquad A|_E\simeq f^*B,\qquad
\mathcal L^a|_D\simeq h^*B.
\]
Assume also that, on the full subscheme $D^\nu$, the composite
\[
 (\nu^*\mathcal L^a)|_{D^\nu}
 \simeq u^*(A|_E)\simeq(fu)^*B
 =(h\nu_D)^*B\simeq\nu_D^*(\mathcal L^a|_D)
\]
is the conductor identification of $\mathcal L^a$.
Then $\mathcal L$ is semiample.
No reducedness assumption on $D,D^\nu,E$ or $T$ is required.
\end{theorem}
\begin{proof}
Since $f$ is finite surjective and $f^*B\simeq A|_E$ is ample,
$B$ is ample. For $n\gg0$, Serre vanishing gives
\[
\begin{gathered}
B^n\text{ generated},\qquad \cI_EA^n\text{ generated},\\
H^0(Y^\nu,A^n)\twoheadrightarrow H^0(E,A^n|_E).
\end{gathered}
\]
We lift generators $t_j\in H^0(T,B^n)$ to
\[
s_j\in H^0(Y^\nu,A^n),\qquad s_j|_E=f^*t_j.
\]
Together with sections $v_k$ generating $\cI_EA^n$, they give
compatible pairs
\[
\bigl((q^\nu)^*s_j,h^*t_j\bigr),\qquad
\bigl((q^\nu)^*v_k,0\bigr).
\]
The first system generates over $E$; the second over $Y^\nu\setminus E$.
Their pullbacks generate $\mathcal L^{an}$ on $X$ by
Proposition~\ref{prop:conductor-sections}.
\end{proof}

\begin{corollary}[A conductor criterion for compactification]
\label{cor:conductor-compactification}
In the notation of this subsection, assume that $X$ is projective
and $Y\subset X$ is a dense open. Suppose
\[
q^\nu:X^\nu\to Z,\qquad (q^\nu)_*\cO_{X^\nu}=\cO_Z,
\qquad (\nu^*\mathcal L)^a\simeq(q^\nu)^*A
\]
with $Z$ normal projective and $A$ ample. Assume $Z$ contains $Y^\nu$
as an open, $(q^\nu)^{-1}(Y^\nu)=\nu^{-1}Y$, and $q^\nu$ is the
identity on this open under the normalization identification.
If the conductor data satisfy Theorem~\ref{sg:conductor-target}, then
the connected semiample contraction of $\mathcal L$ is a projective
compactification $q:X\to\overline Y$ of the original $Y$. Its
normalization is $Z$, and a positive power of $\mathcal L$ is the
pullback of an ample line on $\overline Y$.
\end{corollary}
\begin{proof}
Theorem~\ref{sg:conductor-target} gives the Stein contraction
$q:X\to\overline Y$, with $q_*\cO_X=\cO_{\overline Y}$ and
$\mathcal L^m=q^*B$ for an ample $B$. A curve $C\subset X$ meeting
$Y$ lifts to a curve meeting $Y^\nu$, which $q^\nu$ does not contract.
Hence $\mathcal L\cdot C>0$. A connected $q$-fibre meeting $Y$
is therefore supported at one point. Thus $Y$ is saturated and
$Y\to q(Y)$ is finite, with direct image $\cO_{q(Y)}$; it is an
isomorphism onto an open subset.

The map $q\nu$ contracts the connected projective $q^\nu$-fibres,
so it factors through a proper map $\lambda:Z\to\overline Y$ with
\[
 q\nu=\lambda q^\nu,\qquad \lambda^*B^a\simeq A^m.
\]
Ampleness of $A$ makes $\lambda$ finite. Its restriction is
$Y^\nu\to Y$, so the normal scheme $Z$ is the normalization
of $\overline Y$.
\end{proof}

Theorem~\ref{sg:conductor-target} requires $T$, $B$ and compatible
identifications on the full conductor; it does not construct these
data. The next examples show why the normalized line alone does not suffice.

\begin{example}[Normalization alone is insufficient]
\label{ex:nodal-descent}
For an integral rational nodal curve $C$, let $p,q\in\mathbb P^1$ lie
over its node. We glue the trivial line by $\lambda\in\C^*$ of infinite order:
\[
\nu^*L_\lambda=\cO_{\mathbb P^1},\qquad
H^0(C,L_\lambda^m)=\{c\in\C:c=\lambda^m c\}=0\quad(m>0).
\]
Thus semiampleness does not descend from the normalization: the
gluing character has infinite order.
\end{example}

\begin{example}[A cuspidal obstruction invisible to nodal curves]
\label{sg:cuspidal-obstruction}\label{ex:cuspidal-line-descent}
Let $C=\{y^2z=x^3\}\subset\mathbb P^2$. Its normalization is
$\nu:\mathbb P^1\to C$, with
\[
D=\Spec\C,\qquad D^\nu=\Spec\C[\epsilon]/(\epsilon^2).
\]
We glue the trivial lines on $\mathbb P^1$ and $D$ by
$1+\alpha\epsilon$, with $\alpha\ne0$. The resulting $L_\alpha$ satisfies
\[
\nu^*L_\alpha=\cO_{\mathbb P^1},\qquad
H^0(C,L_\alpha^m)
=\{(c,d)\in\C^2:c=d(1+m\alpha\epsilon)\}=0\quad(m>0).
\]
Every morphism from a reduced nodal curve $T$ factors through $\nu$, since
$\mathbb P^1=C^{\mathrm{sn}}$ and $T$ is seminormal. Thus
\[
g:T\to C,\quad T\text{ nodal}
\Longrightarrow g^*L_\alpha\simeq\cO_T.
\]
Even the combination of trivial normalized Schmid line and all nodal
torsion tests therefore does not imply semiampleness on a nonseminormal
space. The missing condition is the infinitesimal conductor gluing.
\end{example}

\begin{theorem}[Reduced curves]\label{sg:curve-criterion}
Let $C$ be a reduced projective curve and $L\in\operatorname{Pic}(C)$.
Let $C_0$ be the reduced union of the irreducible components on which
$L$ has degree zero. Then
\[
L\text{ semiample}\quad\Longleftrightarrow\quad
\left\{\begin{array}{l}
\deg_{C_i}L\ge0\text{ for every irreducible component }C_i,\\
L|_{C_0}\text{ is torsion}.
\end{array}\right.
\]
The second condition is empty when $C_0=\varnothing$.
No nodality or seminormality hypothesis is needed.
\end{theorem}
\begin{proof}
If $L^m=\varphi^*\cO(1)$, every component has nonnegative degree,
and $\varphi$ contracts each connected component of $C_0$ to a point.
Hence $L^m|_{C_0}\simeq\cO_{C_0}$.
For the converse, we write $C=C_0\cup C_+$ with reduced structures
and put $Z=C_0\cap C_+$ scheme-theoretically. Then
\[
\cO_C=\cO_{C_0}\mathop{\times}_{\cO_Z}\cO_{C_+},
\qquad\dim Z=0,\qquad L|_{C_+}\text{ ample}.
\]
The cases $C_0=\varnothing$ or $C_+=\varnothing$ are immediate.
For a sufficiently divisible $m$,
\[
L^m|_{C_0}\simeq\cO_{C_0},\quad
H^1(C_+,\cI_ZL^m)=0,\quad\cI_ZL^m\text{ generated on }C_+.
\]
We lift $1|_Z$ to $s_+\in H^0(C_+,L^m)$. The section $(1,s_+)$,
together with $(0,v_k)$ for generators of $\cI_ZL^m$, generates $L^m$
on $C$. The matching takes place on the full finite scheme $Z$,
including its nilpotents.
\end{proof}

\begin{corollary}\label{sg:nodal-hodge}
Let $C,L,C_0$ be as in Theorem~\ref{sg:curve-criterion}, with $C$ nodal,
and write $\nu:\coprod_i C_i\to C$ for its normalization. Suppose
each $C_i$ carries an integral polarizable pure CY variation on a
dense open, with unipotent boundary monodromy and Schmid highest
line $\mathcal H_i$. Fix isomorphisms
\[
(\nu^*L)|_{C_i}\simeq\mathcal H_i^{r_i}
\]
with $r_i>0$. If the canonical finite connections on the torsion
component lines, glued by the node identifications of $L$, have
finite monodromy on every connected component of $C_0$, then $L$
is semiample.
\end{corollary}
\begin{proof}
Semipositivity gives nonnegative degree on every component.
On a component of degree zero, the flat highest-line theorem
makes the normalized line torsion. The node hypothesis gives finite holonomy on $C_0$; hence
$L^m|_{C_0}\simeq\cO_{C_0}$ for some $m>0$.
Theorem~\ref{sg:curve-criterion} proves semiampleness.
\end{proof}

\begin{proposition}[Seminormal spaces with finite conductor image]
\label{prop:seminormal-finite-conductor}
Let $X$ be a reduced seminormal proper algebraic space over $\C$,
let $\nu:X^\nu\to X$ be its finite normalization, and let $L$ be a line
bundle on $X$. Suppose $\nu^*L$ is semiample. Write its connected
contraction, componentwise, as
\[
 q^\nu:X^\nu\longrightarrow Y^\nu,\qquad
 (\nu^*L)^a\simeq(q^\nu)^*A,
\]
where $Y^\nu$ is a finite disjoint union of normal projective varieties
and $A$ is ample. Let $D\subset X$ be the conductor subscheme and
$D^\nu=X^\nu\times_XD$. Assume:
\begin{enumerate}
\item The set $E=q^\nu(|D^\nu|)$ is finite.
\item For every morphism $g:C\to X$ from a connected proper nodal curve,
if the pullback of $L$ to every normalized irreducible component of $C$
is torsion, then $g^*L$ is torsion on $C$.
\end{enumerate}
Then $L$ is semiample.
\end{proposition}

\begin{proof}
We give $E$ its reduced structure. On each connected component
$R_\alpha$ of $(D^\nu\times_DD^\nu)_{\red}$, the two maps to
$Y^\nu$ have constant values $e,e'\in E$. The conductor
comparison is therefore a constant isomorphism
\[
 \lambda_\alpha:A_e\xrightarrow{\sim}A_{e'},
\]
since a proper connected reduced space has only constant units.
These data form a finite labelled graph on $E$.

Fix an oriented graph cycle $\gamma$ based at $e\in E$. Let
$\lambda_\gamma\in\operatorname{Aut}(A_e)=\C^*$ be the composite
of its edge comparisons, using inverses when an edge is traversed
in the opposite direction. We choose conductor points representing
these edges and join consecutive lifts by curve chains in the
connected $q^\nu$-fibres, using a projective resolution of $X^\nu$.
We normalize the chains and glue markings with the same image,
inserting constant rational components to separate repeated marks.
This gives a connected proper nodal curve $g:C\to X$ realizing
$\gamma$.

On each normalized component of $C$, the pullback of $L^a$ is
trivial, since $(\nu^*L)^a\simeq(q^\nu)^*A$. Thus the pullback
of $L$ is torsion on each such component. The componentwise
trivializations identify $\lambda_\gamma$ with the holonomy of
$g^*(L^a)$ around the corresponding closed chain in $C$.
Assumption~(2) makes $g^*L$ torsion, so $g^*(L^a)$ is torsion
and $\lambda_\gamma$ is a root of unity. The graph
has finitely generated first homology, so there is one $n_0>0$ with
\[
 \lambda_\gamma^{n_0}=1\qquad\text{for every cycle }\gamma.
\]

We now construct sections generating $L^{am}$.
For a sufficiently large multiple $m$ of $n_0$, Serre vanishing gives
\[
 H^0(Y^\nu,A^m)\twoheadrightarrow\bigoplus_{e\in E}A_e^m,
 \qquad\cI_EA^m\text{ globally generated}.
\]
For each connected component of the graph, trivial holonomy gives
a tuple of nonzero vectors satisfying the edge comparisons
$\lambda_\alpha^m$. We extend it by zero at the other vertices and
lift it to a global section. Together with sections generating
$\cI_EA^m$, these lifts generate $A^m$ and satisfy every comparison.

Their pullbacks agree at every pair of points over the same point
of $X$. We take an \'etale affine chart $\Spec R\to X$ on which
$L^{am}$ is trivial, and write $R^\nu$ for the normalization of $R$.
Each pulled-back section then gives an element $s\in R^\nu$ and
\[
 R\subset R[s]\subset R^\nu.
\]
This finite extension is subintegral: in characteristic zero, failure
of geometric injectivity would occur at a closed point. Seminormality
therefore gives $R[s]=R$. The descended sections generate $L^{am}$
because $\nu$ is surjective. If $E=\varnothing$, the same descent
applies to any generating power of $A$.
\end{proof}

\begin{remark}
In assumption~(2), it suffices to test the finitely many nodal curves
realizing generators of the conductor graph. Seminormality is essential:
pointwise matching does not recover nilpotent conductor conditions.
\end{remark}

\section{The canonical bundle formula}\label{sec:cbf}

We first prove b-semiampleness over algebraic bases. We pass to a
locally stable family and identify its relative
pluricanonical line with a power of a Hodge line. To apply Hodge
semiampleness, we verify numerical integrability using moduli
b-abundance on smooth test spaces, and torsion combinatorial
monodromy using finite pluricanonical representations on nodal
test curves.
We then descend semiampleness to the original moduli b-divisor.
In Section~\ref{sec:analytic-cbf}, we extend the result to analytic
bases by constructing algebraic parameter families and comparing
their moduli lines with that of the original fibration.

\subsection{Locally stable preparation}\label{sec:stable-preparation}

We use lc-trivial and klt-trivial fibrations in the sense of
\cite[Definitions 3.1--3.2]{FG14}. Thus an lc-trivial fibration
$f:(V,\Delta)\to B$ is a proper surjective morphism of normal varieties
with connected fibres. Here $\Delta$ is a $\Q$-divisor on $V$,
the pair $(V,\Delta)$ is sub-lc over the generic point of $B$, and
we require
\[
 \operatorname{rank}f_*\cO_V\bigl(\lceil\mathbf A^*(V,\Delta)\rceil\bigr)=1,
 \qquad K_V+\Delta\sim_{\Q}f^*D
\]
for a $\Q$-Cartier $\Q$-divisor $D$ on $B$. Here the discrepancy
b-sheaf is computed on a log resolution and $\mathbf A^*$ is obtained
from the discrepancy b-divisor by deleting its coefficients equal
to $-1$. For a klt-trivial fibration, sub-lc is replaced by sub-klt
and $\mathbf A^*$ by $\mathbf A$. Effectiveness of the induced
geometric generic boundary is not part of these definitions; we
impose it explicitly where needed.

A \emph{log smooth Ambro model} is a proper birational morphism
$B'\to B$, with $B'$ smooth, on which the moduli b-divisor descends and is nef,
the support of $B_{B'}$ together with the exceptional divisors is SNC,
and, on every higher model $\nu:B''\to B'$,
\[
 M_{B''}=\nu^*M_{B'},\qquad
 K_{B''}+B_{B''}=\nu^*(K_{B'}+B_{B'}).
\]
Such models exist by \cite[Theorem 3.6]{FG14} and resolution.
Thus it is $\mathbf K+\mathbf B$, rather than the discriminant
b-divisor $\mathbf B$ alone, that descends by Cartier pullback.

For a locally stable family $f:(Y,\Delta)\to W$, the morphism $f$
is flat, $\Delta\ge0$ is a relative $\Q$-divisor, every fibre pair
is slc, and $K_{Y/W}+\Delta$ is $\Q$-Cartier. We use the conventions of
\cite[Definition--Theorems 2.3 and 4.7]{KolFamilies}.
Following the preparation strategy of \cite[Proposition 6.16]{BFMT},
we construct such a model by weak semistable reduction and relative
minimal model theory. The original boundary may be negative vertically.

\begin{lemma}[Locally stable log Calabi--Yau models]
\label{cbf:locally-stable-preparation}
Let $h:(V,\Delta)\to T$ be a projective lc-trivial fibration over a
normal projective base. Assume that $\Delta_{\bar\eta}\ge0$, where
$\bar\eta$ is a geometric generic point of $T$.
There are a smooth projective alteration $q:W\to T$ and a projective
locally stable family $g:(Y,\Gamma)\to W$ with connected
fibres such that
\[
 Y\text{ normal},\qquad \Gamma\ge0,\qquad
 K_Y+\Gamma\sim_{\Q,W}0.
\]
Its geometric generic pair is crepant birational to the original one.
The moduli b-divisor is the induced pullback of that of $h$.
Any prescribed finite base extension may be made before the construction.
\end{lemma}
\begin{proof}
\emph{Step 1. A locally stable family with relatively semiample log canonical divisor.}
After the prescribed finite extension and a further finite extension,
we choose a crepant $\Q$-factorial dlt model $(F,\Delta_F)$ of the
geometric generic pair. We take a log resolution $\mu:\widetilde F\to F$
whose exceptional discrepancies are $>-1$, and set
\[
 \widetilde\Delta=\mu_*^{-1}\Delta_F+\sum_E E,
 \qquad G=\sum_E(1+a(E,F,\Delta_F))E\ge0.
\]
Then
\[
 K_{\widetilde F}+\widetilde\Delta
       =\mu^*(K_F+\Delta_F)+G.
\]
The divisor $G$ is exceptional. Since $K_F+\Delta_F\sim_{\Q}0$,
$(F,\Delta_F)$ is a good minimal model of
$(\widetilde F,\widetilde\Delta)$.

We spread this resolution and equality over a nonempty open of the
altered base. We apply weak semistable reduction \cite{AK00}, with the
horizontal boundary, as in \cite[Theorem 4.59]{KolFamilies}.
We use the locally stable intermediate family in that proof.
It extends our log smooth family; we do not take the canonical model
used there for generically ample log canonical divisors.
Therefore its generic pair still has the good minimal model
$(F,\Delta_F)$. We resolve the altered base and the graph of the map
to $T$ to obtain a smooth projective $W$.
Reduced base change preserves local stability
\cite[Definition--Theorem 4.7]{KolFamilies}, giving
$g_1:(Y_1,\Gamma_1)\to W$ projective and locally stable.

We need $Y_1$ normal before taking a dlt modification.
Flatness, $S_2$ fibres and smoothness of $W$ give $S_2$.
Horizontal codimension-one points are regular because the generic
fibre is normal. At a vertical codimension-one point the reduced
fibre makes a base uniformizer generate the maximal ideal.
Hence $R_1$ also holds, and Serre's criterion gives normality.

We now take a crepant $\Q$-factorial dlt modification
$\mu:(Y_2,\Gamma_2)\to(Y_1,\Gamma_1)$ which is an isomorphism
over the log smooth open, and write $g_2=g_1\mu$.
For every SNC divisor $D\subset W$, we have
\[
 K_{Y_2}+\Gamma_2+g_2^*D
       =\mu^*(K_{Y_1}+\Gamma_1+g_1^*D).
\]
The right side is the crepant pullback of an lc pair, so $g_2$ is
locally stable by \cite[Corollary 4.55]{KolFamilies}.
Every lc stratum dominates $W$ by Corollary 4.56 there.
The generic good minimal model and its pluricanonical trivialization
spread to an open $W^\circ$; every lc stratum meets this open.
We may therefore apply \cite[Theorem 1.1]{HX13} to obtain a good
minimal model $g:(Y,\Gamma)\to W$. It remains locally stable by
\cite[Corollary 4.57(1)]{KolFamilies}, since it is a weak canonical
model of the locally stable dlt family.

\emph{Step 2. Relative $\Q$-linear triviality.}
For divisible $r$, let the relative semiample contraction be
\[
 Y\xrightarrow{p}Z\xrightarrow{\rho}W,
 \quad p_*\cO_Y=\cO_Z,
 \quad\mathcal N=\cO_Y(r(K_Y+\Gamma))\simeq p^*A,
\]
with $A$ $\rho$-ample. The generic log Kodaira dimension is zero,
so $\rho$ is birational and $\mathcal N_\eta\simeq\cO$.
Let $H$ be $g$-ample and suppose $n=\dim(Y/W)>0$.
Flatness gives $\mathcal N_w\cdot H_w^{n-1}=0$ for all $w\in W$.
Each component contributes nonnegatively. A component with
positive-dimensional image under $p$ would contribute positively.
Hence $\rho$ has finite fibres; this is immediate for $n=0$ as well.
Normality now gives $Z=W$, and thus
\[
 K_Y+\Gamma\sim_{\Q,W}0.
\]
Stein factorization and reduced fibres give
$H^0(Y_w,\cO_{Y_w})=\C$. Cohomology and base change show that the
fibres are connected. The generic pair is lc with
effective boundary. On a log resolution, $\lceil\mathbf A^*\rceil$
is effective and exceptional and its sheaf pushes forward to $\cO$.
This establishes the rank-one condition for $g$.

\emph{Step 3. Pullback of the moduli b-divisor.}
For a prime $D\subset W$, local stability gives
$(Y,\Gamma+g^*D)$ lc over its generic point. Since the fibre is
reduced, the threshold cannot exceed one. Therefore
\[
 \operatorname{lct}_{\eta_D}(Y,\Gamma;g^*D)=1,
 \qquad B_W(g)=0.
\]
Let $p:Z\to W$ be a common resolution of the induced original
family and $g$, with crepant boundaries $\Delta_Z,\Gamma_Z$.
They agree generically, so $E=\Delta_Z-\Gamma_Z$ is vertical.
Relative $\Q$-triviality writes it as
\[
 E=p^*D_0+s^{-1}\operatorname{div}(v),
 \qquad v\in\C(Z)^*,\quad D_0\text{ a }\Q\text{-divisor on }W.
\]
The function $v$ has no zero or pole on the proper normal
geometrically connected generic fibre, and hence lies in $\C(W)^*$.
Hence $E$ itself is pulled back from $W$. Adding a divisor from the
base changes the discriminant and total divisor in the canonical
bundle formula equally. It does not change the moduli part;
the same argument holds on every higher model.

It remains to compare after finite base change.
For $\tau:W'\to W$, the threshold formula
\cite[Lemma 5.1]{Amb04} applies to lc generic pairs. At a divisorial
valuation with ramification index $e$, it gives
\[
 t'=et,\qquad b'=1-e+eb,
 \qquad K_{W'}+B_{W'}=\tau^*(K_W+B_W).
\]
The canonical ramification term cancels the discriminant correction,
so the moduli part pulls back. We factor an alteration through a finite
map and a birational map to obtain the required equality of b-divisors.
\end{proof}

\begin{lemma}\label{cbf:relative-line}
Let $g:(Y,\Gamma)\to W$ be a projective surjective locally stable
morphism, where $Y$ is normal, $W$ is a smooth connected variety,
and $\Gamma\ge0$ is a $\Q$-divisor on $Y$. Assume that $g$ has
connected fibres and $K_Y+\Gamma\sim_{\Q,W}0$.
Then $g$ is lc-trivial. We write $\mathbf M_g$ for its moduli
b-divisor and $M_W$ for its trace on $W$.
For sufficiently divisible $r>0$, there are isomorphisms
\[
 \mathcal N:=\omega_{Y/W}^{[r]}(r\Gamma)\simeq g^*P,\qquad
 P=g_*\mathcal N\simeq\cO_W(rM_W).
\]
For every morphism $u:T\to W$, where $T$ is smooth or a nodal curve,
we put $Y_T=Y\times_WT$. These isomorphisms induce
\[
 \omega_{Y_T/T}^{[r]}(r\Gamma_T)\simeq g_T^*u^*P,\qquad
 (g_T)_*\omega_{Y_T/T}^{[r]}(r\Gamma_T)\simeq u^*P.
\]
These identifications commute with normalization and fibre restriction.
Moreover, $\mathbf M_g=\overline{M_W}$.
\end{lemma}
\begin{proof}
The generic pair $(Y_\eta,\Gamma_\eta)$ is lc, with $\Gamma_\eta\ge0$.
On a log resolution, $\lceil\mathbf A^*\rceil$ is effective and
exceptional, so its sheaf pushes forward to $\cO$. Geometric connectedness therefore gives the
rank-one condition. Local stability and reduced fibres give
divisorial thresholds one, hence $B_W=0$. Clearing denominators
in the canonical bundle formula gives an isomorphism
$\mathcal N\simeq g^*\cO_W(rM_W)$; the projection formula and
$g_*\cO_Y=\cO_W$ identify its pushforward and evaluation map.

To check compatibility with base change, we fix $r$ so that
$r\Gamma$ is integral, $\mathcal N$ is
invertible, and the displayed comparison holds, and use this same
$r$ for every base change. Local stability survives reduced
base change.
The divisorial pluricanonical sheaves have the fibre restrictions
and the $S_2$ property of
\cite[Corollary 4.33 and Proposition 2.79(3)]{KolFamilies}.
Let $\rho:Y_T\to Y$ be the projection, and let $U_T\subset Y_T$
be the pullback of the relatively Gorenstein locus where $r\Gamma$
is Cartier. Relative duality and Cartier pullback give
\[
 (\rho^*\mathcal N)|_{U_T}
   \simeq\omega_{U_T/T}^{\otimes r}(r\Gamma_T|_{U_T}).
\]
This locus contains the codimension-one
points of every fibre: an slc fibre is nodal in codimension one and
its boundary avoids the nodes. The complement has codimension at
least two in each fibre. Since $Y_T$ is flat with $S_2$ fibres over
the $S_2$ base $T$, it satisfies $S_2$. We therefore extend the
isomorphism uniquely to all of $Y_T$ as an isomorphism of divisorial
sheaves. The comparisons already commute with composition,
normalization with the conductor, and fibre restriction in
codimension one. Uniqueness of divisorial extension preserves these
compatibilities. Proper flatness and
geometrically reduced connected fibres give
$(g_T)_*\cO_{Y_T}=\cO_T$ by cohomology and base change; projection
formula proves the pushforward assertion.

Finally, we prove $\mathbf M_g=\overline{M_W}$.
We take a smooth birational base change
$\mu:W'\to W$ and the relative locally stable model $(Y',\Gamma')$.
The same relative-line comparison, with compatible canonical
representatives, identifies the absolute crepant boundary with
$\Gamma'-g'^*K_{W'/W}$. The locally stable family has discriminant
zero, whereas its absolute crepant model has
\[
 B_{W'}=-K_{W'/W},\qquad M_{W'}=\mu^*M_W.
\]
Every divisorial valuation occurs on such a model, proving
$\mathbf M_g=\overline{M_W}$. Hence only the trace $B_W$ vanishes;
the discriminant b-divisor need not vanish.
\end{proof}

The normalization assertion in the following lemma also appears in
\cite[Lemma 6.13]{BFMT}. Here a \emph{klt source} is the pair obtained
by adjunction to a minimal lc stratum on a crepant dlt model; its map
to the base is factored through its Stein factorization.

\begin{lemma}\label{sg:smooth-branch}
Let $f:(Y,\Delta)\to W$ be a projective surjective locally stable
family over a smooth connected base, with $\Delta\ge0$ a
$\Q$-divisor on $Y$.
Normalization with the conductor, crepant dlt modification of each
normal component, and adjunction to any lc stratum give locally
stable families dominating $W$.

Suppose in addition that, for a sufficiently divisible $r>0$,
$\omega_{Y/W}^{[r]}(r\Delta)\simeq f^*P$ for a line $P$ on $W$.
There are a finite \'etale surjection $\sigma:T\to W$ and a locally
stable klt-trivial fibration $g:(S,B_S)\to T$, where $S$ is normal,
$B_S\ge0$ is a $\Q$-divisor, and $g$ has connected fibres,
such that, after a common
replacement $(r,P)\mapsto(ar,P^a)$,
\[
 \omega_{S/T}^{[r]}(rB_S)\simeq g^*\sigma^*P,
 \qquad B_T(g)=0,\qquad
 \mathbf M_g=\overline{M_T(g)},\quad
 \cO_T(rM_T(g))\simeq\sigma^*P.
\]
The line comparison is induced by pluricanonical residue and respects
multiplication.
\end{lemma}
\begin{proof}
For every local SNC divisor $D\subset W$, the pair
$(Y,\Delta+f^*D)$ is slc by \cite[Corollary 4.55]{KolFamilies}.
Its normalization with the conductor is lc, so the same criterion
proves local stability of each normalized component.
For a crepant dlt modification $\mu:(Y',\Delta')\to(Y^\nu,\Delta^\nu)$,
the identity
\[
 K_{Y'}+\Delta'+(f\mu)^*D
 =\mu^*(K_{Y^\nu}+\Delta^\nu+f^*D)
\]
gives local stability again. Every lc stratum dominates $W$ by
\cite[Corollary 4.56]{KolFamilies}. For a coefficient-one component
$S$, adjunction commutes with adding the pullback Cartier divisor:
\[
 \operatorname{Diff}_{S}(\Delta'-S+(f\mu)^*D)
 =\operatorname{Diff}_{S}(\Delta'-S)+(f\mu|_S)^*D.
\]
We prove the first assertion by repeating adjunction and applying
the same criterion; compare \cite[Lemma 2.11]{KolFamilies}.

For the second assertion, we take a minimal lc stratum $S$, or a
whole normal component when its induced pair is klt.
The pair on $S$ obtained by adjunction
is klt over the geometric generic point of $W$. Its Stein
factorization is
\[
 S\xrightarrow{g}T\xrightarrow{\sigma}W,\qquad
 g_*\cO_S=\cO_T.
\]
Since $S\to W$ is proper and flat with geometrically reduced fibres,
$\sigma$ is finite \'etale \cite[Tag 0BUN]{Stacks}.
Local stability is \'etale-local, so $g:S\to T$ is locally stable
with connected fibres. After replacing $(r,P)$ by
$(ar,P^a)$ for a common multiple $ar$, normalization, crepant pullback and ordered residues identify
$\omega_{S/T}^{[r]}(rB_S)$ with $g^*\sigma^*P$.
We obtain its pushforward by the projection formula;
these identifications preserve products.
Lemma~\ref{cbf:relative-line} gives the rank-one condition,
$B_T=0$, $\cO_T(rM_T)\simeq\sigma^*P$, and
$\mathbf M_g=\overline{M_T}$. Since the geometric generic pair
is klt, $g$ is klt-trivial.
\end{proof}

\begin{lemma}[Cyclic Hodge realization]\label{cbf:cyclic-realization}
Let $g:(Y,\Gamma)\to T$ be a projective klt-trivial fibration over
a projective base. Assume that $\Gamma_{\bar\eta}\ge0$, where
$\bar\eta$ is a geometric generic point of $T$. After a smooth
projective alteration and logarithmic preparation, its moduli
b-divisor descends to $T$ and there is an integral polarizable pure
CY variation on $T\setminus E$, where $E$ is SNC, with unipotent
monodromy and Schmid highest line $\overline{\mathcal H}$ satisfying
\[
 \overline{\mathcal H}^{\otimes r}\simeq\cO_T(rM_T)
 \qquad(r>0\text{ sufficiently divisible}).
\]
The comparison is preserved under further generically finite base
change between smooth log pairs with unipotent monodromy.
\end{lemma}
\begin{proof}
We take the minimal generic index $d$ and a representative
\[
 K_Y+\Gamma+\tfrac1d\operatorname{div}(\varphi)
       =g^*(K_T+B_T+M_T),\qquad \psi^d=\varphi.
\]
We take an equivariant resolution of the normalization in the cyclic
algebra $\C(Y)[z]/(z^d-\varphi)$, with deck action $z\mapsto\zeta z$.
Under further base change on the smooth open, we retain this entire
algebra and all normalization components if it splits. The character
$\chi(\zeta)=\zeta$ is unchanged.
By the logarithmic preparation in \cite[Proposition 5.4]{Amb04},
we may assume that the resulting cover $\bar\pi:\widetilde Y\to T$
is semistable in codimension one, with SNC degeneration. We write
$\pi=\bar\pi|_{T\setminus E}$, $n=\dim(Y/T)$, and
\[
 \mathbb W_{\Z}=R^n\pi_*\Z/\mathrm{tors},\qquad
 \mathcal H=F^n\mathbb W_\chi,\qquad\operatorname{rk}\mathcal H=1.
\]
On this model, $M_T$ is integral by \cite[Lemma 5.2(5)]{Amb04}.
For the induced log resolution $g':(Y',\Gamma')\to T$, parts (4)--(5)
of that lemma identify the eigensheaf:
\[
 (\bar\pi_*\omega_{\widetilde Y/T})_\chi
 =g'_*\cO_{Y'}\!\left(
   \left\lceil-\Gamma'+g'^*B_T+g'^*M_T\right\rceil\right)\psi
 \simeq\cO_T(M_T)\psi.
\]
The eigensheaf on the left is the Schmid extension of $\mathcal H$
\cite[proof of Theorem 3.3]{Amb05}. We obtain the displayed
tensor-power comparison by clearing the chosen principal representatives.
Its base-change compatibility follows from the natural eigensheaf
comparison in \cite[proof of Proposition 5.5]{Amb04}.

We make $\mathcal H$ the highest line of an integral variation,
using the shift of \cite[Remark 6.26]{BFMT}. A relatively ample
class and the Lefschetz decomposition give a polarization $Q$ of
$\mathbb W_\Q$. After clearing denominators, it is integral on
$\mathbb W_\Z$. We replace $Q$ by the sum of its $\mu_d$-translates;
the resulting integral polarization is invariant under the deck group.
For nonreal $\chi$, we retain the entire lattice and put
\[
\begin{gathered}
 \mathbb W'_{\Z}=\mathbb W_{\Z},\qquad
 \mathbb W'_{\C}=\mathbb W_\chi(-2,2)
  \oplus\!\!\bigoplus_{\eta\ne\chi,\bar\chi}\mathbb W_\eta
  \oplus\mathbb W_{\bar\chi}(2,-2),\\
 (-2,2):(p,q)\longmapsto(p+2,q-2),\qquad
 F^{n+3}\mathbb W'=0,\quad F^{n+2}\mathbb W'=\mathcal H.
\end{gathered}
\]
The opposite shifts on conjugate summands preserve the real
structure and weight. Since $Q$ pairs only inverse characters and
$i^{p-q}$ is unchanged, $Q$ still polarizes the shifted structure.
Constant shifts on flat eigensummands preserve Griffiths transversality.
For real $\chi=\pm1$, we use instead
\[
 \mathbb W'_{\Q}=\mathbb W_{\chi,\Q},\qquad
 \mathbb W'_{\Z}=\mathbb W_{\Z}\cap\mathbb W'_{\Q},\qquad
 F^{n+1}\mathbb W'=0,\quad F^n\mathbb W'=\mathcal H.
\]
This lattice is saturated and the restricted polarization is
nondegenerate. In both cases the connection and highest subline are
unchanged, so their unipotent Schmid extension remains
$\overline{\mathcal H}$.
\end{proof}

\subsection{Numerical integrability from abundance}

For a smooth projective $T$ and a morphism $T\to W$, we apply
Ambro's abundance theorem to the family over $T$.
Lemma~\ref{cbf:relative-line} identifies its pluricanonical line
even when $T$ maps into the boundary. We then compare the
fibres of the abundant divisor with the fibres of the period map.

\begin{lemma}\label{cbf:test-abundance}
Let $g:(Y,\Gamma)\to W$ and $P$ be as in
Lemma~\ref{cbf:relative-line}. For every morphism $u:T\to W$ with
$T$ smooth connected projective, there are a finite \'etale
surjection $\tau:T_1\to T$, a smooth projective birational model
$p:T_2\to T_1$, a contraction $\rho:T_2\to Q$ to a normal
projective variety, and an integer $\ell>0$ such that
\[
 (p^*\tau^*u^*P)^\ell\simeq\rho^*H,
 \qquad H\in\operatorname{Pic}(Q)\text{ nef and big}.
\]
\end{lemma}
\begin{proof}
We apply Lemmas~\ref{cbf:relative-line} and~\ref{sg:smooth-branch} to
the family pulled back to $T$. After a finite \'etale cover
$\tau:T_1\to T$, we obtain a locally stable klt source. For a
suitable $a>0$, its relative pluricanonical line at index $ar$ is
$(\tau^*u^*P)^a$. If $\mathbf M$ denotes its moduli b-divisor, then
\[
 \mathbf M=\overline{M_{T_1}},\qquad
 \cO_{T_1}(arM_{T_1})\simeq(\tau^*u^*P)^a.
\]
By Ambro's theorem \cite[Theorem 3.3 and Definition 3.2]{Amb05},
there are $p:T_2\to T_1$ birational and $\rho:T_2\to Q$ a
contraction such that
\[
 p^*M_{T_1}\sim_{\Q}\rho^*H_0,
\]
where $Q$ is normal proper and $H_0$ is a nef and big $\Q$-divisor.
We can make $Q$ projective by taking a projective modification,
and make $T_2$ smooth projective by resolving the graph.
The pullback of $H_0$ is still nef and big.
For a sufficiently divisible $\ell>0$ with $a\mid\ell$, put
$H=\cO_Q(r\ell H_0)$. Then
\[
 (p^*\tau^*u^*P)^\ell
 \simeq\cO_{T_2}(r\ell p^*M_{T_1})\simeq\rho^*H,
\]
as required.
\end{proof}

\begin{lemma}
\label{cbf:abundant-period}
Let $T$ be a smooth projective variety and $E\subset T$ a reduced SNC
divisor. Let $\mathbb V$ be a polarizable integral pure CY variation
on $T\setminus E$, minimal
with respect to its highest line. Let $A$ be its unipotent Schmid
highest line. Suppose a smooth projective alteration $a:T'\to T$
and a surjective morphism $\rho:T'\to Q$ to a normal projective
variety satisfy, for an integer $\ell>0$,
\[
 (a^*A)^\ell\simeq\rho^*H,\qquad
 H\in\operatorname{Pic}(Q)\text{ nef and big}.
\]
If $\Phi_{\mathbb V}$ is generically immersive, then $A$ is big.
\end{lemma}
\begin{proof}
Suppose that $A$ is not big. Since $H$ is nef and big and
$(a^*A)^\ell\simeq\rho^*H$, we have $\dim Q<\dim T$.
We take $x\in T'$ such that $a$ is \'etale at $x$, $\rho$ is smooth
at $x$, $d(\Phi_{\mathbb V}\circ a)_x$ is injective, and
$\mathbb V_{a(x)}$ is pointwise minimal. Such a point exists by
Lemma~\ref{prep:generic-minimal}.

The fibre of $\rho$ through $x$ has positive dimension.
We take an integral complete-intersection curve $C$ in that fibre,
smooth at $x$, and replace it by its normalization. We have
\[
 \ell(a^*A\cdot C)=\rho^*H\cdot C=0.
\]
By Schmid functoriality and Lemma~\ref{prep:line}, the highest
Hodge line on the open curve is a flat subsystem.
Its character has finite image
\cite[Corollaire 4.2.8(iii)(b)]{DeligneHII}.
After a finite \'etale cover of the open curve, the fixed-part
theorem therefore gives a constant rational subvariation
containing this line.

This subvariation must be the whole restricted variation.
Indeed, its orthogonal complement has zero highest Hodge piece,
so its fibre at a point over $x$ is zero by pointwise minimality.
The complement is a local system, and hence is zero everywhere.
Hence the period map is constant along $C$, and
\[
 0\ne T_xC\subset\ker d(\Phi_{\mathbb V}\circ a)_x.
\]
This contradicts the choice of $x$, proving that $A$ is big.
\end{proof}

\subsection{Pluricanonical gluing on nodal curves}
\label{sec:nodal-holonomy}

As in \cite[proof of Theorem 7.2]{BFMT}, we pass to klt sources
and use pluricanonical representations. Pluricanonical sheaves are
divisorial extensions from the relatively Gorenstein locus, with
comparison maps fixed before restriction to a curve. We compare
their residues at each node: componentwise torsion alone is
insufficient, as Example~\ref{ex:nodal-descent} shows.

We use the compatibility of residue maps in
\cite[Proposition 14, Remark 15 and Theorem 17]{KolSources}.
For a $\mathbb P^1$-link the two residues of $(dt/t)\wedge\eta$
are $\eta$ and $-\eta$; an even pluricanonical power removes the sign.

\begin{lemma}\label{nh:etale-preparation}
Let $C$ be a connected projective nodal curve, with normalization
$\coprod_i C_i$. Finitely many finite \'etale covers of the $C_i$
split after a finite \'etale surjection $\pi:C'\to C$. If $P$ is a
line bundle on $C$ and $\pi^*P$ is torsion, then $P$ is torsion.
\end{lemma}
\begin{proof}
On each $C_i$, we choose a common Galois refinement of the prescribed
covers. By taking disjoint copies, we arrange that all degrees equal
one integer $d$. We then glue the two $d$-element fibres over each
node by a bijection. The resulting cover is finite \'etale: at a node it is,
\'etale locally, a disjoint union of copies of the original node.
Finally,
\[
 (\pi^*P)^N\simeq\cO_{C'}
 \quad\Longrightarrow\quad
 \operatorname{Nm}_{C'/C}((\pi^*P)^N)=P^{dN}\simeq\cO_C.
\]
\end{proof}

\begin{lemma}[Residues on an slc fibre]\label{nh:fibre-residues}
Let $(F,\Delta_F)$ be a connected projective slc pair, and fix
\[
 A\otimes\cO_F\xrightarrow{\sim}
 \omega_F^{[r]}(r\Delta_F),\qquad \dim A=1,
\]
for a sufficiently divisible even $r>0$. Minimal lc centres on crepant
dlt models of the normalization components form a connected graph of
crepant birational comparisons. A component whose induced pair
is klt is itself included as a vertex. For a vertex $(Z_v,B_v)$,
residue gives an isomorphism
\[
 R_v:A\xrightarrow{\sim}H^0(Z_v,r(K_{Z_v}+B_v)).
\]
Each edge $\phi:Z_w\dashrightarrow Z_v$ satisfies
$\phi^*R_v=R_w$. These comparisons also commute with successive
adjunction from an ambient lc pair.
\end{lemma}
\begin{proof}
We normalize $F$, include the conductor in the boundary, and take
crepant dlt models $(F_j',\Delta_j')$. The given section pulls back
to a nowhere-zero section of $r(K_{F_j'}+\Delta_j')$. Within one
component, minimal lc centres are $\mathbb P^1$-linked by
\cite[Theorem 10]{KolSources}. Proposition 14 of that paper
identifies their residue maps by crepant birational pullback.
The even index removes the sign of a $\mathbb P^1$-link.

Since $F$ is connected and $S_2$, its normalization components are
connected through the codimension-one conductor.
The normalization of the conductor carries an involution $\iota$,
and the residues of the given pluricanonical section satisfy
\[
 \operatorname{Res}(\sigma_j)
   =\iota^*\operatorname{Res}(\sigma_k)
 \qquad(r\text{ even}).
\]
Successive adjunction along the conductor, followed by
\cite[Theorem 10]{KolSources} within the component, connects the
minimal sources on the two sides. Theorem 17 of that paper makes
these maps compatible with the original residues. All vertices are
klt log Calabi--Yau pairs; their displayed spaces of sections are
one-dimensional. This establishes the assertions.
\end{proof}

\begin{theorem}[Pluricanonical torsion on a nodal base]
\label{nh:prepared-holonomy}
Let $f:(Y,\Delta)\to C$ be a projective locally stable family with
connected fibres over a connected projective nodal
curve. Suppose that, for a sufficiently divisible $r>0$, there is
an isomorphism
\begin{equation}\label{eq:nh-relative-line}
 \omega_{Y/C}^{[r]}(r\Delta)\simeq f^*P,
 \qquad P\in\operatorname{Pic}(C),
\end{equation}
compatible with restriction to the fibres and to the normalization
branches of $C$. If $P|_{C_i}$ is torsion for every normalized
component $C_i$, then $P$ is torsion. Its canonical componentwise
flat connections have finite monodromy with the given node gluing.
\end{theorem}
\begin{proof}
We compare transport around a graph cycle with the pluricanonical
representation of one fixed klt source.

\emph{Step 1.}
By Lemma~\ref{nh:etale-preparation}, we may make finite \'etale
base changes on $C$ and descend by the norm at the end. We include
the cyclic covers trivializing $P|_{C_i}$ and replace $(r,P)$ by
$(ar,P^a)$, so that $r$ is even and $P|_{C_i}\simeq\cO_{C_i}$.
Lemma~\ref{sg:smooth-branch}, followed by another cover as in
Lemma~\ref{nh:etale-preparation}, gives locally stable klt sources
$h_i:(S_i,B_i)\to C_i$ with connected fibres, discriminant zero
and relative line $h_i^*P|_{C_i}$. Their moduli b-divisors are
$\Q$-linearly trivial.

By \cite[Theorem 4.7]{Amb05}, each becomes a product over a dense
open after a finite cover \'etale in codimension one. This cover
is \'etale since $C_i$ is smooth. Proposition 4.4 of that paper
extends the product away from the singular locus of the base,
the discriminant support and the image of the negative vertical
boundary. Here all three are empty. Hence
\[
 (S_i,B_i)\simeq(Z_i,B_{Z_i})\times C_i,
 \qquad K_{Z_i}+B_{Z_i}\sim_{\Q}0.
\]
Residues give frames
\[
 \beta_i:H^0(Z_i,r(K_{Z_i}+B_{Z_i}))\otimes\cO_{C_i}
       \xrightarrow{\sim}P|_{C_i}.
\]
Their transport is the canonical transport of the trivial line,
since $C_i$ is complete.

\emph{Step 2.}
At a node $p$, the two branches have the same slc fibre
$(F,\Delta_F)$ and the same identification
$P_p\otimes\cO_F\simeq\omega_F^{[r]}(r\Delta_F)$.
By Lemma~\ref{sg:smooth-branch}, conductor strata of a normalized
branch are locally stable and horizontal, so contain no fibre
component. Normalization is therefore birational on fibre
components, and adjunction keeps the induced conductor boundary.

Let $g_i:(Y_i,\Theta_i)\to C_i$ be the normal dlt branch model
from which $(S_i,B_i)$ was obtained, and let $p_i$ lie over $p$.
We take a crepant dlt model of $(Y_i,\Theta_i+g_i^*p_i)$ dominating
it. For a coefficient-one horizontal component $D$, we have
\[
 \operatorname{Diff}_{D^\nu}(\Theta_i-D+g_i^*p_i)
 =\operatorname{Diff}_{D^\nu}(\Theta_i-D)
       +(g_i|_{D^\nu})^*p_i.
\]
We factor the map from this model to $C_i$ through its Stein
factorization and work over the point met by the chosen horizontal
source. The fibre is now connected, its residue field is $\C$,
and the log canonical divisor remains relatively $\Q$-trivial.
We compare two orders of adjunction in this same dlt pair.
Following the horizontal adjunction chain first gives the plt
product $(Z_i,B_{Z_i})\times(C_i,p_i)$ up to crepant comparison;
adjunction to its fibre gives the klt pair $(Z_i,B_{Z_i})$.
Thus the terminal centre is minimal: by transitivity of adjunction,
a smaller lc centre would induce an lc centre of this klt pair.

In the same connected fibre, choose a vertical component dominating
a component of the normalized slc fibre. Adjunction to it is
crepant birational to that fibre component with its conductor
boundary. By \cite[Theorem~10]{KolSources}, it contains a minimal
lc centre linked to the horizontal one. Its induced pair gives
a vertex $(Z_{p,v},B_{p,v})$ in Lemma~\ref{nh:fibre-residues}.
We compare residues on the common sufficiently divisible even
power $r$. Proposition~14 of \cite{KolSources} identifies them
along the links, and Theorem~17 makes the comparison compatible
with these successive adjunctions and crepant models.
We use the model only for adjunction, so do not need local stability
of its vertical strata.

With $H_r(Z,B)=H^0(Z,r(K_Z+B))$, the two chains give
\[
 \begin{array}{ccc}
 P_p&\xrightarrow{\ R_{\mathrm{hor}}\ }&H_r(Z_i,B_{Z_i})\\
 \Vert&&\Big\downarrow\psi^*\\
 P_p&\xrightarrow{\ R_{\mathrm{vert}}\ }&H_r(Z_{p,v},B_{p,v}),
 \end{array}
\]
where $\psi:(Z_{p,v},B_{p,v})\dashrightarrow(Z_i,B_{Z_i})$ is
crepant birational and $(Z_{p,v},B_{p,v})$ is a minimal source of
a normalized fibre component. After the crepant comparison,
$R_{\mathrm{vert}}$ is the map $R_v$ of
Lemma~\ref{nh:fibre-residues} on the section induced by $P_p$.

To check the identity on the left, let $N_i$ be the preimage of
the nodes on $C_i$. Then
\[
 \nu_C^*\omega_C|_{C_i}=\omega_{C_i}(N_i).
\]
For a branch parameter $t$, the residue $\operatorname{Res}(dt/t)=1$
identifies the logarithmic dualizing fibre with $\C$, independently
of multiplying $t$ by a unit. The two branch residues differ by a
sign, which disappears in the even power $r$. Cancelling the
base dualizing factors therefore gives the identity on the original
$P_p$.

The connected graph in Lemma~\ref{nh:fibre-residues} joins the
sources reached by the two branches. Its edge maps and the
adjunction comparisons give a crepant birational map
$\phi_{ij}:(Z_i,B_{Z_i})\dashrightarrow(Z_j,B_{Z_j})$ with
\[
 g_{ij}\circ\beta_j(p_j)=\beta_i(p_i)\circ\phi_{ij}^*,
\]
where $g_{ij}$ is the given gluing of $P$ at $p$.

\emph{Step 3.}
For a closed path $\gamma$ based on $C_i$, we compose these maps,
using the fixed pair on each component and joining the adjacent
sources at every node, including the return to the initial pair.
We obtain a B-birational self-map $\phi_\gamma$ of $(Z_i,B_{Z_i})$.
By the residue square,
\[
 \beta_i(p_i)^{-1}\operatorname{Hol}_\gamma(P)\beta_i(p_i)
       =\rho_r(\phi_\gamma),
\]
where $\rho_r$ is the representation on $H_r(Z_i,B_{Z_i})$.
Its image is finite by \cite[Theorem 1.1]{FGBrep}.
A power killing this image makes all graph holonomies trivial;
the constant component frames then glue to a nowhere-zero section
on $C$. Hence $P$ is torsion with its given node gluing, and norm
descent finishes the proof.
\end{proof}

\subsection{B-semiampleness and the base pair}

\begin{corollary}\label{cor:cbf}
Let $f:(V,\Delta)\to B$ be a projective algebraic lc-trivial fibration,
where $B$ is a normal algebraic variety and $\Delta$ is a $\Q$-divisor
on $V$. Assume that $\Delta_\eta\ge0$, where $\eta$ is the generic
point of $B$. Then its moduli b-divisor is b-semiample.
\end{corollary}

\begin{proof}
\emph{Step 1.}
We reduce to a locally stable klt-trivial fibration over a smooth
projective base.
By Chow's lemma and resolution, we replace $B$ by a smooth
quasi-projective proper birational model $B_1$ and take the induced
lc-trivial fibration $f_1:(V_1,\Delta_1)\to B_1$. We take a smooth
projective compactification $\overline B_1$ and a normal projective
compactification $\overline f_1:\overline V_1\to\overline B_1$.
We extend a canonical-bundle-formula representative by setting
\[
K_{\overline V_1}+\overline\Delta_1+\tfrac1r(\varphi)
=\overline f_1^*\overline D,
\qquad
\overline\Delta_1|_{V_1}=\Delta_1.
\]
Here $r>0$, $\varphi\in\C(V_1)^*$, and $\overline D$ extends the
corresponding divisor on $B_1$. The generic pair and the rank-one
condition are unchanged, so the extended fibration is lc-trivial
and $\overline\Delta_1$ is effective over the generic point of
$\overline B_1$. Its moduli b-divisor restricts to the original one, so semiampleness
on a model of $\overline B_1$ restricts to a model over $B$.
We may therefore assume $B$ projective.

We take a log smooth Ambro model $B'\to B$ and apply
Lemma~\ref{cbf:locally-stable-preparation} to obtain
$f_0:(Y_0,\Gamma_0)\to T_0$ with relative pluricanonical line $P_0$.
Lemma~\ref{sg:smooth-branch} gives a locally stable klt source
$g:(W,\Delta_W)\to T$ over a finite \'etale $T\to T_0$, with
the same relative pluricanonical line after pullback. For the composite
alteration $\tau:T\to B'$, it gives
\[
 \mathbf M_g=\overline{(M_g)_T},\qquad
 (M_g)_T\sim_{\Q}\tau^*M_{B'}.
\]
\emph{Step 2.}
We apply Lemma~\ref{cbf:cyclic-realization} and make the resulting
alteration. Reduced base change preserves local stability and the
pluricanonical-line comparisons of Lemma~\ref{cbf:relative-line}.
We retain the notation $g:(W,\Delta_W)\to T$ and
$\tau:T\to B'$ for the resulting family and composite alteration.
For one common sufficiently divisible $r$, the lemma and
Lemma~\ref{cbf:relative-line} give
\[
 P\simeq\overline{\mathcal H}^{\otimes r}
   \simeq\cO_T\bigl(r\tau^*M_{B'}\bigr),\qquad
 \omega_{W/T}^{[r]}(r\Delta_W)\simeq g^*P.
\]

\emph{Step 3. Numerical integrability and torsion combinatorial monodromy.}
We take a smooth projective test $(Z,E_Z)$ and a morphism
$u:Z\to\overline{T_I}$ as in Definition~\ref{prep:assumptions}, and set
\[
 Z^\circ=Z\setminus E_Z,\qquad u^\circ=u|_{Z^\circ},\qquad
 \mathbb V_Z=((u^\circ)^*P_I)^{\mathrm{tr}}.
\]
Lemmas~\ref{prep:restricted-pure-part} and~\ref{prep:schmid-comparison}
identify its Schmid highest line with $u^*\overline{\mathcal H}$.
We apply Lemma~\ref{cbf:test-abundance} to $u:Z\to T$.
On a smooth projective alteration $a:Z'\to Z$, it gives a contraction
$\rho:Z'\to Q$ and an integer $\ell>0$ with
\[
 (a^*u^*\overline{\mathcal H})^{r\ell}
 \simeq(a^*u^*P)^\ell\simeq\rho^*H,
 \qquad H\in\operatorname{Pic}(Q)\text{ nef and big}.
\]
Lemma~\ref{cbf:abundant-period} then gives
\[
 \Phi_{\mathbb V_Z}\text{ generically immersive}
 \quad\Longrightarrow\quad u^*\overline{\mathcal H}\text{ big}.
\]
Zero-dimensional tests are immediate. This establishes numerical
integrability, including all boundary tests.

For a nodal test $c:C\to T$, Lemma~\ref{cbf:relative-line}
identifies the relative line, including its node fibres,
with $c^*P\simeq(c^*\overline{\mathcal H})^r$.
If $c^*\overline{\mathcal H}$ is torsion on each normalized
component, Theorem~\ref{nh:prepared-holonomy} gives an integer $N>0$
such that
\[
 (c^*\overline{\mathcal H})^{rN}\simeq(c^*P)^N\simeq\cO_C.
\]
This establishes torsion combinatorial monodromy with the given node gluing.

By Theorem~\ref{thm:main}, $\overline{\mathcal H}$ is semiample
on $T$, and hence so is $\tau^*M_{B'}$. We descend semiampleness
to $M_{B'}$ by Lemma~\ref{sg:generically-finite-descent}. Since $\mathbf M=\overline{M_{B'}}$, the moduli b-divisor
is b-semiample.
\end{proof}

\Needspace{8\baselineskip}
The following consequence of b-semiampleness also appears in
\cite[Theorem 7.6]{BFMT}.

\begin{corollary}\label{cor:lc-image}
Let $f:(V,\Delta)\to B$ be a projective algebraic lc-trivial fibration
over a normal algebraic base. Assume that $\Delta\ge0$ is a
$\Q$-divisor on $V$ and the pair $(V,\Delta)$ is lc. Then there
is a $\Q$-divisor $\Delta_B\ge0$ on $B$ such that the pair
$(B,\Delta_B)$ is lc and
\[
K_V+\Delta\sim_{\Q}f^*(K_B+\Delta_B).
\]
If $(V,\Delta)$ is klt, $\Delta_B$ may be chosen so that
$(B,\Delta_B)$ is klt.
\end{corollary}

\begin{proof}
On a log smooth Ambro model $\mu:B'\to B$, Corollary~\ref{cor:cbf} gives
\[
M_{B'}\text{ semiample},\qquad
G\in|mM_{B'}|\text{ general},\qquad
M'=\frac1mG\sim_{\Q}M_{B'}.
\]
For sufficiently divisible $m$, Bertini's theorem gives
\[
(B',B_{B'}+M')\text{ sub-lc},
\qquad\text{sub-klt if }(V,\Delta)\text{ is klt}.
\]
We write $M'=M_{B'}+m^{-1}\operatorname{div}(u)$ for
$u\in\C(B)^*$ and set $\Delta_B=B_B+\mu_*M'\ge0$.
Then $K_B+\Delta_B=(K_B+B_B+M_B)+m^{-1}\operatorname{div}(u)$
is $\Q$-Cartier, and
\[
K_{B'}+B_{B'}+M'=\mu^*(K_B+\Delta_B).
\]
The discrepancy comparison makes $(B,\Delta_B)$ lc, respectively
klt, and the canonical bundle formula gives the asserted
$\Q$-linear equivalence.
\end{proof}

The lc Fano-type assertion below was stated by Druel
\cite[Lemma 4.5]{Druel16}. Hashizume--Hu
\cite[Theorem 1.6]{HH20} proved it when either the original variety
or its image is $\Q$-factorial. We give a proof from
Corollary~\ref{cor:lc-image}, which also yields the lc Calabi--Yau case.
A normal projective variety $X$ is \emph{of lc Calabi--Yau type}
(respectively, \emph{of lc Fano type}) if there is a $\Q$-divisor
$B\ge0$ on $X$ such that the pair $(X,B)$ is lc and
$K_X+B\sim_{\Q}0$ (respectively, $-(K_X+B)$ is ample).
For lc Fano type, allowing real boundaries gives the same class
by rational approximation \cite[proof of Theorem 1.6]{HH20}.

\begin{corollary}[Images of lc Calabi--Yau and lc Fano type varieties]\label{cor:lc-fano-image}
Let $f:X\to Y$ be a contraction of normal projective varieties.
If $X$ is of lc Calabi--Yau type (respectively, lc Fano type),
then so is $Y$.
\end{corollary}
\begin{proof}
We take a $\Q$-divisor $B\ge0$ on $X$ as above. By the projection
formula, pullback by $f$ is injective on $\Q$-Cartier $\Q$-divisors
modulo $\Q$-linear equivalence. If $K_X+B\sim_{\Q}0$,
Corollary~\ref{cor:lc-image} gives a $\Q$-divisor $B_Y\ge0$ on $Y$
such that $(Y,B_Y)$ is lc and $K_Y+B_Y\sim_{\Q}0$.

Now suppose that $-(K_X+B)$ is ample.
We take an ample $\Q$-Cartier $\Q$-divisor $H$ on $Y$
sufficiently small that $-(K_X+B)-f^*H$ is ample. For a sufficiently divisible
$m$, a general member of its $m$-fold linear system defines
\[
 A\ge0,\qquad A\sim_{\Q}-(K_X+B)-f^*H,\qquad
 (X,B+A)\text{ lc}.
\]
Hence $K_X+B+A\sim_{\Q}-f^*H$. Since $B+A\ge0$, the pair
$(X,B+A)$ is lc, and $f_*\cO_X=\cO_Y$, the rank-one condition
holds for an lc-trivial fibration. Corollary~\ref{cor:lc-image}
provides a $\Q$-divisor $B_Y\ge0$ on $Y$ such that $(Y,B_Y)$ is lc and
\[
 f^*(K_Y+B_Y)\sim_{\Q}K_X+B+A\sim_{\Q}-f^*H.
\]
The same injectivity gives $K_Y+B_Y\sim_{\Q}-H$, as required.
\end{proof}

Replacing ampleness by nefness and bigness in the definition of
lc Fano type gives \emph{lc weak Fano type}. The image statement
fails in this case.

\begin{example}[Failure for lc weak Fano type]\label{ex:weak-fano-image}
There is a birational contraction $f:X\to Y$ of normal projective
threefolds such that $X$ is smooth and of lc weak Fano type, whereas
$Y$ admits no lc boundary near one point.

Let $S$ be the smooth ruled surface over an elliptic curve used in
\cite[Example 2.1]{LX23}: $-K_S$ is nef but not big, and $S$ has
no $\R$-complement. We choose a sufficiently ample divisor $A$ giving
a projectively normal embedding of $S$. The affine cone over
$(S,A)$ is not of lc type at its vertex, by the same example.
Let $Y$ be the corresponding normal projective cone and set
\[
 p:X=\mathbb P_S(\cO_S\oplus\cO_S(A))\longrightarrow S.
\]
We use the quotient convention. Let $H$ be the tautological divisor
and $E$ the section corresponding to the quotient onto $\cO_S$.
The linear system $|H|$ defines the birational contraction
$f:X\to Y$, contracting $E$ to the vertex. Since
\[
 E\sim H-p^*A,\qquad K_X\sim-2H+p^*(K_S+A),
 \qquad -(K_X+E)\sim H+p^*(-K_S),
\]
the last divisor is nef and big: $H$ is globally generated with
$H^3=A^2>0$, and $-K_S$ is nef. The pair $(X,E)$ is log smooth
and lc, so $X$ is of lc weak Fano type. The vertex of $Y$ has an
affine neighbourhood isomorphic to the affine cone above; it has
no lc boundary, and hence $Y$ is not of lc weak Fano type.
\end{example}

\subsection{A uniform multiple for big boundaries}
\label{sec:effective-bsa}

For uniform statements we use the minimal generic index $r$:
the least positive integer with $r(K_F+\Delta_F)\sim0$ on the
generic fibre. Changing its rational trivialization changes the
moduli divisor by $r^{-1}\operatorname{div}(u)$ for a base rational
function. All multiples below are divisible by $r$.

\begin{proposition}[Uniformity for finite parameter families]
\label{prop:finite-family-uniformity}
Fix finitely many projective locally stable families
$h_i:(X_i,\Delta_i)\to S_i$ over smooth quasi-projective varieties,
whose geometric fibre pairs are klt log Calabi--Yau pairs.
Suppose that for fixed
integers $r_i>0$, each $r_i\Delta_i$ is integral and the relative
pluricanonical sheaves satisfy
\[
 \omega_{X_i/S_i}^{[r_i]}(r_i\Delta_i)\simeq h_i^*A_i,
 \qquad A_i\in\operatorname{Pic}(S_i).
\]
There is an integer $c>0$, depending only on these families, with
the following property. Let $f:(V,\Delta)\to B$ be a projective
klt-trivial fibration over a normal algebraic variety and set
$K=\C(B)$. If its induced generic pair $(V_\eta,\Delta_\eta)$
is crepant birational over $K$ to the fibre of one of these
families over a point $\operatorname{Spec}K\to S_i$, then
$c\mathbf M_f$ is b-Cartier and b-free.
\end{proposition}
\begin{proof}
\emph{Step 1.}
Before choosing $f$, we compactify the parameter families as in
Corollary~\ref{cor:cbf} and apply
Lemma~\ref{cbf:locally-stable-preparation}. Each alteration is
finite \'etale over a dense open. We spread a common resolution
and its crepant equality, shrinking so that both birational maps
and the resolution restrict to the fibres. Repeating on smooth
strata of the complement gives, by Noetherian induction, finitely
many diagrams
\[
 T_j\longrightarrow\overline U_j,\qquad
 g_j:(Y_j,\Gamma_j)\longrightarrow T_j.
\]
Here $T_j$ is smooth projective, the first map has degree $d_j$
and is finite \'etale over $U_j$, and $g_j$ is locally stable
with $K_{Y_j}+\Gamma_j\sim_{\Q,T_j}0$. Above $U_j$, its fibre
pairs are crepant birational to the corresponding parameter pairs.

We take $\rho_j>0$ with $r_i\mid\rho_j$ for the relevant $i$, so that
Lemma~\ref{cbf:relative-line} gives
\[
 \omega_{Y_j/T_j}^{[\rho_j]}(\rho_j\Gamma_j)\simeq g_j^*P_j,
 \qquad P_j\simeq\cO_{T_j}(\rho_jM_j).
\]
The lemma also gives descent of the moduli b-divisor to $T_j$.
Hence Corollary~\ref{cor:cbf} and birational descent make $P_j$
semiample on $T_j$. We take $m_j>0$ with $P_j^{m_j}$ generated.
All of $\rho_j,m_j,d_j$ now depend only on the parameter families.

\emph{Step 2.}
Suppose first that $B$ is projective. The classifying map has
generic image in some $U_j$. A residue field of its pullback
from $T_j|_{U_j}\to U_j$ gives an extension of $\C(B)$ of degree
$e\le d_j$. Resolving the induced graph and taking a common model
with a smooth Ambro model $B'\to B$, we obtain
\[
 q:W\to B',\qquad a:W\to T_j,\qquad \deg q=e,
\]
with $W$ smooth and $q$ proper. Further modifications are
birational. The pullback $a^*g_j$ is locally stable and its
$\rho_j$-pluricanonical line is $a^*P_j$ by
Lemma~\ref{cbf:relative-line}. We claim that
\begin{equation}\label{eq:uniform-actual-line}
 \cO_W(\rho_j q^*M_{B'})\simeq a^*P_j.
\end{equation}
On a common resolution, the given generic crepant comparison
makes the difference of the two log canonical divisors vertical.
It is relatively $\Q$-trivial, hence pulled back from $W$ by the
argument of Lemma~\ref{cbf:locally-stable-preparation}.
The discriminant and the total divisor in the canonical bundle
formula change by this same divisor.

Since the parameter point and comparison are defined over
$\C(B)$, the minimal generic index $r$ satisfies
$r\mid r_i\mid\rho_j$. We take compatible canonical divisors on
the common resolution and a rational frame of $P_j$ nonvanishing
at the generic image of $a$. Let $\varphi_f$ be the function in
the index-$r$ formula, and $\varphi_g$ the function in the induced
index-$\rho_j$ formula for $a^*g_j$. Their ratio has zero divisor
on the proper normal geometrically connected generic fibre, so
\[
 u=\varphi_f^{\rho_j/r}/\varphi_g\in\C(W)^*.
\]
Subtracting the formulas cancels the equal discriminant changes
and gives
\[
 \rho_j\bigl(q^*M_{B'}-M_{a^*g_j}\bigr)=\operatorname{div}(u).
\]
The line for $\rho_jM_{a^*g_j}$ is $a^*P_j$; thus this integral
principal difference proves \eqref{eq:uniform-actual-line}, also
along the boundary.

\emph{Step 3. Global generation by the norm.}
We set $D=\rho_jm_jM_{B'}$.
The divisor $q^*D$ is integral Cartier and
globally generated. It has degree zero on curves in a $q$-fibre.
Within each connected component of such a fibre, any two points
are joined by a chain of complete curves. Hence the associated
morphism is constant on each connected component.
For $x\in B'$, we choose a section nonvanishing on the whole fibre
and let $E\in|q^*D|$ be its divisor. The pushforward
$eD=q_*q^*D$ is integral, hence Cartier since $B'$ is smooth.
Writing $E=q^*D+\operatorname{div}(v)$ gives
\[
 q_*E=eD+\operatorname{div}\operatorname{Nm}_{\C(W)/\C(B')}(v).
\]
This divisor avoids $x$ by properness. Hence $e\rho_jm_jM_{B'}$
is globally generated, and we may take
\[
 c=\operatorname{lcm}_{j,\,1\le e\le d_j}(e\rho_jm_j).
\]
The moduli b-divisor descends to $B'$, proving the projective case.
The compactification in Corollary~\ref{cor:cbf} preserves both
the generic pair and its parameter point. Restriction therefore
gives the same integer over any algebraic base.
\end{proof}

\begin{corollary}[Effective b-semiampleness with big boundary]
\label{cor:effective-big-boundary}
Fix a positive integer $n$ and a finite set $I\subset(0,1)\cap\Q$. There is an
integer $c=c(n,I)>0$ such that the following holds. Let
$f:(V,\Delta)\to B$ be a projective algebraic lc-trivial fibration
of relative dimension $n$. Assume that the induced geometric generic
pair $(V_{\bar\eta},\Delta_{\bar\eta})$ is klt, that
$\Delta_{\bar\eta}\ge0$ has all its nonzero coefficients in $I$,
and that $\Delta_{\bar\eta}$ is big on $V_{\bar\eta}$.
Then $c\mathbf M_f$ is
b-Cartier and b-free.
\end{corollary}
\begin{proof}
\emph{Step 1. A crepant birational generic pair with ample boundary.}
We restrict to a nonempty smooth affine open of $B$ on which
$(V,\Delta)$ is klt and $\Delta\ge0$, and take a small
$\Q$-factorialization. Openness of klt in the coefficients of
this fixed $\Q$-Cartier boundary gives a rational $t>0$ with
$(V,(1+t)\Delta)$ klt. The number $t$ may depend on $f$.
Since
\[
 K_V+(1+t)\Delta\sim_{\Q,B}t\Delta
\]
is big over $B$, its relative log canonical model exists by
\cite[Theorem 1.2 and Corollary 1.4.3]{BCHM}.
As $K_V+\Delta\sim_{\Q,B}0$, the model is crepant for the
original pair and keeps its moduli b-divisor.
On the resulting generic pair $(F,\Delta_F)$, both
$K_F+\Delta_F\sim_{\Q}0$ and $K_F+(1+t)\Delta_F$ are $\Q$-Cartier,
and the latter is ample. Their difference $t\Delta_F$ is thus
$\Q$-Cartier and ample. Hence $-K_F\sim_{\Q}\Delta_F$ is ample
and $F$ is klt.

\emph{Step 2.}
By \cite[Theorem 1.3]{HX15}, these pairs are bounded.
Its section statement also gives a uniform integer $q>0$ with
$q(K_F+\Delta_F)\sim0$: after clearing the coefficient
denominators, a nonzero section gives an effective numerically
trivial divisor, which must be zero.
Global ACC \cite[Theorem 1.5]{HMX14}, as used in
\cite[Proposition 3.1]{HX15}, gives
$\varepsilon=\varepsilon(n,I)>0$ such that these pairs and their
underlying varieties are $\varepsilon$-klt.
On the bounded families, we apply \cite[Proposition 2.4]{HX15}
to both the pairs and their underlying varieties. Taking closures
of the required parameter sets and using Noetherian induction,
we obtain finitely many smooth strata with flat families and
$\Q$-Cartier relative canonical divisor and boundary.
We enlarge $q$ to clear the relative Cartier indices and retain
the loci where the fibres are klt, the boundary is ample with
nonzero coefficients in $I$, and $q(K_F+\Delta_F)\sim0$.
The last condition is constructible:
an invertible sheaf on an integral projective fibre is trivial
exactly when it and its inverse have nonzero sections.

Choose $c_0>0$ so that $H=c_0\Delta$ is Cartier and relatively
very ample on these families, with locally free direct image
commuting with base change. For $b\in I$, let $D_b$ be the
reduced union of boundary components of coefficient $b$.
After flattening the $D_b$, only finitely many Hilbert polynomials
occur. Frame bundles of the direct images of $\cO(H)$ parametrize
all embeddings by complete linear systems $|H|$. Together with
the $D_b$, they map to finitely many Hilbert parameter spaces
for embedded varieties and their marked subschemes.
Let $P$ be the union of these images. It is constructible by
Chevalley's theorem and includes every choice of projective frame.
For every embedded pair parametrized by $P$, we have
\[
 \Delta=\sum_{b\in I}bD_b,\qquad
 \cO_F(1)\simeq\cO_F(c_0\Delta).
\]

We restrict the universal families to locally closed strata of
$P$. Applying \cite[Proposition 2.4]{HX15} again to the pairs
and their underlying varieties, and using Noetherian induction,
we refine these into smooth strata $S_i$ with flat families
$h_i:(X_i,\Delta_i)\to S_i$ and the same $\Q$-Cartier properties.
Choose a common relative Cartier index $r_0$ divisible by $q$.
All geometric fibre pairs are klt log Calabi--Yau, so cohomology
and base change give
\[
 \omega_{X_i/S_i}^{[r_0]}(r_0\Delta_i)
 \simeq h_i^*h_{i*}\omega_{X_i/S_i}^{[r_0]}(r_0\Delta_i).
\]
The pushforward is a line bundle because its fibres have dimension
one. These families are locally stable, as their fibre pairs are klt.

\emph{Step 3. The parameter point over $\C(B)$.}
We put $K=\C(B)$.
Descending the pair to a finitely generated
characteristic-zero field, embedding that field in $\C$, and
using invariance under algebraically closed field extension give
the same bounds and Hilbert polynomials over $\overline K$.
Cartierness and very ampleness descend, so $H=c_0\Delta_F$ is
defined over $K$ and is very ample. Each $D_b$ is defined over
$K$ as well: Galois permutes geometric components only among
those with equal coefficients. A $K$-basis of
$H^0(F,\cO_F(H))$, together with these divisors, gives a point
of the Hilbert parameter space. We spread this embedded pair over
a nonempty open subset of $B$, shrinking so that its closed fibres
are klt log Calabi--Yau pairs with ample boundary and the prescribed
coefficients and polarization. All these closed points belong to
$P$. Since $P$ is constructible, the generic point belongs to $P$
as well: otherwise the inverse image of its complement would
contain a nonempty open subset of $B$ and hence a closed point.
The $K$-point therefore
lies in one of the locally closed strata $S_i$, giving
\[
 \operatorname{Spec}K\longrightarrow S_i,\qquad
 (X_i,\Delta_i)\times_{S_i}\operatorname{Spec}K\simeq(F,\Delta_F).
\]
Proposition~\ref{prop:finite-family-uniformity} now applies over
$K$ itself, with a multiple depending only on $n$ and $I$.
\end{proof}

The integer $c(n,I)$ is uniform, but the proof does not give an
explicit numerical bound: it uses generating powers of the moduli
lines on finitely many prepared parameter families.

\subsection{Extension to analytic bases}\label{sec:analytic-cbf}

We follow the Hilbert-scheme reduction in the proof of
\cite[Theorem 7.3]{BFMT}. After constructing an algebraic parameter
family, we apply Corollary~\ref{cor:cbf} and compare the moduli
lines across the boundary. Changes of projective frame induce
an algebraic action on the parameter space. We use this action to
descend the semiample contraction and its line bundle.

Analytic spaces are assumed Hausdorff and second countable.
An $f$-ample holomorphic line bundle means a line
bundle whose positive powers give relative projective embeddings
locally on the base.

\begin{lemma}\label{an:parameters}
Let $f:X\to B$ be proper and surjective between irreducible normal
complex analytic spaces, and let $L$ be an $f$-ample holomorphic line
bundle. Let $\Delta$ be a $\Q$-divisor on $X$. Suppose that, on a dense
analytic open $B^\circ\subset B$, the fibres are normal and connected,
$\Delta_b\ge0$, and
\[
 (X_b,\Delta_b)\text{ is lc},\qquad
 r(K_{X_b}+\Delta_b)\sim0
\]
for one integer $r>0$ with $r\Delta$ integral.
After a proper modification of $B$, there are a vector bundle $E$,
its projective-frame bundle $\mathcal F\to B$, and a $G$-equivariant
holomorphic map
\[
 \psi:\mathcal F\longrightarrow S^{\an},\qquad
 G=\operatorname{PGL}(\operatorname{rk}E),
\]
with Zariski-dense image, where $S$ is a normal projective variety.
A dense $G$-invariant smooth open $H\subset S$ carries a $G$-equivariant
locally stable algebraic family $h:(\mathcal Y,\Xi)\to H$ of normal
connected lc pairs, with $\Xi\ge0$ and
\[
 \omega_{\mathcal Y/H}^{[\rho]}(\rho\Xi)\simeq h^*P,
 \qquad P=h_*\omega_{\mathcal Y/H}^{[\rho]}(\rho\Xi),\qquad r\mid\rho.
\]
Here $P$ is a line bundle with its natural $G$-linearization.
On a dense open of $\mathcal F$, the pair
$\psi^*(\mathcal Y,\Xi)$ is isomorphic to
$(X,\Delta)\times_B\mathcal F$.
Over the whole frame bundle, the universal embedded family has a
meromorphic bimeromorphic correspondence with the dominating
transform of $X$.
\end{lemma}

\begin{proof}
\emph{Step 1.}
We take $k>0$ such that $L^k$ embeds the fibres over a nonempty
open of $B$, and put
\[
 \mathcal E=f_*L^k,\qquad N=\operatorname{rank}\mathcal E.
\]
We replace $B$ by the normalization of the dominating component
of $\operatorname{Gr}_N(\mathcal E)$, and $X$ by its normal
dominating transform. This is a proper modification of $B$;
the pulled-back $\mathcal E$ has a locally free rank-$N$ quotient
$E$ with torsion kernel. Evaluation therefore factors through
$f^*E\to L^k$. Blowing up its ideal gives
\[
 \widetilde X\to\mathcal Z\hookrightarrow\mathbb P_B(E).
\]
Here $\mathcal Z$ is proper over $B$ and agrees with the embedded
original family on the chosen open.

\emph{Step 2.}
We write $\sum_{j=1}^s a_j\Delta_j$ for the horizontal boundary,
with $\Delta_j$ reduced and irreducible. Each meets a fixed
compact fibre, so local finiteness of a divisor implies $s<\infty$.
Let $\mathcal D_j$ be its image in $\mathcal Z$.
By generic flatness \cite[IV.9 and IV.14]{Frisch67}, $\mathcal Z$
and the $\mathcal D_j$ have fixed Hilbert polynomials on a dense open.

In the projective flag-Hilbert bundle for these polynomials,
we consider the incidence conditions
\[
 Z'_b\subset\mathcal Z_b,\qquad
 D'_{j,b}\subset\mathcal D_{j,b}\quad(1\le j\le s).
\]
Locally on $B$, after a large twist, each condition is the
vanishing of a map of universal direct-image bundles, so the
incidence space is closed analytic and proper over $B$.
On the flat open, equality of Hilbert polynomials makes the
inclusions equalities. Thus the graph of the Hilbert map
\cite[Section 9]{Douady66} is a dense open of the dominating
component. We normalize and resolve this component.

Frames of $E$ now give a holomorphic map from $\mathcal F$ to a
fixed projective flag-Hilbert scheme. Let $S$ normalize the
Zariski closure of its image. Normality of $\mathcal F$ lifts the
map to $\psi:\mathcal F\to S^{\an}$. The construction commutes
with changes of frame, so $\psi$ is $G$-equivariant.
The resolved evaluation map and incidence construction give the
meromorphic correspondence on the whole frame bundle.

\emph{Step 3.}
On a dense open of $S$, the marked supports are relative Mumford
divisors \cite[Definition 4.68]{KolFamilies}, and the fibres are
normal and connected. The original family gives a Zariski dense
set of lc pairs. By the finite locally closed stratification in
\cite[Lemma 4.44]{KolFamilies}, the universal family is locally
stable on a dense open. Its $G$-saturation and smooth locus give $H$.

We take $\rho$ divisible by $r$ with $\rho\Xi$ integral and
$\mathcal N=\omega_{\mathcal Y/H}^{[\rho]}(\rho\Xi)$ invertible.
For an integral proper fibre,
\[
 \mathcal N_b\simeq\cO_{\mathcal Y_b}
 \quad\Longleftrightarrow\quad
 h^0(\mathcal N_b)>0\ \text{and}\ h^0(\mathcal N_b^{-1})>0.
\]
The inequalities hold on the dense set given by the original
family, hence everywhere on $H$ by semicontinuity.
Hence $\mathcal N_b$ is trivial with $h^0(\mathcal N_b)=1$.
Cohomology and base change give
\[
 P=h_*\mathcal N\in\operatorname{Pic}(H),\qquad
 h^*P\xrightarrow{\sim}\mathcal N.
\]
The action on the family induces the natural linearization of $P$,
and evaluation is equivariant.
\end{proof}

\begin{lemma}\label{an:line-comparison}
Let $g:(Y,\Gamma)\to W$ be a projective algebraic locally stable
family over a smooth base, with connected fibres and $\Gamma\ge0$.
Fix an isomorphism of invertible sheaves
\[
 \omega_{Y/W}^{[s]}(s\Gamma)\simeq g^*P,
 \qquad P\in\operatorname{Pic}(W),\quad s\Gamma\text{ integral}.
\]
Let $u:T\to W^{\an}$ be holomorphic, with $T$ smooth and connected,
and suppose that the pulled-back fibres are normal on a dense open.
Then $Y_T$ is normal and
\[
 \omega_{Y_T/T}^{[s]}(s\Gamma_T)\simeq g_T^*u^*P,
 \qquad (g_T)_*\omega_{Y_T/T}^{[s]}(s\Gamma_T)\simeq u^*P.
\]
The fibration $g_T$ has discriminant zero. On a smooth modification
$\pi:T'\to T$, its absolute crepant transform has
\[
 B_{T'}=-K_{T'/T},\qquad
 \mathcal O_{T'}(sM_{T'})\simeq\pi^*u^*P.
\]
Suppose that $f:(X,\Delta)\to T$ is another proper fibration with
$X$ normal, $f_*\cO_X=\cO_T$, $s\Delta$ integral, and
$\omega_{X/T}^{[s]}(s\Delta)\simeq f^*D$ invertible.
If a meromorphic bimeromorphic correspondence over $T$ identifies
the two pairs crepantly on a dense open, then
\[
 D\otimes\cO_T(-sB_T(f))\simeq u^*P.
\]
The comparison uses the given generic identification and commutes
with absolute crepant pullback to higher smooth models.
\end{lemma}

\begin{proof}
\emph{Step 1.}
We first verify normality of $Y_T$. It is $S_2$ by flatness, the
$S_2$ condition on the fibres and smoothness of $T$.
At horizontal codimension-one points it is regular since the
general fibre is normal. At a vertical codimension-one point,
a base uniformizer generates the maximal ideal because the
fibre is generically reduced. Hence $Y_T$ is $R_1$ as well.

Relative duality gives the comparison where the morphism is
Gorenstein and $s\Gamma$ is Cartier. This locus contains the
codimension-one points of every fibre: slc fibres are nodal
there and their boundaries avoid the nodes.
As in Lemma~\ref{cbf:relative-line}, the comparison extends to
an isomorphism of divisorial sheaves. Reduced connected fibres
give $(g_T)_*\cO_{Y_T}=\cO_T$, and the projection formula gives
the direct image.

\emph{Step 2. Vanishing of $B_T$ and the discriminant on higher models.}
At a general point $t$ of a prime divisor $\Sigma\subset T$,
we choose coordinates $(t_1,\ldots,t_e)$ with $\Sigma=(t_1=0)$
and a smooth flag
\[
 T_1\subset\cdots\subset T_e=T,\qquad
 T_j=(t_{j+1}=\cdots=t_e=0),\qquad\Sigma_j=\Sigma\cap T_j,
\]
where every $T_j$ meets the locus of normal fibres. Step~1 makes
the corresponding total spaces normal. Starting with the slc
fibre at $t$, analytic inversion of adjunction
\cite[Theorem 1.1]{Fuj24IOA}, with the different on its
normalization, makes $(Y_{T_1},\Gamma_{T_1}+g_{T_1}^*t)$ lc
near the fibre. We apply it successively to
\[
 (Y_{T_j},\Gamma_{T_j}+g_{T_j}^*\Sigma_j+Y_{T_{j-1}}),
 \qquad 2\le j\le e.
\]
Adjunction agrees with the relative sheaves in Step~1.
Removing $Y_{T_{j-1}}$ from the boundary at each stage shows
that $(Y_T,\Gamma_T+g_T^*\Sigma)$ is lc over a general point
of $\Sigma$. Hence its threshold is at least one; the reduced
fibre gives the reverse inequality. We obtain
\[
 t_\Sigma=1,\qquad\operatorname{coeff}_\Sigma B_T=0.
\]
For a smooth modification $\pi:T'\to T$, the relative boundary
$\Gamma_{T'}$ differs from the absolute crepant boundary by
$g_{T'}^*K_{T'/T}$. Applying the same calculation gives
$B_{T'}=-K_{T'/T}$. The relative pluricanonical line is
$\pi^*u^*P$, which proves the moduli-line identity on $T'$.

\emph{Step 3.}
Let $p:Z\to T$ be a common model resolving the correspondence.
The two crepant boundaries agree over a dense open, so
$E=\Delta_Z-\Gamma_Z$ is vertical, and
\[
 \cO_Z(sE)\simeq p^*\bigl(D\otimes(u^*P)^{-1}\bigr).
\]
On an open cover $\{U_i\}$ trivializing the line, the comparison
is a meromorphic function $v_i$ on $p^{-1}(U_i)$ with divisor
$sE$. Off the image of $E$, the equality $p_*\cO_Z=\cO_T$
gives $v_i=p^*w_i$. Taking the proper image of the meromorphic
graph extends $w_i$ across the remaining set.
Their transition functions are those of $D\otimes(u^*P)^{-1}$,
so their divisors define a Cartier divisor $C$ with
\[
 sE=p^*C,\qquad\cO_T(C)\simeq D\otimes(u^*P)^{-1}.
\]
Adding $p^*(C/s)$ to the boundary decreases the threshold at
$\Sigma$ by $\operatorname{coeff}_\Sigma(C)/s$. Step~2 gives
\[
 B_T(f)=B_T(g_T)+C/s=C/s,
\]
which proves the comparison. The calculation on a higher smooth
model is the same, giving compatibility with pullback.
\end{proof}

\begin{theorem}[B-semiampleness over analytic bases]
\label{thm:analytic-cbf}
Let $f:(V,\Delta)\to B$ be proper and surjective between
irreducible normal complex analytic spaces, with
$f_*\cO_V=\cO_B$. Assume that $V$ carries an $f$-ample holomorphic
line bundle and that $\Delta$ is a $\Q$-divisor. Suppose that,
for some $r>0$ with $r\Delta$ integral and some
$D\in\operatorname{Pic}(B)$, there is an isomorphism
\[
 \omega_V^{[r]}(r\Delta)\simeq f^*D
\]
of invertible sheaves. On a dense analytic open $B^\circ\subset B$,
assume that $V_b$ is normal and connected, $\Delta_b\ge0$, and
$(V_b,\Delta_b)$ is lc for every $b\in B^\circ$.
Then there exist a proper modification $\mu:B'\to B$ with $B'$
smooth, a projective variety $Z$, an ample line bundle $A$ on $Z$,
a holomorphic map $q:B'\to Z^{\an}$, and an integer $m>0$ divisible
by $r$ such that
\begin{equation}\label{eq:analytic-moduli-line}
 (\mu^*D)^{m/r}\otimes\omega_{B'}^{-m}
       \otimes\cO_{B'}(-mB_{B'})\simeq q^*A.
\end{equation}
Here $B_{B'}$ is the discriminant of the induced absolute crepant
fibration and $mB_{B'}$ is integral. The moduli $\Q$-line in
\eqref{eq:analytic-moduli-line} pulls back on every higher smooth
model. In particular, the moduli part is b-semiample.
\end{theorem}

\begin{proof}
\emph{Step 1.}
We first apply Corollary~\ref{cor:cbf} over an algebraic parameter space.
We shrink $B^\circ$ so that it is smooth and $f$ is flat there.
Adjunction gives $r(K_{V_b}+\Delta_b)\sim0$, so we apply
Lemma~\ref{an:parameters} and use its notation $S,H,P,\rho$.
We take a smooth projective modification $S_1\to S$ which is an
isomorphism over $H$, compactify the universal family, and extend a rational
trivialization of its generic log canonical power, as in Step~1
of the proof of Corollary~\ref{cor:cbf}. The resulting fibration
is lc-trivial; its vertical boundary may be negative.
By Corollary~\ref{cor:cbf}, some projective modification
$\nu:\widetilde S\to S$ carries
\[
 \widetilde c:\widetilde S\to Z,\qquad
 \cO_{\widetilde S}(\rho aM_{\widetilde S})
       \simeq\widetilde c^*A_0,\qquad a>0,
\]
where $Z$ is projective and $A_0$ is ample.
Over $H$ this line is $\nu^*P^a$ by Lemma~\ref{cbf:relative-line},
so $\widetilde c$ contracts the $\nu$-fibres. Since
$\nu_*\cO_{\widetilde S}=\cO_S$, the map and comparison descend to
\[
 c:H\to Z,\qquad P^a\simeq c^*A_0.
\]

Applying Lemma~\ref{cbf:locally-stable-preparation} over $\widetilde S$
gives a smooth projective alteration $W\to S$ of degree $d$, a
locally stable family on $W$, and $c_W:W\to Z$. Taking corresponding
multiples and powers of $a,A_0$, we fix
\[
 n=\rho a,\qquad r\mid n,\qquad
 P_W=\cO_W(nM_W)\simeq c_W^*A_0,
\]
with $n$ divisible by the prepared family's index.
We spread a common crepant resolution of the generic pairs.
Over a suitable dense open $H^\circ\subset H$, the alteration
is finite \'etale, its prepared fibres are normal, and this
resolution gives crepant birational comparisons on the fibres.
Using the same rational pluricanonical trivialization makes the
two comparisons with $A_0$ agree on $W\times_S H^\circ$.

\emph{Step 2.}
We show that $c(gh)=c(h)$ and that the comparison
$P^a\simeq c^*A_0$ respects the $G$-action.
Let $U\simeq\mathbb A^1$ be a root subgroup of $G$ and let
$h\in H^\circ$ be general. The pair over $Uh$ is constant.
Completing the orbit parameter to $\mathbb P^1$, we choose a
finite complete curve cover with a lift $\ell:C\to W$.
Its family is generically crepant birational to the product.
Lemma~\ref{cbf:relative-line} and the vertical comparison in
Lemma~\ref{cbf:locally-stable-preparation} identify their moduli
parts. Hence
\[
 \deg(c_W\ell)^*A_0=\deg\ell^*P_W=0.
\]
Ampleness makes $c$ constant on $Uh$. By density this holds for
all $h\in H$, and, since root subgroups generate $G$, we obtain
$c(gh)=c(h)$.

To prove that $P^a\simeq c^*A_0$ is equivariant, we compare the
two linearizations. Their ratio is a regular unit $\varepsilon$
on $G\times H$ with $\varepsilon(1,h)=1$.
For fixed $h$, it is constant on every translate of a root subgroup,
since units on $\mathbb A^1$ are constant. These subgroups generate
$G$, so $\varepsilon=1$.

We take connected frame charts $U_i\subset B$ with sections
$s_i:U_i\to\mathcal F$. The full torsor over $U_i$ has
Zariski-dense image in $S$ by the identity theorem. After a
constant change of frame, $\phi_i=\psi s_i$ meets $H^\circ$,
with dense inverse image. On overlaps the $\phi_i$ differ by
holomorphic $G$-valued functions, so $c\phi_i=c\phi_j$ wherever
defined. We pull back the proper graph of $S\dashrightarrow Z$,
glue the dominating components and resolve. This gives
\[
 \mu:B'\to B,\qquad q:B'\to Z^{\an}.
\]
The equivariant comparison descends on a dense open as an
identification of the $n$-th moduli line with $q^*A_0$.

\emph{Step 3.}
It remains to prove \eqref{eq:analytic-moduli-line} on all of $B'$
for one integer $m$.
We set $D'=\mu^*D$. Over $U\subset B'$ above a frame chart, we resolve
a dominating component of $U\times_S W$ to obtain
\[
 \tau:T\to U,\qquad\kappa:T\to W,\qquad c_W\kappa=q\tau,
\]
where $T$ is smooth and $\tau$ is proper and generically finite
of degree at most $d$. The incidence construction and algebraic
common model give a meromorphic correspondence on all of $T$.
The original fibration's relative log power is pulled back from
$D_T=\tau^*D'\otimes\omega_T^{-r}$. For its absolute crepant
discriminant $B_T$, Lemma~\ref{an:line-comparison} gives
\[
 D_T^{n/r}\otimes\cO_T(-nB_T)
       \simeq\kappa^*P_W\simeq(q\tau)^*A_0;
\]
in particular, $nB_T$ is integral.
If a prime $\Theta\subset T$ dominates a prime $\Sigma\subset U$, with
ramification index $e\le d$ in the finite Stein factor, then
\[
 b_\Theta=1-e+eb_\Sigma,\qquad
 ne b_\Sigma=nb_\Theta+n(e-1)\in\Z.
\]
Hence the same integers
\[
 \lambda_d=\operatorname{lcm}(1,\ldots,d),\qquad
 m=n\lambda_d,\qquad A=A_0^{\lambda_d}
\]
clear denominators on all charts. Step~2 gives a generically
defined comparison of line bundles
\[
 \sigma:\mathcal M_m\dashrightarrow q^*A,\qquad
 \mathcal M_m=(D')^{m/r}\otimes\omega_{B'}^{-m}
                   \otimes\cO_{B'}(-mB_{B'}).
\]
Upstairs it extends meromorphically as $\tau^*\sigma$.
The proper graph image therefore extends $\sigma$ meromorphically
on $U$. The ramification term $e-1$ cancels $1-e$ in the
discriminant formula, and the isomorphism upstairs gives
\[
 e\,\operatorname{ord}_\Sigma(\sigma)
       =\operatorname{ord}_\Theta(\tau^*\sigma)=0.
\]
There is no divisorial zero or pole, so normality extends $\sigma$
to an isomorphism. These extensions agree on overlaps by their
agreement on the dense open. The same argument on any smooth
modification of $B'$ gives the same $m$ and the pullback of $q^*A$.
Finally, finitely many sections generate a power of $A$; their
pullbacks generate the corresponding moduli line on all of $B'$.
\end{proof}

\appendix
\section{Sections and finite descent}\label{sec:section-descent}

We prove the two descent statements used in
Theorem~\ref{thm:frame-criterion}: lifting sections through nilpotent
thickenings and gluing them along finite diagrams. In each case,
a first-order operator records the obstruction. We construct
generating sections in its kernel, so that they lift or glue as required.

\subsection{Horizontal sections and nilpotent thickenings}

For a derivation $\delta:\cO_B\to E$, a $\delta$-connection on a
line bundle $A$ is a $\C$-linear map
\[
 \nabla:A\longrightarrow A\otimes E,
 \qquad \nabla(fs)=f\nabla(s)+s\otimes\delta(f).
\]
We write $\nabla_m$ for the induced operator on $A^m$.
In this appendix, a section of $A^m$ is called horizontal if it
lies in the kernel of $\nabla_m$.
We require the operators to be algebraic, but allow the local
horizontal frames to be analytic.

\begin{lemma}[Horizontal generating sections]\label{lem:horizontal-generation}
Let $B$ be a projective complex scheme. Let $A$ be an ample line
bundle on $B$, and let $E$ be a coherent $\cO_B$-module.
Fix a derivation $\delta:\cO_B\to E$ and a $\delta$-connection
$\nabla$ on $A$. Suppose that every point of $B^{\an}$ has an
analytic neighbourhood on which $A^{\an}$ admits a frame $e$
with $\nabla e=0$. Then there is an integer $q>0$ such that, for
every positive multiple $m$ of $q$, the sections in
$\ker H^0(\nabla_m)$ generate $A^m$.
\end{lemma}

\begin{proof}
Let $N=\operatorname{Nil}(\cO_B)$. We induct on $\dim B$ and then
on $r$ such that $N^rE=0$. In dimension zero, the local analytic
horizontal frames are algebraic and glue on the finitely many
open-and-closed components.

\emph{Step 1. The case $NE=0$: reduction.}
Put $Y=B_{\red}$ and $F=\delta(N)$. Since $NE=0$, the restriction
of $\delta$ to $N$ is $\cO_B$-linear, so $F$ is coherent and we have
\[
 \bar\delta:\cO_Y\to E/F,\qquad
 \bar\nabla:A_Y\to A_Y\otimes(E/F),\qquad A_Y=A|_Y.
\]
Let $M=\im(\Omega_Y^1\to E/F)$. We claim that $M$ vanishes at
the generic point of every component $Y_0$. For $d=\dim Y_0>0$,
we take a projective resolution $\pi:\widetilde Y\to Y_0$ and set
$\widetilde A=\pi^*(A_Y|_{Y_0})$. Suppose that $M$ has positive
generic rank on $Y_0$. On a dense smooth open $V\subset Y_0$ where
$\pi$ is an isomorphism and $M|_V$ is locally free, the kernel
\[
 \mathcal K=\ker(\Omega_V^1\longrightarrow M|_V)
 \quad\text{has rank }<d.
\]
We pull the analytic horizontal frames back to $\widetilde Y$.
Their transition units satisfy $\bar\delta(g_{ij})=0$, so
$d\log g_{ij}|_V\in\mathcal K^{\an}$. Hence
\[
 d\log g_{i_0i_1}\wedge\cdots\wedge d\log g_{i_{d-1}i_d}=0.
\]
Each product is holomorphic on the whole overlap in $\widetilde Y^{\an}$
and vanishes on its dense intersection with $\pi^{-1}V$; it therefore
vanishes identically. These \v{C}ech products represent
$c_1(\widetilde A)^d$ up to the nonzero normalization constant. Thus
$\widetilde A^d=0$, contradicting
$\widetilde A^d=(A_Y|_{Y_0})^d>0$. For $d=0$, the generic symbol is zero.

The composite of $\bar\nabla$ with
$A_Y\otimes((E/F)/M)$ has zero symbol and is $\cO_Y$-linear.
It vanishes in analytic horizontal frames, hence is zero by
faithful analytification. Thus we regard $\bar\nabla$ as taking
values in $A_Y\otimes M$. Put $I=\operatorname{Ann}(M)$.
Since $I^2M=0$, the sheaf $M$ is an $\cO_T$-module for
$T=V(I^2)\subset Y$, and the Leibniz rule gives
\[
 \bar\delta(I^2)=0,\qquad\bar\nabla(I^2A_Y)=0.
\]
It therefore induces a connection on $A_Y|_T$ with values in
$(A_Y|_T)\otimes M$. The support of $T$ has smaller dimension.
Dimension induction applies to this possibly nonreduced scheme.
For sufficiently divisible $m$, Serre vanishing gives
\[
 H^0(Y,A_Y^m)\twoheadrightarrow H^0(T,A_Y^m|_T).
\]
As $I^2M=0$, any lift of a horizontal section is horizontal on $Y$.
Sections of $I^2A_Y^m$ generate away from $T$ when $m\gg0$.
If $M=0$, all sections on $Y$ are horizontal.

\emph{Step 2. The case $NE=0$: lifting.}
Serre vanishing also lifts these sections from $Y$ to $B$.
For a lift $s$, the obstruction $\nabla_ms$ lies in $F\otimes A_Y^m$.
The restriction of $\nabla_m$ gives a coherent surjection
\[
 N\otimes A^m\longrightarrow F\otimes A_Y^m,
 \qquad n s\longmapsto s\otimes\delta(n).
\]
The formula is independent of the frame because $NE=0$.
For $m\gg0$, the map is surjective on global sections, so we can
choose $\eta\in H^0(B,NA^m)$ with
\[
 \nabla_m\eta=\nabla_ms,\qquad
 \nabla_m(s-\eta)=0,\qquad(s-\eta)|_Y=s|_Y.
\]
Nakayama's lemma shows that the corrected sections generate.
This establishes the assertion for $r=1$.

\emph{Step 3. Induction on $r$.}
For $r>1$, the preceding case applied modulo $NE$ gives sections
$s_0,\ldots,s_t$ generating $A^q$ and horizontal modulo $NE$,
for some $q>0$. The resulting map $h:B\to H\subset\mathbb P^t$
is finite because $A$ is ample. On $s_i\ne0$, the ratios $s_j/s_i$
are killed by $\bar\delta=\delta\bmod NE$. They are the pullbacks
of affine coordinates on $H$, so $h_*\bar\delta$ is
$\cO_H$-linear. Hence
\[
 C=\ker\bigl(h_*\bar\delta:h_*\cO_B\to h_*(E/NE)\bigr)
\]
is a coherent finite $\cO_H$-algebra. Exactness of analytification
and its compatibility with finite pushforward identify $C^{\an}$
with the kernel of the corresponding analytic map. Set $B_0=\Spec_H C$. We have
\[
 B\xrightarrow{g}B_0\xrightarrow{\pi}H,
 \qquad A_0=\pi^*\cO_H(1),\qquad g^*A_0\simeq A^q.
\]
Both maps are finite, and $A_0$ is ample.

We will use that $g$ is a homeomorphism analytically.
The characteristic idempotents on disjoint neighbourhoods of
points in an $h$-fibre are killed by every derivation. They belong
to $C^{\an}$ and separate the points of a $g$-fibre. Moreover,
$g$ is surjective by the integral inclusion $C\subset h_*\cO_B$
and is closed by finiteness. This establishes the assertion.

Restricting $\delta$ to $C$ and $\nabla_q$ to $A_0$ gives a
connection with values in $E_0=g_*(NE)$: its coefficients in the
frames $s_i$ lie in $NE$. The inclusion $C\subset h_*\cO_B$
sends $N_0=\operatorname{Nil}(\cO_{B_0})$ into $N$. Thus
$N_0^{r-1}(NE)\subset N^rE=0$, and hence
\[
 \operatorname{Nil}(\cO_{B_0})^{r-1}E_0=0.
\]
If $e$ is an original horizontal frame, then $e^q/s_i\in C^{\an}$.
Since $g$ is a homeomorphism, we may shrink to a saturated
neighbourhood; the ratio then descends $e^q$ to a horizontal
frame of $A_0$. Induction gives horizontal generators of $A_0^m$.
Pulling back gives horizontal generators of $A^{qm}$, and their
products give all positive multiples of this degree.
\end{proof}

\begin{lemma}[Lifting through nilpotent thickenings]\label{lem:nilpotent-frame-lifting}
Let $X$ be a proper complex algebraic space, and let $L$ be a line
bundle on $X$. Assume that $L|_{X_{\red}}$ is semiample, and let
$f_{\red}:X_{\red}\to B_{\red}$ be its connected contraction.
Suppose that there is an integer $n>0$ such that, for every fibre
$F$ of $f_{\red}$, the line $L^n$ has a nowhere-vanishing analytic
section on an open neighbourhood of $|F|$ in $X^{\an}$.
Then $L$ is semiample.
\end{lemma}

\begin{proof}
We replace $L$ by that power and induct through the powers of
the nilradical. At each step we have
\[
 0\longrightarrow J\longrightarrow\cO_X
 \longrightarrow\cO_S\longrightarrow0,
 \qquad J^2=0,
  \qquad L_S^a\simeq f^*A.
\]
Here $f:S\to B$ is the connected contraction constructed at the
preceding step, and $A$ is an ample line bundle on $B$.
The underlying fibres of $f$ are those of $f_{\red}$.
We identify the small \textup{\'etale} topoi of $S$ and $X$
\cite[Tag 05ZJ]{Stacks}, without assuming a morphism $X\to B$.
The additive connecting maps and the projection formula define
\[
 \delta:\cO_B\longrightarrow E:=R^1f_*J,
 \qquad \nabla:A\longrightarrow A\otimes E.
\]
Differences of local lifts of functions give \v{C}ech cocycles in $J$.
Since $J^2=0$, the product rule gives
$\delta(uv)=u\delta(v)+v\delta(u)$ and
$\nabla(us)=u\nabla(s)+s\otimes\delta(u)$.
For a local frame $e$ of $A$, lift $f^*e$ to frames $\ell_i$ of $L^a$.
Write $\ell_j=(1+\epsilon_{ij})\ell_i$ and
$\theta_e=[\epsilon_{ij}]\in E$. Then
\[
 \nabla_m(ue^m)=e^m\otimes\bigl(\delta(u)+mu\theta_e\bigr).
\]
Indeed, $(1+\epsilon_{ij})^m=1+m\epsilon_{ij}$, so this is exactly
the connecting obstruction for the actual line $L^{am}$.
The Leibniz rule also shows that these algebraic operators factor
through the corresponding sheaves of first principal parts.

A frame of $L^a$ near a whole fibre restricts to a frame of $f^*A$.
By properness, after shrinking an analytic neighbourhood $U\subset B$,
the frame is defined on the open thickening of $f^{-1}(U)$ in $X$.
Proper analytic comparison and $f_*\cO_S=\cO_B$ descend its restriction
and its inverse to $A|_U$ and $A^{-1}|_U$.
Their product is $1$, so the descended section is a frame.
Its analytic connecting obstruction is zero because the frame lifts to $X$.
We may therefore apply Lemma~\ref{lem:horizontal-generation} to
obtain horizontal generating sections of sufficiently divisible
powers of $A$. For $m\gg0$, Serre vanishing gives
\[
 H^1(B,f_*J\otimes A^m)=0.
\]
The Leray sequence makes the map
\[
 H^1(S,J\otimes f^*A^m)
 \longrightarrow H^0(B,R^1f_*J\otimes A^m)
\]
injective. For $s\in\ker H^0(\nabla_m)$, the obstruction to
lifting $f^*s$ maps to $\nabla_ms=0$ and hence vanishes.
Hence $f^*s$ lifts to a section of $L^{am}$ on $X$.
Nakayama's lemma shows that the lifted sections generate $L^{am}$.
They define $g:X\to\mathbb P^t$ with $g|_S=h\circ f$,
where $h:B\to\mathbb P^t$ is finite because $h^*\cO(1)=A^m$ is ample.
Since $|X|=|S|$, each $g$-fibre is the disjoint union of the old
connected $f$-fibres over the finite set $h^{-1}(p)$.
The Stein factorization $X\to B'\xrightarrow{\pi}\mathbb P^t$
therefore has exactly the old underlying fibres
\cite[Tag 0A1B]{Stacks}. Moreover, $B'$ is projective,
$\pi^*\cO(1)$ is ample, and its pullback is $L^{am}$.
This supplies the next connected contraction; the original neighbourhood
frames restrict to every thickening, so the induction continues.
\end{proof}

\subsection{Finite diagrams with compatible frames}

\begin{lemma}[Sections on a finite diagram]\label{lem:finite-polarized-diagram}
Let $B$ and $T$ be projective complex schemes, and let
$p_1,p_2:T\to B$ be finite morphisms. Let $A$ be an ample line
bundle on $B$, equipped with an isomorphism
\[
 \alpha:p_1^*A\xrightarrow{\sim}p_2^*A.
\]
Suppose that the equivalence relation generated by the pairs
$(p_1(t),p_2(t))$, for $t\in T(\C)$, is a closed reduced subspace
$R\subset B_{\red}\times B_{\red}$ whose projections are finite.
Assume also that each $R$-class has a saturated analytic open
neighbourhood $U\subset B^{\an}$ on which $A$ admits a frame $e$
satisfying
\[
 \alpha(p_1^*e)=p_2^*e
\]
on $p_1^{-1}(U)=p_2^{-1}(U)$ with its full analytic structure,
including its nilpotents.
Then sections invariant under $\alpha$ generate a positive power
of $A$. The categorical quotient is a projective scheme, and the
quotient map is finite with the $R$-classes as its fibres.
The same power of $A$ descends to an ample line bundle on the quotient.
\end{lemma}

\begin{proof}
\emph{Step 1. Reduction to invariant sections.}
We prove invariant generation by induction on $\dim B$, first for
reduced diagrams and then for their nilpotent thickenings.
Once invariant sections generate a power of $A$, they define a
finite map $h:B\to H$ with $hp_1=hp_2$. Set
\begin{equation}\label{eq:finite-invariant-algebra}
 \mathcal C=\ker\bigl(h_*\cO_B\xrightarrow{p_1^*-p_2^*}
                    (hp_1)_*\cO_T\bigr).
\end{equation}
This is a finite $\cO_H$-algebra, and $Q=\Spec_H\mathcal C$ is the
categorical quotient. Locally near a finite fibre, characteristic
idempotents constant on each class belong to the kernel and
separate the classes, even with nilpotents; compare
\cite[Lemma 17]{KolFinite}. For $\pi:Q\to H$, the ample line
$\pi^*\cO_H(1)$ pulls back to the chosen power of $A$.
Hence it suffices to construct the invariant sections.
For every finite quotient used below, its coherent equalizer
of functions remains an equalizer after analytification:
analytification is exact on coherent sheaves and commutes with
finite direct image. Thus invariant holomorphic functions,
and in particular the transition units of compatible frames,
descend analytically. We retain the descent isomorphisms with
their prescribed comparisons on the diagrams.

In dimension zero the given analytic frames are algebraic.
Assume the result in smaller dimension, and suppose first that
$B,T$ are reduced. We may replace $T$ by the generated reduced
relation: the frames extend $\alpha$ uniquely along each class,
proper GAGA algebraizes this comparison, and invariant sections
are unchanged.

\emph{Step 2. Normal $B$.}
Let $d=\dim B$, with $B$ normal, and let $B_d$ be the union of
its $d$-dimensional components. The $d$-dimensional part $R_d$
of the relation has universally open projections, so each
component of $R_d\times_{B_d}R_d$ dominates a component of $B_d$.
Transitivity extends from generic points by closure; hence $R_d$
is an equivalence relation \cite[Lemma 28]{KolFinite}. By
\cite[Lemma 21]{KolFinite}, it has a finite normal quotient.
Together with the unchanged smaller-dimensional components this
gives $q_0:B\to Q_0$ finite.

The transition units of the compatible frames are invariant.
They descend to an analytic line on $Q_0$, algebraic by proper
GAGA, such that
\[
 q_0^*A_0\simeq A.
\]
Finite descent makes $A_0$ ample and $Q_0$ projective.
The reduced image $R_0$ of the entire relation in $Q_0\times Q_0$
is a finite equivalence relation, since every $q_0$-fibre lies in
one original class. All nondiagonal components have dimension
$<d$. We take a reduced saturated closed $D\subset Q_0$ of
dimension $<d$ containing their projections and the
smaller-dimensional components. This is possible because $R_0$
is finite. The descended frames give the comparison on $R_0$;
it is algebraic by proper GAGA.

Induction gives invariant generators on $D$. For large common
powers, they lift by
\[
 H^0(Q_0,A_0^m)\twoheadrightarrow H^0(D,A_0^m|_D).
\]
Every lift is invariant: it agrees on $D$, the relation off $D$
is the diagonal, and $R_0$ is reduced. Sections vanishing on $D$
are also invariant and generate off $D$ for $m\gg0$. Pulling back
to $B$ proves the normal case.

\emph{Step 3. Reduced $B$.}
Let $\nu:\widetilde B\to B$ be the normalization, let $C\subset B$
be the full conductor, and put $\widetilde C=C\times_B\widetilde B$.
We pull the relation back to $\widetilde B\times\widetilde B$,
reduce it, and include all pairs over the same point of $B$.
Step~2 gives a finite quotient $q:\widetilde B\to Q$ and an
ample $A_Q$ whose pullback is a fixed power of $\nu^*A$.
Let $Z\subset Q$ be the full scheme-theoretic image of
$\widetilde C$. On $Z\sqcup C$ use the ample lines
$A_Q|_Z$ and the corresponding power of $A|_C$. Consider the
smaller-dimensional diagram
\[
 \widetilde C\rightrightarrows Z\sqcup C
\]
with the comparison induced by the chosen isomorphism
$q^*A_Q\simeq\nu^*A^k$, where $k$ is the chosen power. Every lift of a point of $C$ has the
same image in $Z$, since all pairs over that point occur in the
relation on $\widetilde B$. Consequently, the generated reduced
relation consists of the diagonal, the reduced image of $\widetilde C$
in $Z\times C$ and its transpose, and the reduced image of
$\widetilde C\times_Z\widetilde C$ in $C\times C$.
These images are closed and finite over both factors. Each class
therefore contains one point of $Z$ and finitely many points of
$C$. The original frames descend to $Q$ by the analytic equalizer
comparison and restrict to $C$, with agreement on the full
scheme $\widetilde C$.

Induction supplies compatible generating pairs $(s_Z,s_C)$.
For large common powers, Serre vanishing extends $s_Z$ to a
section $s_Q$ on $Q$. The pair $(q^*s_Q,s_C)$ agrees on
$\widetilde C$, and hence descends to $B$ by the conductor
identity
\[
 \cO_B=\nu_*\cO_{\widetilde B}
       \times_{\nu_*\cO_{\widetilde C}}\cO_C,
\]
tensored with the chosen power of $A$. Sections on $Q$ vanishing on $Z$, paired with zero on
$C$, generate over $Q\setminus Z$. All these sections pull back
from $Q$ on $\widetilde B$, so are invariant on the pulled-back
relation. They are invariant on $B$ as well, since the original
relation is reduced.

\emph{Step 4. Nilpotent thickenings.}
We use the simultaneous filtrations
$B_n=V(\operatorname{Nil}(\cO_B)^{n+1})$ and
$T_n=V(\operatorname{Nil}(\cO_T)^{n+1})$, which the maps $p_i$
preserve. Suppose the preceding stage $B',T'$ has a finite
quotient $q:B'\to Q$, with $Q$ projective and an equivariant
descent of the chosen power to an ample line $A_Q$.
Let $J_B,J_T$ be the ideals of $B',T'$ in the next stage.
They are square-zero ideals and are modules on $B',T'$;
no morphism from the next stage to $Q$ is assumed. Set
$r=qp_1=qp_2:T'\to Q$ and
\[
 \begin{split}
 F&=\ker\bigl(q_*J_B\xrightarrow{d}r_*J_T\bigr),\\
 E&=\coker\bigl(q_*J_B\xrightarrow{d}r_*J_T\bigr),
 \qquad d=p_1^*-p_2^*.
 \end{split}
\]
The two maps agree on $\cO_Q$ at the preceding stage, so $d$
is $\cO_Q$-linear and $F,E$ are coherent.
For a section of $A_Q^m$, $E\otimes A_Q^m$ records the local
obstruction to an invariant lift, and $F\otimes A_Q^m$ records
the difference between two such lifts.
For an affine open $V\subset Q$, the inverse images in $B',T'$
and their nilpotent thickenings are affine. If $c\in\cO_Q(V)$,
choose a lift $\widetilde c$ of $q^*c$ on the thickened inverse
image and define
\[
 \delta(c)=[p_1^*\widetilde c-p_2^*\widetilde c]\in E(V).
\]
Changing the lift changes this difference by $d(J_B)$, so
$\delta$ is well-defined. The square-zero product rule gives
$\delta(cc')=c\delta(c')+c'\delta(c)$.
Similarly, lift the pullback of a local section of $A_Q$ and
take the difference of its two pullbacks using the prescribed
comparison. Modulo $d(J_B)$ this defines
\[
 \nabla:A_Q\longrightarrow A_Q\otimes E,
 \qquad\nabla(cs)=c\nabla(s)+s\otimes\delta(c).
\]
The preceding equivariant descent identifies the target with
$A_Q\otimes E$ and makes this construction independent of the
lift. These constructions are algebraic and commute with
restriction; tensor powers give the induced operators
$\nabla_m$. A compatible frame on the thickened diagram,
restricted to the preceding stage, descends by the analytic
equalizer comparison. Its lift has zero difference, so the
descended frame is horizontal.

By Lemma~\ref{lem:horizontal-generation}, horizontal sections
generate $A_Q^m$ for sufficiently divisible $m$.
Their $E$-valued obstructions vanish, so they have local invariant
lifts. The lifts of a fixed section form a torsor under
$F\otimes A_Q^m$. For $m\gg0$, Serre vanishing gives
$H^1(Q,F\otimes A_Q^m)=0$, and the lifts glue. They generate on
$B$. Equation~\eqref{eq:finite-invariant-algebra} gives the new
quotient and its ample line. The descended frames keep their
comparison induced by $\alpha$, allowing induction through both
finite nilpotent filtrations.
\end{proof}

\begingroup
\renewcommand{\addcontentsline}[3]{}
\section*{Acknowledgements}

The first named author is supported by Overseas High-Level Young Talent
Recruitment Programs and is partially supported by the Chongqing Natural
Science Foundation Innovation and Development Joint Fund
(Grant No.\ CSTB2023NSCQ-LZX0031).
The second named author is partially supported by the National Key Research
and Development Program of China (Grant No.\ 2024YFA1013201),
the Natural Science Foundation of Chongqing (Grant No.\ CSTB2024NSCQ-LZX0040),
and the Special Project of Chongqing Municipal Science and Technology Bureau
(Grant No.\ 2025CCZ015).

\endgroup

\end{document}